\documentclass[amsart]{article}
\usepackage{amsmath,amssymb,amsthm,amsopn,amsfonts}
\usepackage{color}

\newcommand{\field}[1]{\mathbb{#1}} \newcommand{\rz}{\field{R}}
\newcommand{\R}{\field{R}}
\newcommand{\cz}{\field{C}} \newcommand{\nz}{\field{N}}

\newtheorem{theoreme}{Theorem}[section]
\newtheorem{proposition}[theoreme] {Proposition}
\newtheorem{lemme}[theoreme]{Lemma}

\newtheorem{definition}[theoreme]{Definition}
\newtheorem{remarque}[theoreme]{Remark}

\def\XXint#1#2#3{{\setbox0=\hbox{$#1{#2#3}{\int}$}
     \vcenter{\hbox{$#2#3$}}\kern-.5\wd0}}

\def\begindemonstration{\noindent{\bf Proof: \ }}
\def\enddemonstration{{\hfill $\qed$}

{\medskip}}

\numberwithin{equation}{section}

 \DeclareMathOperator\ad{ad}
 
\DeclareMathOperator\supp{supp}

\DeclareMathOperator\Ran{Ran} \DeclareMathOperator\Real{\rz e}
 
 \DeclareMathOperator\Id{Id}

  \def\12{\frac{1}{2}}

\begin{document}
\title{Adiabatic approximation for a two-level atom in a light beam}
\author{Amandine Aftalion\footnote{CNRS \& Universit{\'e} Versailles-Saint-Quentin-en-Yvelines,
Laboratoire de Math{\'e}matiques de Versailles, CNRS UMR 8100, 45
avenue des {\'E}tats-Unis, 78035 Versailles Cedex, France}\ \ and
Francis Nier\footnote{1) IRMAR, Universit{\'e} de Rennes 1, 35042 Rennes
  Cedex, France. 2) CERMICS, INRIA project-team MICMAC.}}

\date{\today}

\maketitle
\noindent\textbf{Abstract:}
Following the recent experimental realization of synthetic gauge potentials, Jean Dalibard addressed the question whether the adiabatic ansatz could be mathematically justified for a model of an atom in 2 internal
  states, shone by a quasi resonant laser beam. In this paper, we
  derive rigorously the asymptotic model guessed by the physicists,
  and show that this asymptotic analysis contains the information
  about the presence of vortices. Surprisingly, the main difficulties
  do not come from the nonlinear part but from the linear
  Hamiltonian. More precisely, the analysis of the nonlinear
  minimization problem and its asymptotic reduction to simpler ones,
  relies on an accurate partition of low and high frequencies (or momenta).
This requires to reconsider carefully previous mathematical works
about the adiabatic limit. Although the estimates are not sharp, this
asymptotic analysis provides a good insight about the validity of the
asymptotic picture, with respect to the size of the many parameters
initially put in the complete model.\\
\noindent\textbf{R{\'e}sum{\'e}~:} Suite \`a  la r{\'e}alisation exp{\'e}rimentale de champs
 de jauge artificiels, Jean Dalibard a soulev{\'e} la question de l'approximation
adiabatique pour un mod{\`e}le d'atome {\`a} deux niveaux, {\'e}clair{\'e} par un
faisceau laser r{\'e}sonnant. Dans
cet article, nous d\'erivons rigoureusement le mod{\`e}le asymptotique
devin{\'e} par les physiciens et montrons que cette analyse contient
l'information sur la pr{\'e}sence de vortex. Les difficult{\'e}s, et
c'est une surprise, ne viennent pas du terme non lin{\'e}aire. Plus
pr{\'e}cis{\'e}ment, l'analyse du probl{\`e}me non lin{\'e}aire, et la r{\'e}duction
asymptotique  {\`a} un mod{\`e}le plus simple, reposent sur une s{\'e}paration
pr{\'e}cise des grandes et basses fr{\'e}quences (ou grands et bas
moments). Cela n\'ecessite de reconsid{\'e}rer avec soin les r{\'e}sultats
math{\'e}matiques existants sur la limite adiabatique. Bien que les
estimations ne soient pas optimales, elles fournissent une bonne
intuition sur la validit{\'e} du mod{\`e}le asymptotique, par rapport aux
tailles des diff{\'e}rents param{\`e}tres  initialement mis dans le mod{\`e}le.
\section{Introduction}
 A lot of interest, both in the mathematical and physical community,
 has been devoted in the past 10 years to the study of the rotation of
 a Bose Einstein condensate: experiments \cite{MCWD1,MCWD2},
 theoretical works (see \cite{Coo,Fet} for reviews), mathematical
 contributions \cite{Aft,LiSe}.
In the first experimental production of such rotating Bose Einstein
condensates,  a rotating laser beam was superimposed on the magnetic
trap
holding the atoms in order to spin up the condensate by creating a
harmonic anisotropic rotating potential \cite{MCWD1,MCWD2}.
 Recently, a new experimental device has emerged which consists in realizing artificial or synthetic gauge magnetic forces and leads to the formation of vortex lattices
  at rest in the lab frame \cite{LCGPS}. A colloquium \cite{DGJO} has
  analyzed in detail the artificial gauge fields and their
  manifestations. In order to understand the main ingredients of the
  physics of geometrical gauge fields, \cite{DGJO} (see also
  \cite{GCYRD}) study the case of a single quantum particle
  state with a 2 levels internal structure. More complex systems with
  more than 2 internal levels are
  also discussed in \cite{DGJO} but we stick here to the simpler case of 2 levels,
  which contains all the mathematical difficulties.
A key issue is to determine whether one internal state can be followed adiabatically. A question raised by Jean Dalibard is to analyze in particular  whether vortex formation may break down the adiabatic process. In \cite{DGJO}, some conditions are provided, that we want to analyze from a mathematical point of view.

We are interested in the minimization of the  energy
\begin{eqnarray*}
\mathcal{E}_{\kappa}(\phi)&=&\int_{\rz^{2}}|\nabla \phi|^{2}+
V_{\kappa}(x,y)|\phi|^{2}+\frac{G}{2}|\phi|^{4}~dxdy\\
&&
\hspace{1cm}+
\Omega_{\kappa,\ell_{\kappa}}(x)\langle \phi\,,
\begin{pmatrix}
  \cos(\theta_{\ell_{\kappa}}(x)) & e^{i\varphi_{k}(y)}\sin(\theta_{\ell_{\kappa}}(x))\\
e^{-i\varphi_{k}(y)}\sin(\theta_{\ell_{\kappa}}(x))&-\cos(\theta_{\ell_{\kappa}}(x))
\end{pmatrix}
\phi\rangle_{\cz^{2}},
\end{eqnarray*}
where $\phi=
\begin{pmatrix}
  \phi_{1}(x,y)\\\phi_{2}(x,y)
\end{pmatrix}
\in \cz^{2}$ and
$|\phi(x,y)|^{2}=|\phi_{1}(x,y)|^{2}+|\phi_{2}(x,y)|^{2}$\,. We prescribe  that the $L^2$ norm of $\phi$ is 1.
 Here $\phi_1$ and $\phi_2$ are the internal degree of freedom of a particle: ground and excited state of the atom. It is assumed that the atom is shone by a quasi resonant laser beam.
 The functions $\Omega_{\kappa,\ell_{\kappa}}(x)$, $\varphi_{k}(y)$ and $\theta_{\ell_{\kappa}}(x)$ are
given as in \cite{DGJO} by
\begin{eqnarray*}
&&  \Omega_{\kappa,\ell_{\kappa}}(x)=\kappa\Omega(\frac{x}{\ell_{\kappa}})\hbox{ with }\Omega(x)=\sqrt{1+x^{2}}\,,
\\ && \theta_{\ell_{\kappa}}(x)=\underline{\theta}(\frac{x}{\ell_{\kappa}})\hbox{ with }
\cos(\underline{\theta}(x))=\frac{x}{\sqrt{x^{2}+1}} \quad,\quad
\sin(\underline{\theta}(x))=\frac{1}{\sqrt{x^{2}+1}}\,,\\
&& \varphi_{k}(y)=\underline{\varphi}(ky)\hbox{ with }
\underline{\varphi}(y)=y\,,
\end{eqnarray*} where $\underline\varphi $ is the phase of the propagating laser beam while $\underline\theta$ is the mixing angle.
 We define \begin{eqnarray}\label{eq.defM}
 M(x,y)=\Omega(x)
\begin{pmatrix}
  \cos(\underline{\theta}(x)) & e^{i\underline{\varphi}(y)}\sin(\underline{\theta}(x))\\
e^{-i\underline{\varphi}(y)}\sin(\underline{\theta}(x))&-\cos(\underline{\theta}(x))
\end{pmatrix}\,.
\end{eqnarray}The matrix  $M$ models the coupling between the atom and the laser.
 The $2\times 2$ matrix
$M(\frac{x}{\sqrt{\ell_{\kappa}k}},\sqrt{\ell_{\kappa}k}y)$
can be diagonalized
 in the bases $(\psi_{+},\psi_{-})$ respectively associated with the eigenvalues~$\pm \Omega(x)$,
\begin{equation}\label{psi+-}
\psi_{+}=
\begin{pmatrix}
  C\\Se^{-i\varphi}
\end{pmatrix}
\;,\;
\psi_{-}=
\begin{pmatrix}
  Se^{i\varphi}\\-C
\end{pmatrix}\,,
\end{equation}
with $C=\cos\left(\frac{1}{2}\,\underline{\theta}(\frac{x}{\sqrt{\ell_{\kappa}k}})\right)$,
$S=\sin\left(\frac{1}{2}\,\underline{\theta}(\frac{x}{\sqrt{\ell_{\kappa}k}})\right)$ and
$\varphi=\sqrt{\ell_{\kappa}k}\,y$.  When the particle follows adiabatically the eigenstate $\psi_-$,  this corresponds to set formally $\phi= u(x,y) \psi_-$, where
$ u(x,y)\in \cz$ in $\mathcal{E}_{\kappa}$. Then $u$ minimizes a Gross-Pitaevskii type energy functional
with a modified trapping potential called the geometrical gauge potential.
 The scalar
potential $V_{\kappa}(x,y)=V_{\kappa}(x,y)\Id_{\cz^{2}}$ will be adjusted in order to produce a
harmonic potential after the addition of the geometrical gauge potential from the
adiabatic theory.
 We want to justify the adiabatic approximation for states close to
$\psi_{-}$  and analyze the error term between the
initial and effective Hamiltonians.

After a rescaling, the parameter occurring in the experiments have the
following orders of magnitude
$$
\kappa\sim 10^{6}\quad,\quad G\sim 600\quad,
\quad\ell_{\kappa}\sim 25\quad,\quad k\sim 50\,,
$$
but other values can be discussed. Conditions on the strength and
spatial extent of the artificial potential  have to be prescribed in
order to induce
large circulation. Two cases,
$\ell_{\kappa} k \geq 1$ and $\ell_{\kappa} k\leq 1$, can be
distinguished and the problem has to be rewritten in two different ways
in order to apply semiclassical techniques. In fact, we will focus on
the case $\ell_{\kappa}k\geq 1$, corresponding to the previous
numerical values.
The complete analysis is carried out in the asymptotic regime  $\ell_{\kappa}k\to
+\infty$ but some partial results are also valid for
$\ell_{\kappa}k\leq 1$ or $\ell_{\kappa}k\to 0$\,.\\
A change of scale
$\phi(x,y)=\sqrt{\frac{k}{\ell_{\kappa}}}\psi(\sqrt{\frac{k}{\ell_{\kappa}}}x
, \sqrt{\frac{k}{\ell_{\kappa}}}y)$ yields a new expression for the energy:
\begin{multline*}
\int_{\rz^{2}} \frac{k}{\ell_{\kappa}}(|\partial_{x}\psi|^{2}+
|\partial_{y}\psi|^{2})+ V_{\kappa}(\sqrt{\frac{\ell_{\kappa}}{k}}x, \sqrt{\frac{\ell_{\kappa}}{k}}y)|\psi|^{2}
\\
+
\kappa\Omega(\frac{x}{\sqrt{k\ell_{\kappa}}})
\langle \psi\,,
\begin{pmatrix}
  \cos(\underline{\theta}(\frac{x}{\sqrt{k\ell_{\kappa}}})) & e^{i\underline{\varphi}(\sqrt{k\ell_{\kappa}}y)}\sin(\underline{\theta}(\frac{x}{\sqrt{k\ell_{\kappa}}}))\\
e^{-i\underline{\varphi}(\sqrt{k\ell_{\kappa}}y)}\sin(\underline{\theta}(\frac{x}{\sqrt{k\ell_{\kappa}}}))&-\cos(\underline{\theta}(\frac{x}{\sqrt{k\ell_{\kappa}}}))
\end{pmatrix}\psi\rangle_{\cz^{2}}\\
+ \frac{Gk}{2\ell_{\kappa}}|\psi|^{4}~dxdy\,.
\end{multline*}
According to the two cases $\ell_{\kappa} k\geq 1$ or $\ell_{\kappa} k\leq 1$, we
 define a small parameter $\varepsilon$ that allows to rescale the energy. In fact, we define rather
  the parameter $\varepsilon^{2+2\delta}$, where $\delta$ can be taken as a first step equal
   to $5/2$\,. The exponent $\delta >0$ is a technical trick which
provides the right quantitative estimates for the adiabatic
approximation with a quadratic kinetic energy term\,. The suitable choice of this new parameter $\delta$ is
   discussed further down in this introduction, in
   Subsections~\ref{se.mainres} and~\ref{se.scheme}.
The small parameter $\varepsilon>0$ is thus introduced
according to the two cases:
\begin{description}
\item[if $\ell_{\kappa} k \geq 1$,] then
 \begin{equation}
\label{eq.defeps1}
 \varepsilon^{2+2\delta}= \frac{k^{2}}{\kappa}\;,\;
\delta >0\quad,\quad G_{\varepsilon}=\frac{Gk}{\kappa
  \ell_{\kappa}}=\frac{G}{k\ell_{\kappa}}\varepsilon^{2+2\delta}\,;
\end{equation}
 (This leads in our example  to $\varepsilon^{2+2\delta}= 2.5~10^{-3}$.)
\item[if $\ell_{\kappa} k \leq 1$,] then
\begin{equation}
\label{eq.defeps2}
 \varepsilon^{2+2\delta}= \frac{1}{\ell_{\kappa}^{2}\kappa}\;,\;
\delta >0\quad,\quad G_{\varepsilon}=\frac{Gk}{\kappa
 \ell_{\kappa}}= Gk\ell_{\kappa}\varepsilon^{2+2\delta}.\end{equation}
\end{description} We define \begin{equation}  \label{tau}\tau=(\tau_{x},\tau_{y})=
\left\{
  \begin{array}[c]{ll}
    (\frac{1}{k\ell_{\kappa}}, 1)&\text{if}~k\ell_{\kappa} \geq 1\,,\\
   (1,k\ell_{\kappa})& \text{if}~k\ell_{\kappa} \leq 1\,.
  \end{array}
\right.
\end{equation}
In both cases, this leads to  \begin{equation}
\label{eq.rescEn1+}\kappa^{-1}{\mathcal E}_{\kappa}(\phi)=
  {\mathcal E}_{\varepsilon}(\psi)
\end{equation} where
\begin{eqnarray}
\label{eq.defEeps}
&& \mathcal{E}_{\varepsilon}(\psi)=\langle \psi\,,\, H_{Lin}\psi\rangle
+ \frac{G_{\varepsilon,\tau}}{2}\int |\psi|^{4}~dxdy\,, \label{eq.a+L2}\\
&& G_{\varepsilon,\tau}=G\tau_{x}\tau_{y}\varepsilon^{2+2\delta}\,, \label{Geps}
\end{eqnarray}and  \begin{eqnarray}
H_{Lin}&=&-\varepsilon^{2\delta}\tau_{x}\tau_{y}\varepsilon^{2}\Delta +
  V_{\varepsilon,\tau}(x,y)+
  M(\sqrt{\frac{\tau_{x}}{\tau_{y}}}x,\sqrt{\frac{\tau_{y}}{\tau_{x}}}y)
 \label{Hlin}
\,,\\ \label{eq.defV}
 V_{\varepsilon,\tau}(x,y)&=&\kappa^{-1}V_{\kappa}(\sqrt{\frac{\ell_{\kappa}}{k}}x,
 \sqrt{\frac{\ell_{\kappa}}{k}}y)\quad (\textrm{to~be~fixed})\,.
\end{eqnarray}
The quadratic energy (linear Hamiltonian) is
\begin{eqnarray}
\nonumber
\mathcal{E}_{quad,
  \varepsilon}(\psi)
&=&\int_{\rz^{2}}\varepsilon^{2+2\delta}\tau_{x}\tau_{y}|\nabla\psi|^{2}
+
V_{\varepsilon,\tau}(x,y)|\psi|^{2}
\\
\label{enerquad}
&&\hspace{1cm}+
\langle \psi\,,\, M(\sqrt{\frac{\tau_{x}}{\tau_{y}}}x,\sqrt{\frac{\tau_{y}}{\tau_{x}}}y)
\psi\rangle_{\cz^{2}} ~dxdy
\\
&=& \langle \psi\,, H_{Lin} \psi\rangle\,.\label{Hlinener}
\end{eqnarray}
  At least when $\tau_{x}=\tau_{y}=1$ and $\delta=0$, this problem looks like the standard problem of spatial adiabatic
approximation
studied in \cite{Sor}, \cite{MaSo}, \cite{MaSo2} \cite{PST}, although it requires some
adaptations because the symbols are neither bounded nor elliptic.\\

We shall consider the asymptotic analysis as $\varepsilon\to 0$ with
uniform control with respect to the parameters $G,k\ell_{\kappa}$ which
allow to fix the range of validity of the reduced models. Then we shall
consider the asymptotic behaviour of the reduced model and the whole
system as $\ell_{\kappa}k$ is large.
Specifying the right assumptions on $V_{\varepsilon,\tau}(x,y)$
or $V_{\kappa}$, possibly in a scale depending on $(\ell_{\kappa},k)$, is also an issue.

\subsection{Main result for the Gross-Pitaevskii energy}
\label{se.mainres}
We shall choose the potential $V_{\varepsilon,\tau}$ in (\ref{eq.defV}) such that after
the addition
 of the adiabatic potential in the lower energy band, the effective potential is
  almost harmonic. Namely, we assume
\begin{equation}
  \label{eq.defVepsIntro}
V_{\varepsilon,\tau} (x,y)=
\frac{\varepsilon^{2+2\delta}}{\ell_{V}^{2}}v(\sqrt{\tau_{x}}x,\sqrt{\tau_{x}}y)
+
\sqrt{1+\tau_{x}x^{2}}-\varepsilon^{2+2\delta}
\left[
\frac{\tau_{x}^{2}}{(1+\tau_{x}x^{2})^{2}}+
\frac{1}{1+\tau_{x}x^{2}}
\right]\,,
\end{equation}
with the potential $v$ chosen such that \begin{eqnarray}
\label{eq.defvIntro}
&&v(x,y)=(x^{2}+y^{2})\chi_{v}(x^{2}+y^{2})+(1-\chi_{v}(x^{2}+y^{2}))\,
\\
\nonumber
&&
\chi_{v}\in \mathcal{C}^{\infty}_{0}([0,2))\quad, \quad 0\leq
\chi_{v}\leq 1\,, \quad \chi_{v}\equiv 1~\text{on}~[0,1],
\end{eqnarray}
 and $\ell_V>0$ parametrizes the shape of the quadratic potential
 around the origin.\\
With these assumptions on the potential $V_{\varepsilon,\tau}$,
 if one chooses $\psi=u(x,y) \psi_-$, where $\psi_-$ is
 the eigenfunction corresponding to $-\Omega (x\sqrt {\tau_x/\tau_y})$
of the matrix
$\Omega (x\sqrt {\tau_x/\tau_y})M (x\sqrt {\tau_x/\tau_y}, y\sqrt{\tau_y/\tau_x})$, then the linear part
$H_{Lin}$ in ${\cal E}_\varepsilon$ is formally replaced by the scalar $\varepsilon^{2+2\delta}\tau_x\hat H_-$
where \begin{equation}
    \label{eq.defH-}
\hat{H}_{-}=
-\partial_{x}^{2}-
\left(\partial_{y}- i\frac{x}{2\sqrt{1+\tau_{x}x^{2}}}\right)^{2}
+
\frac{1}{\ell_{V}^{2}\tau_{x}}v(\sqrt{\tau_{x}}x,\sqrt{\tau_{x}}y)\,.
\end{equation} In the limit $\tau_x\to 0$,  a natural $(\varepsilon,\tau)$-independent
 scalar reduced model emerges:
$$
\mathcal{E}_{H}(u)=
\langle
u\,,\,
\left[-\partial_{x}^{2}-(\partial_{y}-\frac{ix}{2})^{2}
  + \frac{x^{2}+y^{2}}{\ell_{V}^{2}}\right]u\rangle
+\frac{G}{2}\int_{\R^2}|u|^{4}\,,\quad u(x,y)\in \cz\,.
$$
 Because $\kappa$ is large, or $\varepsilon$ is small, it is natural to expect that the ground state
 of ${\cal E}_\varepsilon$ is close, up to a unitary transform, to a vector $u(x,y) \psi_-$. This is the aim of the adiabatic theory
and leads to a scalar problem. In order to get good bounds on the energy, we need to study the limit $\tau_x$ small
at the same time.

In all our work $\ell_{V}>0$ and $G>0$ are assumed to be fixed, while
the asymptotic behaviour is studied as $\varepsilon\to 0$ and
$\tau_{x}\to 0$\,.\\
The quadratic part of the above energy is associated with the
Hamiltonian
\begin{equation}
  \label{eq.defHellV}
H_{\ell_{V}}= -\partial_{x}^{2}-(\partial_{y}-\frac{ix}{2})^{2}
  + \frac{x^{2}+y^{2}}{\ell_{V}^{2}}\,, \quad (0<\ell_{V}<+\infty)\,,
\end{equation}
with the domain
\begin{equation}
{\mathcal H}_{2}=\left\{u\in L^{2}(\rz^{2})\,,
 \sum_{|\alpha|+|\beta|\leq 2} \|q^{\alpha}D_{q}^{\beta}u\|_{L^{2}}
 <+\infty\right\}\,,\quad q=(x,y)\in \rz^{2}\,,
\end{equation}
endowed with the norm
$\|f\|_{\mathcal{H}_{2}}^{2}=\sum_{|\alpha|+|\beta|\leq 2}
\|q^{\alpha}D_{q}^{\beta}u\|_{L^{2}}^{2}$ and the corresponding
distance $d_{\mathcal{H}_{2}}$\,.
It is not difficult (see Section~\ref{se.harmapp}) to check that the
minimization of $\mathcal{E}_{H}(\varphi)$ under the constraint
$\|\varphi\|_{L^{2}}=1$, admits solutions and that the set,
$\mathrm{Argmin}~\mathcal{E}_{H}$, of ground states for $\mathcal{E}_{H}$
is a bounded set of $\mathcal{H}_{2}$\,.
\begin{definition}
  \label{de.min} For a functional $\mathcal{E}$ defined on a Hilbert
  space $\mathcal{H}$ (with $+\infty$ as a possible value), we set
$$
\mathcal{E}_{min}=\inf_{\|u\|_{L^2}=1}\mathcal{E}(u)\quad
\text{and}\quad
\mathrm{Argmin}~\mathcal{E}=\left\{u\in \mathcal{H},\;
  \mathcal{E}(u)=\mathcal{E}_{min}\,\text{and}\,\|u\|_{L^2}=1
\right\}\,.
$$
\end{definition}
\begin{theoreme}
\label{th.mintot}
Fix the constant $\ell_{V},G$ and $\delta$ and assume, for some $C_{\delta}=C(\ell_{V},G,\delta)>0$,
\begin{equation}
    \label{eq.condtaueps}
\varepsilon^{2\delta}\leq \frac{\tau_{x}^{\frac{5}{3}}}{C_{\delta}}\,.
\end{equation}Let
$\chi=(\chi_{1},\chi_{2})$ be a pair of cut-off functions such that
$\chi_{1}^{2}+\chi_{2}^{2}\equiv 1$ on $\rz^{2}$, $\chi_{1}\in
\mathcal{C}^{\infty}_{0}(\rz^{2})$ and $\chi_{1}=1$ in a neighborhood
of $0$\,.\\
There exists  $\nu_{0}=\nu_{\ell_{V},G}\in (0,\frac{1}{2}]$
and for any given $\delta>0$, there are constants
  $\tau_{\delta}=\tau(\ell_{V},G,\delta)>0$\,,
$C_{\chi,\delta}=C(\ell_{V},G,\delta,\chi)>0$\,, and a unitary operator
  $\hat{U}=\hat{U}(\varepsilon,\tau,\ell_{V},G,\delta)$ which
guarantee the  following properties
  \begin{itemize}
  \item The energy $\mathcal{E}_{\varepsilon}$ introduced in
    \eqref{eq.defEeps} admits ground states as soon as
    $\tau_{x}\leq \tau_{\delta}$, with
$$
|\mathcal{E}_{\varepsilon,min}-\varepsilon^{2+2\delta}\tau_{x}\mathcal{E}_{H,min}|\leq
C_{\delta}
\varepsilon^{2+2\delta}\tau_{x}^{\frac{5}{3}}\,.
$$
\item For any $\psi\in \mathrm{Argmin}~\mathcal{E}_{\varepsilon}$ written
  in the form $\psi=\hat{U}
  \begin{pmatrix}
    e^{i\frac{y}{2\sqrt{\tau_{x}}}}a_{+}\\
 e^{-i\frac{y}{2\sqrt{\tau_{x}}}}a_{-}
  \end{pmatrix}$, the vector $a= \begin{pmatrix}a_{+}\\ a_-\end{pmatrix}$ satisfies
\begin{eqnarray*}
&&
\|a_{+}\|_{L^{2}}\leq C_{\delta}\varepsilon^{2+2\delta}\tau_{x}
\quad,\quad
\|a\|_{H^{2}}\leq \frac{C_{\delta}}{\tau_{x}}
\quad,\quad
\|a\|_{L^{4}}+\|a\|_{L^{6}}\leq C_{\delta}\,,
\\
  && \|\chi_{2}(\tau_{x}^{\frac{1}{9}}.)a_{-}\|_{L^{2}}\leq
  C_{\chi,\delta}\tau_{x}^{\frac{1}{3}}\,,
\\
&&
d_{\mathcal{H}_{2}}(\chi_{1}(\tau_{x}^{\frac{1}{9}}.){a}_{-},
\mathrm{Argmin}~{\mathcal{E}_{H}})\leq
C_{\chi,\delta}(\tau_{x}^{\frac{2\nu_{0}}{3}}+\varepsilon)\,,\\
\\
&&
\|a_{+}\|_{L^{\infty}}\leq
C_{\delta}\varepsilon^{1+\delta}
\quad,\quad
d_{L^{\infty}}(\chi_{1}(\tau_{x}^{\frac{1}{9}}.)a_{-},
\mathrm{Argmin}~\mathcal{E}_{H})
\leq C_{\chi,\delta}(\tau_{x}^{\frac{2\nu_{0}}{3}}+\varepsilon)^{\frac{1}{2}}\,.
\end{eqnarray*}
  \end{itemize}

All the constants can be chosen uniformly with respect to $\delta\in
(0,\delta_{0}]$ for any fixed $\delta_{0}>0$\,.
\end{theoreme}
\begin{remarque} The proof is made in two steps: 1) the limit
  $\varepsilon\to 0$ corresponds to the adiabatic limit for the linear
  problem and allows to replace the linear part $H_{Lin}$, of
  $\mathcal{E}_{\varepsilon}$,  by the scalar
  $\varepsilon^{2+2\delta}\tau_{x}\hat{H}_{-}$ given by
  \eqref{eq.defH-}\,;
2) the limit $\tau_x \to 0$ allows to reduce the asymptotic
  minimization problem to a simpler one where the linear Hamiltonian
  is exactly $H_{\ell_{V}}$ given by \eqref{eq.defHellV}.\\
 One main point in the proof is to have
  precise energy estimates for the limiting problem. In our case, we obtain them in the limit
  $\tau_x\to 0$, because explicit calculations
  are more easily accessible when the magnetic field is constant, that is
  in the case of ${\mathcal E}_{H,min}$.
 In theory, if one had precise energy estimates for a general $\tau_x$
 for the intermediate adiabatic model with the linear part
 $\hat{H}_{-}$,
complete results could be
 performed for general values of~$\tau_x$\,.
\end{remarque}

\begin{remarque}
 The exponent $\nu_{0}=\nu_{\ell_{V},G}$ is a Lojasiewicz-Simon
  exponent (see Remark~\ref{re.loja} for an explanation and
  a.e. \cite{Loj,BCR} for the definition of Lojasiewicz exponents and
  \cite{Sim} for its extension to PDE problems). It is $\frac{1}{2}$ when the Lagrange multiplier
  associated with $u\in \mathrm{Argmin}~\mathcal{E}_{H}$ is a simple
  eigenvalue of $\hat{H}_{\ell_{V}}+G|\psi|^{2}$\,. There are reasons
  to think that it is the case for generic values $(\ell_{V},G)\in
  (0,+\infty)\times[0,+\infty)$, namely outside a subanalytic subset of dimension
 smaller than $1$ (and possibly $0$). Nevertheless for $G=0$, there is
 a discrete set of values of $\ell_{V}$ for which $H_{\ell_{V}}$ has multiple eigenvalues.
\end{remarque}
\noindent\textbf{Consequences and applications:}
One issue is whether  the presence of vortices (zeroes of
the wave function with circulation around them) might break the
adiabatic approach. The answer contained in our theorem is that vortices of
 $\psi$ in the original problem and vortices of the ground state of
 the Gross-Pitaevskii energy ${\mathcal E}_H$ are close and that the
  minimization of ${\mathcal E}_H$ provides all the information on the defects of
  $\psi$. Indeed, the smallness of $a_+$ in $L^\infty$ indicates that the vortices
   of $\psi$ and $a_-$ are close, and the last estimate of the theorem provides
   that the vortices of $a_-$ are close to that of the Gross-Pitaevskii problem.

 Numerically, in \cite{GCYRD},
 the authors observe vortices in a system
with artificial gauge as presented in \cite{DGJO} and modeled with
the energy $\mathcal{E}_{\varepsilon}$. They check that the
vortex pattern is close to that of the Gross-Pitaevskii energy. If one wants to use
our results, one may process in the
following way: once $G$ is fixed, choose $\ell_{V}$ such that the
minimizers of $\mathcal{E}_{H}$ have vortices (detailed conditions
 will be given in section \ref{se.harmapp}). Then
take  $\tau_{x}>0$ and $\varepsilon>0$ small enough so that
the $L^{\infty}$ norm of $a_+$ and the $L^{\infty}$
distance from $\chi(\tau_{x}^{\frac{1}{9}}.)a_{-}$ to a ground state
of $\mathcal{E}_{H}$ are small.
In the above result, the constants are not explicitly controlled and this
control is worse and worse when $\delta$ increases. So it is not
explicit for given numerical values of the parameters
$\ell_{V},G,\tau_{x},\varepsilon$ or equivalently $\ell_{V},G,
\kappa,k,\ell_{\kappa}$\,. Nevertheless it provides a
framework  for numerical
simulations, where the observation of vortices can be confirmed by
 decreasing $\varepsilon$ and $\tau_{x}$\,. The parameters of the experiments
 are just at the border where our constants may become too large
 to provide a reasonable approximation.

The definitions \eqref{eq.defeps1} of $\varepsilon$ and
\eqref{tau} of $\tau_{x}$, transform the condition
\eqref{eq.condtaueps} into
$$
\left(\frac{k^{2}}{\kappa}\right)^{\frac{2\delta}{1+2\delta}}\ll \left(
\frac{1}{k\ell_{\kappa}}\right)^{\frac{5}{3}}\,.
$$
The larger $\delta$, the better, but $\delta$ must not be too
large because of the value of the constants $C_{\delta}$,
$C_{\delta,\chi}$. A value of a few units for $\delta$
does not affect them too much. Another reason for keeping $\delta$
small is that the initial small parameter is
$\varepsilon^{2+2\delta}$, while the final error estimate of
$d_{\mathcal{H}_{2}}(\chi_{1}(\tau_{x}^{\frac{1}{9}}.)a_{-},
\textrm{Argmin}~\mathcal{E}_{H})$ (or $d_{L^{\infty}}(\chi_{1}(\tau_{x}^{\frac{1}{9}}.)a_{-},
\textrm{Argmin}~\mathcal{E}_{H})$) contains
also an $\mathcal{O}(\varepsilon)$ term.\\
As an example for $\delta=\frac{5}{2}$, the above relation becomes
$$
\frac{1}{(k\ell_{\kappa})^{7/3}}\gg \frac{k^{2}}{\kappa} \quad,\quad
\tau_{x}=\frac{1}{k\ell_{\kappa}}
\quad,\quad \varepsilon=\left(\frac{k^{2}}{\kappa}\right)^{\frac{1}{7}}\,.
$$
A given precision of order $\varepsilon+\tau_{x}^{\frac{1}{3}}$, is
more easily achieved by taking $k$ small and
$\ell_{\kappa}=\frac{1}{\tau_{x}k}$ large  (for example
$k=0.1$, $\ell_{\kappa}=500\times 25$ is better than the values $k=50$,
$\ell_{\kappa}=25$ given in the introduction). Note also that the
external potential $V_{\varepsilon,\tau}$, defined in
\eqref{eq.defVepsIntro}
must be adjusted up to the order $\varepsilon^{2+2\delta}=\frac{k^{2}}{\kappa}$\,.
\subsection{Gist of the analysis}
\label{se.scheme}

 Following the general idea of
the founding articles \cite{BoFo,BoOp} of
Born, Fock and Oppenheimer, it is well known in the physics literature that a Hamiltonian system
$$
(\varepsilon D_{q})^{2} + u_{0}(q)
\begin{pmatrix}
  E_{+}(q)&0\\
0&E_{-}(q)
\end{pmatrix}
u_{0}(q)^{*}
$$
is unitarily equivalent to a diagonal Hamiltonian  plus a remainder term:
$$
\varepsilon^{2}
\left[
\begin{pmatrix}
  (D_{q}- A(q))^{2}&0\\
0&(D_{q}+ A(q))^{2}
\end{pmatrix}
+|X(q)|^{2}\right] +R(\varepsilon)\,,
$$
where $(D_{q_{i}}\mp A_{i}(q))$ are the covariant derivatives ($\mp A(q)$ is the
adiabatic connection) associated with
the fiber bundles $u_{0}(q)
\begin{pmatrix}
  \cz\\0
\end{pmatrix}$ and $u_{0}(q)\begin{pmatrix}
  0\\\cz
\end{pmatrix}$, and $|X(q)|^{2}$  is the Born-Huang potential.\\
Since \cite{Kat} and  until
recently (\cite{NeSo,MaSo,PST,PST2}), this
problem has been widely
studied by mathematicians,
and  the remainder term
is formally:
$$
R(\varepsilon)=\sum_{i,j}\varepsilon^{2}C_{ij}(q)(\varepsilon
D_{q_{i}})(\varepsilon D_{q_{j}})+\mathcal{O}(\varepsilon^{3})\,.
$$
It is smaller than $\varepsilon^{2}$ for bounded frequencies (or
momenta) but it  has the
same size as the main term
for a typical frequency of order $\frac{1}{\varepsilon}$\,.
For a nonlinear problem or without any information about the frequency
localization of the quantum states, it is important to estimate the
error terms in the low and high frequency regimes.\\
Introducing $\delta>0$ allows to obtain at the formal level
\begin{multline*}
\varepsilon^{2+2\delta}
\left[
\begin{pmatrix}
  (D_{q}- A(q))^{2}&0\\
0&(D_{q}+ A(q))^{2}
\end{pmatrix}
+|X(q)|^{2}\right] +
\\
\sum_{i,j}\varepsilon^{2+4\delta}C_{ij}(q)(\varepsilon
D_{q_{i}})(\varepsilon D_{q_{j}})+\mathcal{O}(\varepsilon^{3+4\delta}) \,,
\end{multline*}
where the remainder term is now $\mathcal{O}(\varepsilon^{2+4\delta})$ at
the typical frequency $\frac{1}{\varepsilon}$ and thus
$\mathcal{O}(\varepsilon^{2\delta})$ times the size of the main term.
Such an error estimate can be made in the $L^{2}$-sense when applied
to some wave function $\psi$ lying in $\left\{|p|\leq
  r\right\}$, with $p$ quantized into $\varepsilon D_{q}$, or more precisely fulfilling $\psi=\chi(\varepsilon
D_{q})\psi$ for some $\chi\in \mathcal{C}^{\infty}_{0}(\left\{|p|\leq
  2r\right\})$ with
$$
\left\|R(\varepsilon) \psi\right\|\leq
C_{\chi}(\varepsilon^{2+4\delta}+
C_{N}\varepsilon^{N})\|\psi\|\leq (C'_{\chi,\delta}\varepsilon^{2\delta})\varepsilon^{2+2\delta}\|\psi\|\,.
$$
where $C_{N}$ essentially depends on the estimates of $N$-derivatives
of $u_{0}(q)$, $E_{\pm}(q)$\,. For a fixed $\delta$, we choose $N\geq
2+4\delta$\,.
\begin{remarque}
About the choice $\delta>0$,
 another strategy could be considered in order to optimize the exponent
$\delta$, w.r.t $\varepsilon$:
 under analyticity assumptions or more generally assumptions which lead
to an explicit control of $C_{N}$  in terms of $N$, one could think of
optimizing first
$C_{N}\varepsilon^{N}$ w.r.t to $N$ according to the methods of
\cite{NeSo,Sor,MaSo}. As an example with $C_{N}\leq N!$,
this would lead to
$C_{N}\varepsilon^{N}\leq Ce^{-\frac{1}{\varepsilon}}$ after choosing
$N=N(\varepsilon)=\left[\frac{1}{\varepsilon}\right]$, with some
$C\leq 1$\,.
Then taking
$\delta=\delta(\varepsilon)=-\frac{1}{4\varepsilon\log(\varepsilon)}$
would lead to
$$
CC_{\chi}(\varepsilon^{2+4\delta(\varepsilon)}+
e^{-\frac{1}{\varepsilon}})\leq 2CC_{\chi}e^{-\frac{1}{\varepsilon}}\|\psi\|\,.
$$
When $\delta(\varepsilon)\to \infty$ and as compared with
$\varepsilon^{2+2\delta(\varepsilon)}$ considered as the order $1$
term, an $\mathcal{O}(\varepsilon^{2+4\delta(\varepsilon)})$ remainder term
 is almost of order $2$\,.\\
We do not consider this optimization of $\delta$ w.r.t $\varepsilon$,
with $\lim_{\varepsilon\to 0}\delta(\varepsilon)=+\infty$, because
our initial small parameter is $\varepsilon^{2+2\delta}$ and the
 final result of Theorem~\ref{th.mintot}, (the estimate about $d_{\mathcal{H}_{2}}(\chi_{1}(\tau_{x}^{\frac{1}{9}}.)a_{-},
\textrm{Argmin}~\mathcal{E}_{H})$ in the nonlinear problem),
 contains an $\mathcal{O}(\varepsilon)$\,.  Hence we
keep $\delta>0$ independent of $\varepsilon$.
This is why $\delta$ is still present in the
constants of Theorem~\ref{th.mintot}.
\end{remarque}
We need to perform
a frequency (or momentum) truncation. We decompose a general $\psi$ into $\chi(\varepsilon
D_{q})\psi+(1-\chi(\varepsilon D_{q}))\psi$ and use rough estimates
for the part $(1-\chi(\varepsilon D_{q}))\psi$ which will be
compensated in the minimization problem for
$\mathcal{E}_{\varepsilon}$ by good a priori estimates for the norm
$\|(1-\chi(\varepsilon D_{q}))\psi\|$ of the
high-frequency part.
 We also have to check that the unitary transform
$\hat{U}$, implementing the adiabatic approximation, does not perturb too much the nonlinear part of
$\mathcal{E}_{\varepsilon}(\psi)$\,.

In our case, the limit $\varepsilon\to 0$ leads to the Born-Oppenheimer Hamiltonians
$$
-\partial_{x}^{2}-(\partial_{y}\mp\frac{ix}{2\sqrt{1+\tau_{x}x^{2}}})^{2}
+ \frac{v(\sqrt{\tau_{x}}.)}{\ell_{V}^{2}\tau_{x}}\,,
$$
which, in a second step, in the limit $\tau_{x}\to 0$,
 leads to $H_{\ell_{V}}$ (with the sign $-ix$ for the lower energy band). This means that the
convergence  to the Born-Oppenheimer Hamiltonian as
$\varepsilon\to 0$ has to be uniform w.r.t $\tau_{x}\in (0,1]$\,. This
last point requires to reconsider carefully the work of \cite{PST} by
following the uniformity w.r.t $\tau$ of the estimates given by
Weyl-H{\"o}rmander calculus (\cite{Hor,BoLe}) for $\tau$-dependent
metrics which have uniform structural constants.
\\
This is done for the low frequency part in Section~\ref{se.BO} while the basic tools of
semiclassical calculus are reviewed and adapted in the Appendix~\ref{se.semiclass}.\\
In Section~\ref{se.minipb} the error associated in the high-frequency
part is considered, as well as the effect of the unitary adiabatic
transformation on the nonlinear term.\\
Once the adiabatic approximation is well justified in this rather
involved framework, the accurate analysis, as well as the comparison when
$\tau_{x}$ is small, of the two reduced models (the one with
$H_{\ell_{V}}$ and the one with $\hat{H}_{-}$) is carried out in
Section~\ref{se.explnlmin-}. This follows the general scheme of
comparison of minimization problems: 1) Write energy estimates;
2) Use bootstrap arguments and possibly Lojasiewicz-Simon inequalities
in order to compare the minimizers in the energy space; 3) Use
the Euler-Lagrange equations in order to get a better comparison in
higher regularity spaces.\\
In Section~\ref{se.compmini}, all the information of the previous
sections is gathered in order to prove Theorem~\ref{th.mintot} :
 existence of a minimizer for ${\mathcal E}_\varepsilon$ in Proposition \ref{pr.existEeps}, key energy estimates in Proposition~\ref{pr.compen}
  and bounds for minimizers in Proposition~\ref{pr.adEL}.\\
Some comments and additional results are pointed out in
Section~\ref{se.complements}, namely: 1) the question of the smallness
condition of $\varepsilon$ w.r.t $\tau_{x}$ appearing in
Theorem~\ref{th.mintot}; 2) the possible extension to anisotropic
nonlinearities (one would have to check that the unitary transform
implementing the adiabatic approximation does not perturb the nonlinear part); 3) the minimization problem for excited states,
i.e. locally and approximately carried by $\psi_{+}$ instead of $\psi_{-}$; 4) the
extension to the problem of the time nonlinear dynamics of
adiabatically prepared states.


\section{Adiabatic approximation for the linear problem}
\label{se.BO}

In \cite{PST}, the adiabatic approximation is completely justified for
bounded symbols or when global elliptic properties of the complete
matricial symbol allow to reduce to this case after spectral
truncation. Unfortunately, it is not
the case here, because the eigenvalues of the symbol of the linear part
$H_{Lin}$ are $\varepsilon^{2\delta}\tau_{x}\tau_{y}|p|^{2}+ V_{\varepsilon,\tau}(x,y)\pm
\Omega(\sqrt{\frac{\tau_{x}}{\tau_{y}}}x)$\,. In \cite{Sor}, the adiabatic theory for unbounded symbols
is developed after stopping the complete asymptotic expansion in an
optimal way, under some analyticity assumptions, but this would be
particularly tricky here with divergences occurring both in the momentum
and position directions. We shall see that the sublinear divergence in
position makes no difficulty after using the right Weyl-H{\"o}rmander
class. The quadratic divergence in momentum, with the kinetic energy
$\tau_{x}\tau_{y}|p|^{2}$ is solved by first considering truncated kinetic energies and
using the additional scaling  factor $\varepsilon^{2\delta}$,
$\delta>0$, in front of the kinetic energy term.\\
Our problem shows an anisotropy in the position variables $(x,y)$\,. The analysis of the
linear problem can be treated in $\rz^{d}$\,. Then we split the
position and momentum variables, $q\in\rz^{d}$ and $p\in\rz^{d}$,
into:
$$
q=(q',q'')\quad,\quad p=(p',p'')\quad q',p'\in \rz^{d'}\,,\,
q'',p''\in\rz^{d''}\quad,\quad d'+d''=d\,,
$$
and the pair $(\tau_{x},\tau_{y})$ is accordingly denoted by
$(\tau',\tau'')\in (0,1]^{2}$\,.\\
From this section, some notions and notations related with
semiclassical analysis are used. In particular, the notation
$S_{u}(m_{\tau},g_{\tau})$ refers to classes of
$(\varepsilon,\tau)$-dependent
symbols of which the seminorms are \underline{u}niformly controlled
w.r.t to the parameters $(\varepsilon,\tau)\in
(0,\varepsilon_{0})\times(0,1]^{2}$\,. For accurate definitions, we refer the reader to
appendix~\ref{se.semiclass} where all the necessary material is
reviewed and adapted for our analysis, assuming knowledges about
Fr{\'e}chet spaces and generalized functions.

\subsection{Born-Oppenheimer Hamiltonian}
Consider the Hamiltonian in $\hat{H}_{\varepsilon}=H(q,\varepsilon
D_{q},\varepsilon)$
 with the symbol on $\rz^{2d}_{q,p}$
 \begin{eqnarray*}
&&H(q,p,\varepsilon)=\varepsilon^{2\delta}\tau'\tau''|p|^{2}\gamma(\tau'\tau''|p|^{2})+\mathcal{V}(q,\tau,\varepsilon)
\\
&&\;=
\varepsilon^{2\delta}\tau'\tau''|p|^{2}\gamma(\tau'\tau''|p|^{2})
+u_{0}(q,\tau,\varepsilon)
\begin{pmatrix}
  E_{+}(q,\tau,\varepsilon)& 0\\
0 & E_{-}(q,\tau,\varepsilon)
\end{pmatrix}
u_{0}^{*}(q,\tau,\varepsilon)\,,
\end{eqnarray*}
with $\delta>0$ and $(\tau',\tau'')\in (0,1]^{2}$\,.
The following properties, with the splitting of variables
$q=(q',q'')\in \rz^{d}$, are assumed:
\begin{eqnarray}
\nonumber
&&
  E_{\pm}\in S_{u}\Big(\langle \sqrt{\frac{\tau'}{\tau''}}q'\rangle,
  \frac{\frac{\tau'}{\tau''}d{q'}^{2}}{\langle
    \sqrt{\frac{\tau'}{\tau''}}q'\rangle^{2}}+
\frac{\tau''}{\tau'}d{q''}^{2}\Big)
\;,\\
&&
\label{eq.hypellEpm}
E_{+}(q,\tau,\varepsilon)-E_{-}(q,\tau\varepsilon)\geq C^{-1}\langle
\sqrt{\frac{\tau'}{\tau''}}q'\rangle\,,\\
&&
\nonumber
u_{0}=(u_{0}^{*})^{-1}\in S_{u}(1,
 \frac{\frac{\tau'}{\tau''}d{q'}^{2}}{\langle
    \sqrt{\frac{\tau'}{\tau''}}q'\rangle^{2}}+
\frac{\tau''}{\tau'}d{q''}^{2}; \mathcal{M}_{2}(\cz))\,,
\\
&&
\nonumber
\quad \gamma\in \mathcal{C}^{\infty}_{0}(\rz;\rz_{+})\quad\text{with}\quad
\gamma\equiv 1 \in [0,2r_{\gamma}^{2}]\,.
\end{eqnarray}
With these assumptions we are able to justify the Born-Oppenheimer
adiabatic approximation for $\delta>0$, $\tau',\tau''\in (0,1]$.
We shall work with the $\tau$-dependent metric
$$
g_{\tau}=
 \frac{\frac{\tau'}{\tau''}d{q'}^{2}}{\langle
    \sqrt{\frac{\tau'}{\tau''}}q'\rangle^{2}}+
\frac{\tau''}{\tau'}d{q''}^{2}
+\frac{\tau'\tau''dp^{2}}{\langle
 \sqrt{\tau'\tau''} p\rangle^{2}}\quad,
$$
on the phase-space $\rz^{2d}_{q,p}$, which is checked to have uniform
properties w.r.t $\tau\in (0,1]^{2}$ in Proposition~\ref{pr.metadm}\,.
The exact definition of the
parameter dependent H{\"o}rmander symbol classes,
$S_{u}(m,g_{\tau};\mathcal{M}_{2}(\cz))$, is given in
Appendix~\ref{se.semiclass}
(see in particular its meaning for $\tau$-dependent metrics in
Appendix~\ref{se.variametric}). We shall use the notation
$$
S_{u}\Big(\frac{1}{\langle \sqrt{\frac{\tau'}{\tau''}}q'\rangle \langle
  \sqrt{\tau'\tau''} p
\rangle^{\infty}},
g_{\tau};\mathcal{M}_{2}(\cz)\Big)
$$ with the meaning of the exponent $\infty$ being
the same as in
$\mathcal{C}^{\infty}$ of $\mathcal{O}(\varepsilon^{\infty})$\,.
\begin{theoreme}
\label{th.BornOp}
There exists a unitary operator
$\hat{U}=U(q,\varepsilon D_{q},\tau, \varepsilon)$ with
symbol
\begin{eqnarray*}
&&U(q,p,\tau,\varepsilon)=u_{0}(q,\tau,\varepsilon) + \varepsilon
u_{1}(q,p,\tau,\varepsilon)+ \varepsilon^{2}u_{2}(q,p,\tau,\varepsilon)
\,,\\
&&\varepsilon^{-2\delta}u_{1,2}\in
S_{u}\Big(\frac{1}{\langle \sqrt{\frac{\tau'}{\tau''}}q'\rangle \langle
  \sqrt{\tau'\tau''} p
\rangle^{\infty}},
g_{\tau};\mathcal{M}_{2}(\cz)\Big)\,,
\end{eqnarray*}
such that
\begin{multline}
  \label{eq.BornOpp1}
\hat{U}^{*}\hat{H}\hat{U}
=
\begin{pmatrix}
  h_{BO,+}(q,\varepsilon D_{q},\varepsilon)&0\\
0& h_{BO,-}(q,\varepsilon D_{q},\varepsilon)
\end{pmatrix}\\
+ \varepsilon^{2+4\delta}R_{1}(q,\varepsilon
D_{q},\tau,\varepsilon)+ \varepsilon^{3+2\delta}R_{2}(q,\varepsilon
D_{q},\tau,\varepsilon)\,,
\end{multline}
with $h_{BO\pm}(q,p,\tau,\varepsilon)$ equal to
\begin{eqnarray}
\nonumber
&=&\varepsilon^{2\delta}\tau'\tau''(|p
\mp \varepsilon
A|^{2}+|\varepsilon X|^{2}) + E_{\pm}
\\
\label{eq.hBOpm}
&=&\varepsilon^{2\delta}\tau'\tau''\left[\sum_{k=1}^{d}(p_{k}\mp \varepsilon A_{k})^{2}
+|\varepsilon X_{k}|^{2}\right]
+ E_{\pm}\,,
\quad
\end{eqnarray}
when $\sqrt{\tau'\tau''}|p|\leq r_{\gamma}$, and
$$
\begin{pmatrix}
  +A_{k}&X_{k}\\
\bar{X_{k}}&-A_{k}
\end{pmatrix}=iu_{0}^{*}(\partial_{q^{k}}u_{0})\,.
$$
The remainder terms satisfy $R_{1},R_{2}\in
S_{u}(1,g_{\tau};\mathcal{M}_{2}(\cz))$ and $R_{2}$ vanishes in
$\left\{\sqrt{\tau'\tau''}|p|\leq r_{\gamma}\right\}$ and those
estimates are uniform w.r.t  $\delta\in
(0,\delta_{0}]$ (and $\tau\in (0,1]^{2}$, $\varepsilon\in (0,\varepsilon_{0}]$).
\end{theoreme}
\begin{remarque}
  \begin{itemize}
  \item
The remainder term is really negligible only for
  $\delta>0$. This is explained in Subsection~\ref{se.discuss}.
\item  The (Weyl)-quantization of $\varepsilon^{2\delta}\tau'\tau''(|p\mp \varepsilon
  A|^{2}+|\varepsilon X|^{2})+E_{\pm}$ is nothing but
$$
\varepsilon^{2+2\delta}\tau'\tau''\left[\sum_{k=1}^{d}-(\partial_{q^{k}}\mp
  iA_{k})^{2}+|X_{k}|^{2}\right]
+ E_{\pm}\,.
$$
\end{itemize}
\end{remarque}

\subsection{Second order computations for space adiabatic approximate projections of
  the reduced Hamiltonian}
\label{se.secondord}

We shall consider the matricial symbol, on $\rz^{2d}_{q,p}$,
$$
H(q,p,\tau,\varepsilon)=
f_{\varepsilon}(p,\tau)
+ E_{+}(q,\tau,\varepsilon)\Pi_{0}(q,\tau,\varepsilon)
+ E_{-}(q,\tau,\varepsilon)(1- \Pi_{0}(q,\tau,\varepsilon))
$$
where $\Pi_{0}(q,\tau,\varepsilon)=\Pi_{0}(q,\tau,\varepsilon)^{*}=\Pi_{0}(q,\tau,\varepsilon)^{2}\in {\mathcal
  M}_{2}(\cz)$\,,
$\gamma,E_{\pm}$ real-valued,
with $\delta \geq 0$ and the following properties:
\begin{eqnarray}
&&
\label{eq.hypEPi1}
E_{\pm}\in S_{u}(\langle \sqrt{\frac{\tau'}{\tau''}}q'\rangle,
  g_{q\tau})
\quad,\quad
\Pi_{0}\in S_{u}(1, g_{q,\tau}; \mathcal{M}_{2}(\cz))\,,
\\
\nonumber
\text{with}&&
g_{q,\tau}=\frac{\frac{\tau'}{\tau''}d{q'}^{2}}{\langle
    \sqrt{\frac{\tau'}{\tau''}}q'\rangle^{2}}+
\frac{\tau''}{\tau'}d{q''}^{2},
\\
&&
\label{eq.hypEPi2}
E_{+}(q,\tau,\varepsilon)-E_{-}(q,\tau,\varepsilon)\geq C^{-1}\langle
\sqrt{\frac{\tau'}{\tau''}}
q'\rangle\,,\\
&&
\nonumber
f_{\varepsilon}(p,\tau)=\varepsilon^{2\delta}f_{1}(\sqrt{\tau'\tau''}p)\quad,
\quad f_{1}\in \mathcal{C}^{\infty}_{0}(\rz^{d};\rz)
\quad,\quad \delta \geq 0\,.
\end{eqnarray}
For conciseness, the arguments $p$, $q$ and the parameters
$\tau,\varepsilon$,
 will often be omitted in
$\partial_{p}^{\alpha}f_{\varepsilon}(p)$ and
$\partial^{\alpha}_{q}E_{\pm}(q)$ or $\partial^{\alpha}_{q}\Pi_{0}$\,.\\
According to Appendix~\ref{se.variametric}, the metric
$$
g_{\tau}=g_{q,\tau}
+ \frac{\tau'\tau''dp^{2}}{\langle
     \sqrt{\tau'\tau''}p\rangle^{2}}
=
\frac{\frac{\tau'}{\tau''}d{q'}^{2}}{\langle
   \sqrt{\frac{\tau'}{\tau''}}q'\rangle^{2}}
+
\frac{\tau''}{\tau'}d{q''}^{2}
+
\frac{\tau'\tau''dp^{2}}{\langle
     \sqrt{\tau'\tau''}p\rangle^{2}}
$$
has the gain function
$$
\lambda(q,p)=\min\left(\frac{\langle \sqrt{\tau'\tau''}p\rangle\langle
    \sqrt{\frac{\tau'}{\tau''}}q\rangle}{\tau'},\frac{\langle
    \sqrt{\tau'\tau''}p\rangle}
{\tau''}\right)\geq \langle \sqrt{\tau'\tau''}p\rangle\,.
$$
The  $\varepsilon$-quantized version of the symbol $A_{\varepsilon}(q,p)$ will be denoted
$$
\hat A= A(q,\varepsilon D_{q})\,.
$$
Note that the symbols $E_{+}-E_{-}$ is elliptic in
 its class $S_{u}(\langle \sqrt{\frac{\tau'}{\tau''}}q'
\rangle, g_{\tau})$\,,
and
$$
\varepsilon^{-2\delta}f_{\varepsilon}\in
S_{u}(\frac{1}{\langle \sqrt{\tau'\tau''}p\rangle^{\infty}}, g_{\tau})\,.
$$
Our aim is to compute accurately the adiabatic projection
$\Pi^{(n)}(q,p)=\Pi_{0}(q,\varepsilon)+\varepsilon \Pi_{1}(q,p,\varepsilon)+\cdots +
\varepsilon^{n}\Pi_{n}(q,p,\varepsilon)$ such that
$$
\hat \Pi^{(n)}\circ \hat\Pi^{(n)}= \hat \Pi^{(n)}+{\mathcal O}(\varepsilon^{n+1})
\quad,\quad
\left[
\hat H, \hat \Pi^{(n)}\right]= {\mathcal O}(\varepsilon^{n+1})\,.
$$
The general theory presented in \cite{PST}, tells us
that the asymptotic expansion can be
pushed up to $n=\infty$, but we will do here
 accurate calculations up to
$n=2$ (with additional information for $n=3$) and then discuss the influence of the factor
$\varepsilon^{2\delta}$\,.
 Those are feasible and rather easy because the kinetic
energy term and the two-level potential are simple.
This allows to
reconsider accurately the arguments sketched in \cite{PST} for the
Born-Oppenheimer
case with all the technical new peculiarities of our example. Note
also that in \cite{MaSo2}~chapter~10 explicit calculations have been made up to
order $n=5$ but in the slightly different framework of time-dependent
Born-Oppenheimer approximation oriented to polyatomic molecules~:~no use of the
exponent $\delta>0$, no divergence as $q\to \infty$ and no ellipticity
problem, no
extra-parameter $\tau$ and the techniques
are slightly different although still relying on semiclassical calculus.\\
Like in \cite{PST}, \cite{Sor}, \cite{MaSo}, \cite{MaSo2}, the calculations are first
done at the symbolic level and we write
$$
\left(A(\varepsilon)=\mathcal{O}_{S}(\varepsilon^{\nu})\right)
\Leftrightarrow
\left(
\varepsilon^{-\nu}A\in
S_{u}(1,g_{\tau};\mathcal{M}_{2}(\cz))
\right)\,.
$$
For a matricial symbol $A(q,p)$ (possibly depending on $(\varepsilon,\tau)$), it is convenient to introduce the
diagonal and off-diagonal parts
\begin{eqnarray}
  \label{eq.diag}
A^{D}(q,p)&:=&\Pi_{0}(q)A(q,p)\Pi_{0}(q)+ (1-\Pi_{0}(q))A(q,p)(1-\Pi_{0}(q))
\\
\label{eq.offdiag}
A^{OD}(q,p)&:=&\Pi_{0}(q)A(q,p)(1-\Pi_{0}(q))+(1-\Pi_{0}(q))
A(q,p)\Pi_{0}(q)\,,
\end{eqnarray}
where we recall $\Pi_{0}(q)=\Pi_{0}(q,\tau,\varepsilon)$\,.
Note the equalities
\begin{equation}
  \label{eq.diagprod}
(AB)^{D}=A^{D}B^{D}+A^{OD}B^{OD}\quad, \quad (AB)^{OD}=A^{D}B^{OD}+A^{OD}B^{D}\,.
\end{equation}
The ``Pauli matrix''
\begin{equation}
  \label{eq.sigz}
  \sigma_{3}(q,\tau,\varepsilon)=2\Pi_{0}(q,\tau,\varepsilon)-\Id_{\cz^{2}}\,,
\end{equation}
will also be used,
with the relations
\begin{eqnarray*}
&&\sigma_{3}^{2}(q)=\Id_{\cz^{2}}\,,\quad
\sigma_{3}(q) A^{D}(q,p)\sigma_{3}(q)= A^{D}(q,p)\,,\\
&&\sigma_{3}(q)A^{OD}(q,p)\sigma_{3}(q)=-A^{OD}(q,p)\,,\\
&&
\sigma_{3}(q)A^{OD}(q,p)=\Pi_{0}(q)A(q,p)(1-\Pi_{0}(q))-(1-\Pi_{0}(q))A(q,p)\Pi_{0}(q)\,.
\end{eqnarray*}
We are looking for
\begin{eqnarray*}
&&\Pi^{(n)}(q,p,\tau,\varepsilon)
=
\sum_{j=0}^{n}\varepsilon^{j}\Pi_{j}(q,p,\tau,\varepsilon)\quad\in S_{u}(1,g_{\tau};\mathcal{M}_{2}(\cz))\,,
\\
\text{with}
&&
\Pi_{j}\in
S_{u}(1,g_{\tau};\mathcal{M}_{2}(\cz))\,,
\end{eqnarray*}
and such that
\begin{eqnarray}
  \label{eq.proj}
&& \Pi^{(n)}\sharp^{\varepsilon}\Pi^{(n)}-\Pi^{(n)}=
  \mathcal{O}_{S}(\varepsilon^{n+1})\,,\\
\label{eq.adj}
&&
{\Pi^{(n)}}^{*}=\Pi^{(n)}\,,
\\
\label{eq.comm}
&&
H\sharp^{\varepsilon}\Pi^{(n)}-\Pi^{(n)}\sharp H=\mathcal{O}_{S}(\varepsilon^{n+1})
\,.\end{eqnarray}
Like in \cite{MaSo}, \cite{PST}, this system is solved by induction by
starting from
$\Pi^{(0)}(q,p,\tau,\varepsilon)=\Pi_{0}(q,\tau,\varepsilon)$, with
\begin{eqnarray}
  \label{eq.defGn1}
  G_{n+1}&:=&
\varepsilon^{-(n+1)}\left[\Pi^{(n)}\sharp^{\varepsilon}\Pi^{(n)}-\Pi^{(n)}\right]
\mod
\mathcal{O}_{S}(\varepsilon)\,,\\
\label{eq.Pin1D}
\Pi_{n+1}^{D}&:=&-\sigma_{3}G_{n+1}^{D}\,,\\
\nonumber
F_{n+1}&:=&\varepsilon^{-(n+1)}
\left[H\sharp^{\varepsilon}(\Pi^{(n)}+\varepsilon^{n+1}\Pi_{n+1}^{D})\right.\\
\nonumber
&&
\hspace{3cm}
\left.-(\Pi^{(n)}+\varepsilon^{n+1}\Pi_{n+1}^{D})\sharp^{\varepsilon}H\right]\mod\mathcal{O}_{S}(\varepsilon)
\\
\label{eq.defFn1}
&=&
\varepsilon^{-(n+1)}
\left[H\sharp^{\varepsilon}\Pi^{(n)}-\Pi^{(n)}\sharp^{\varepsilon}H\right]
\mod\mathcal{O}_{S}(\varepsilon)\,,
\\
\label{eq.Pin1OD}
\Pi_{n+1}^{OD}&:=&
-\frac{1}{E_{+}(q)-E_{-}(q)} \sigma_{3}F_{n+1}^{OD}=\frac{1}{E_{+}(q)-E_{-}(q)} F_{n+1}^{OD}\sigma_{3}\,.
\end{eqnarray}

The general theory says that the principal symbol of
$F_{n+1}$ is off-diagonal,
$F_{n+1}=F_{n+1}^{OD}\mod\mathcal{O}_{S}(\varepsilon)$\,,
 and $F_{n+1}$ can be chosen so that
 \begin{equation}
   \label{eq.FOD}
F_{n+1}=F_{n+1}^{OD}\,.
\end{equation}
Below are the computations up to $n=2$ in our specific case. In these
computations, we shall use Einstein's summation rule
$s_{k}t^{k}=\sum_{k}s_{k}t^{k}$ with
the coordinates $(p_{k},q^{k})$ or $(p^{k},q^{k})$ with
$p^{k}=p^{\ell}\delta_{\ell,k}=p_{k}$ like in the examples
$$
|p|^{2}=p^{k}p_{k}=p_{k}p_{\ell}\delta^{\ell,k}=p^{k}p^{\ell}\delta_{k,\ell}\quad,\quad
(\partial_{p}f_{\varepsilon}).\partial_{q}=(\partial_{p^{k}}f_{\varepsilon})\delta^{k,\ell}\frac{\partial}{\partial_{q^{\ell}}}=(\partial_{p_{k}}f_{\varepsilon})\partial_{q^{k}}\,.
$$
\noindent{$\mathbf{n=0}$:} Start with $\Pi^{(0)}=\Pi_{0}(q,\tau,\varepsilon)$ and notice
$\partial_{p}\Pi_{0}\equiv 0$ and $\partial_{q}\Pi_{0}\equiv (\partial_{q}\Pi_{0})^{OD}$\,.\\
\noindent{$\mathbf{n=1}$:} Take
\begin{eqnarray*}
&&G_{1}=\Pi_{0}\circ \Pi_{0}-\Pi_{0}=0\,,
\\
\text{and}&&
\Pi_{1}^{D}=0\,.
\end{eqnarray*}
Next compute
\begin{eqnarray*}
\varepsilon^{-1}\left[H\sharp^{\varepsilon}\Pi_{0}(q,\varepsilon)-\Pi_{0}(q,\varepsilon)\sharp^{\varepsilon}H\right]
&=&\varepsilon^{-1}\left[f_{\varepsilon}\sharp^{\varepsilon}\Pi_{0}(q,\varepsilon)-\Pi_{0}(q,\varepsilon)
\sharp^{\varepsilon}f_{\varepsilon}\right]\\
&=& -i\partial_{p_{k}}f_{\varepsilon}\partial_{q^{k}}\Pi_{0}
\mod \mathcal{O}_{S}(\varepsilon)\,,
\end{eqnarray*}
and take
\begin{eqnarray*}
&&
F_{1}=-i \partial_{p_{k}}f_{\varepsilon}\partial_{q^{k}}\Pi_{0}=
-i\partial_{p_{k}}f_{\varepsilon}(\partial_{q^{k}}\Pi_{0})^{OD}\,,
\\
\text{and}&&
\Pi^{(1)}=\Pi_{0}+
\frac{\varepsilon i\partial_{p_{k}}f_{\varepsilon}}{E_{+}-E_{-}}
\sigma_{3}(\partial_{q^{k}}\Pi_{0})^{OD}\,,
\end{eqnarray*}
with
$$
\varepsilon^{-2\delta}\Pi_{1}=\frac{i\partial_{p_{k}}f_{1}}{E_{+}-E_{-}}
\sigma_{3}(\partial_{q^{k}}\Pi_{0})^{OD}\in
S_{u}\Big(\frac{1}{\langle \sqrt{\frac{\tau'}{\tau''}}q'\rangle
    \langle \sqrt{\tau'\tau''}p\rangle^{\infty}},g_{\tau};\mathcal{M}_{2}(\cz)\Big)\,.
$$
\noindent$\mathbf{n=2}$: Consider now
\begin{eqnarray*}
\varepsilon^{2}G_{2}&=&\Pi^{(1)}\sharp^{\varepsilon}
\Pi^{(1)}-\Pi^{(1)}\mod\mathcal{O}_{S}(\varepsilon^{3})\\
&=&\Pi_{0}\sharp^{\varepsilon} \Pi_{0}-\Pi_{0}+\varepsilon(\Pi_{0}\sharp^{\varepsilon}
\Pi_{1}+\Pi_{1}\sharp^{\varepsilon} \Pi_{0})-\Pi_{1}+\varepsilon^{2}\Pi_{1}\sharp^{\varepsilon}\Pi_{1}
\mod\mathcal{O}_{S}(\varepsilon^{3})
\,.
\end{eqnarray*}
According
to \eqref{eq.Wexpanexp}, with
$\Pi_{0}\sharp \Pi_{0}=\Pi_{0}$ and with
$$
\Pi_{0}\Pi_{1}+\Pi_{1}\Pi_{0}=
\Pi_{0}^{D}\Pi_{1}^{OD}+\Pi_{1}^{OD}\Pi_{0}^{D}
=
\Pi_{1}\,,
$$
 we
can take
$$
G_{2}=
\frac{1}{2i}\left[-\partial_{q^{k}}\Pi_{0}\partial_{p_{k}}\Pi_{1}+\partial_{p_{k}}\Pi_{1}\partial_{q^{k}}\Pi_{0}\right]
+ \Pi_{1}^{2}\,.
$$
The first term $-\partial_{q^{k}}\Pi_{0}\partial_{p_{k}}\Pi_{1}$, with
Einstein's summation rule,
equals
\begin{eqnarray*}
-\partial_{q^{k}}\Pi_{0}\partial_{p_{k}}\Pi_{1}
&=&-\frac{i\partial^{2}_{p_{k}p_{\ell}}f_{\varepsilon}}{(E_{+}-E_{-})}
(\partial_{q^{k}}\Pi_{0})^{OD}\sigma_{3}(\partial_{q^{\ell}}\Pi_{0})^{OD}
\\
&=&\frac{i\partial^{2}_{p_{k}p_{\ell}}f_{\varepsilon}}{(E_{+}-E_{-})}
\sigma_{3}(\partial_{q^{k}}\Pi_{0})(\partial_{q^{\ell}}\Pi_{0})\,,
\end{eqnarray*}
while the second term
$+(\partial_{p_{k}}\Pi_{1})(\partial_{q^{k}}\Pi_{0})$ gives the same result.\\
The third term $\Pi_{1}^{2}$ is given by
\begin{eqnarray*}
\Pi_{1}^{2}&=&-\frac{(\partial_{p_{k}}f_{\varepsilon})(\partial_{p_{\ell}}f_{\varepsilon})}{(E_{+}-E_{-})^{2}}
\sigma_{3}(\partial_{q^{k}}\Pi_{0})^{OD}\sigma_{3}
(\partial_{q^{\ell}}\Pi_{0})^{OD}
\\
&=&
\frac{(\partial_{p_{k}}f_{\varepsilon})(\partial_{p_{\ell}}f_{\varepsilon})}{(E_{+}-E_{-})^{2}}
(\partial_{q^{k}}\Pi_{0})(\partial_{q^{\ell}}\Pi_{0})\,.
\end{eqnarray*}
Hence the diagonal second order correction is given by
\begin{eqnarray*}
&&(E_{+}-E_{-})^{2}\Pi_{2}^{D}=-(E_{+}-E_{-})^{2}\sigma_{3}G_{2}
\\
&&=
-
\left[(\partial_{p_{k}}f_{\varepsilon})(\partial_{p_{\ell}}f_{\varepsilon})
+
  (E_{+}-E_{-})(\partial^{2}_{p_{k}p_{\ell}}f_{\varepsilon})\sigma_{3}\right]
(\partial_{q^{k}}\Pi_{0})(\partial_{q^{\ell}}\Pi_{0})
\end{eqnarray*}
and satisfies
$$
\varepsilon^{-2\delta}\Pi_{2}^{D}\in
S_{u}\Big(\frac{1}{\langle \sqrt{\frac{\tau'}{\tau''}}q'\rangle\langle \sqrt{\tau'\tau''}p\rangle^{\infty}}, g_{\tau};
\mathcal{M}_{2}(\cz)\Big)\,.
$$
Consider now $\Pi_{2}^{OD}$: By referring to \eqref{eq.FOD}, $F_{2}$ can be chosen as the
off-diagonal part of
$\varepsilon^{-2}\left[H\sharp^{\varepsilon}\Pi^{(1)}-\Pi^{(1)}\sharp^{\varepsilon}H\right]$
which, according to the previous steps and \eqref{eq.Wexpanexp},
equals
\begin{multline}
  \label{eq.F2}
-\frac{1}{8}\left[\partial^{2}_{p_{k},p_{\ell}}H\partial^{2}_{q^{k},q^{\ell}}\Pi_{0}
-
\partial^{2}_{q^{k},q^{\ell}}
 \Pi_{0}\partial^{2}_{p_{k},p_{\ell}}H\right]
\\
+
\frac{1}{2i}\left[\partial_{p_{k}}H\partial_{q^{k}}\Pi_{1}-\partial_{q^{k}}H\partial_{p_{k}}\Pi_{1}
- \partial_{p_{k}}\Pi_{1}\partial_{q^{k}}H
+\partial_{q^{k}}\Pi_{1}\partial_{p_{k}}H\right]
\mod\mathcal{O}_{S}(\varepsilon)\,.
\end{multline}
Since
$\partial^{2}_{p_{k},p_{\ell}}H=(\partial^{2}_{p_{k},p_{\ell}}f_{\varepsilon})$
as a scalar symbol
commutes with  $\partial^{2}_{q^{k}q^{\ell}}\Pi_{0}$ , the first term of \eqref{eq.F2}
 vanishes.\\
Similarly the factor $\partial_{q^{k}}\Pi_{1}$ appearing in the second
term of \eqref{eq.F2}
contains three terms
\begin{multline*}
\partial_{q^{k}}\Pi_{1}=
-i\frac{\partial_{q^{k}}(E_{+}-E_{-})}{(E_{+}-E_{-})^{2}}
\sigma_{3}(\partial_{p_{\ell}}f_{\varepsilon})(\partial_{q^{\ell}}\Pi_{0})+
\frac{i(\partial_{p_{\ell}}f_{\varepsilon})}{(E_{+}-E_{-})}
(\partial_{q^{k}}\sigma_{3})(\partial_{q^{\ell}}\Pi_{0})\\
+
\frac{i(\partial_{p_{\ell}}f_{\varepsilon})}{(E_{+}-E_{-})}
\sigma_{3}(\partial^{2}_{q^{k}q^{\ell}}\Pi_{0})\,,
\end{multline*}
where the second one is diagonal (Remember $\sigma_{3}(q,\tau,\varepsilon)=2\Pi_{0}(q,\tau,\varepsilon)-\Id_{\cz^{2}}$)\,.
Since $\partial_{p_{k}}H= \partial_{p_{k}}f_{\varepsilon}$ is
diagonal, we get
\begin{multline*}
\frac{1}{2i}\left[\partial_{p_{k}}H\partial_{q^{k}}\Pi_{1}
+\partial_{q^{k}}\Pi_{1}\partial_{p_{k}}H\right]^{OD}
\\
=
(\partial_{p_{k}}f_{\varepsilon})(\partial_{p_{\ell}}f_{\varepsilon})\sigma_{3}
\left(\partial_{q^{k}}\left[(E_{+}-E_{-})^{-1}\partial_{q^{\ell}}\Pi_{0}\right]\right)^{OD}
\,.
\end{multline*}
In the quantity
$-\partial_{q^{k}}H \partial_{p_{k}}\Pi_{1}-\partial_{p_{k}}\Pi_{1}\partial_{q^{k}}H$,
the derivatives
$$
\partial_{p_{k}}\Pi_{1}=\frac{i(\partial^{2}_{p_{k}p_{\ell}}f_{\varepsilon})}{(E_{+}-E_{-})}
\sigma_{3}(\partial_{q^{\ell}}\Pi_{0})^{OD}
$$
are off-diagonal factors, while
$$
\partial_{q^{k}}H= (\partial_{q^{k}}E_{+})
\Pi_{0}+(\partial_{q^{k}}E_{-})(1-\Pi_{0})
+(E_{+}-E_{-})(\partial_{q^{k}}\Pi_{0})^{OD}\,.
$$
With the two equalities,
\begin{eqnarray*}
&&(\partial_{q^{k}}E_{+})\Pi_{0}\sigma_{3}(\partial_{q^{\ell}}\Pi_{0})
+ \sigma_{3}(\partial_{q^{\ell}}\Pi_{0}) (\partial_{q^{k}}E_{+})\Pi_{0}
=(\partial_{q^{k}}E_{+})\sigma_{3}(\partial_{q^{\ell}}\Pi_{0})\,,
\\
&&
\hspace{-0.7cm}
 (\partial_{q^{k}}E_{-})(1-\Pi_{0})\sigma_{3}(\partial_{q^{\ell}}\Pi_{0})
+ \sigma_{3}(\partial_{q^{\ell}}\Pi_{0}) (\partial_{q^{k}}E_{-})(1-\Pi_{0})
= (\partial_{q^{k}}E_{-})\sigma_{3}(\partial_{q^{\ell}}\Pi_{0})\,,
\end{eqnarray*}
we get
$$
-\frac{1}{2i}\left[\partial_{q^{k}}H \partial_{p_{k}}\Pi_{1}+\partial_{p_{k}}\Pi_{1}\partial_{q^{k}}H\right]^{OD}
=
-\frac{\partial^{2}_{p_{k}p_{\ell}}f_{\varepsilon}}{2(E_{+}-E_{-})}(\partial_{q^{k}}(E_{+}+E_{-}))
\sigma_{3}(\partial_{q^{\ell}}\Pi_{0})\,.
$$
This leads to
\begin{multline*}
  F_{2}=
(\partial_{p_{k}}f_{\varepsilon})(\partial_{p_{\ell}}f_{\varepsilon})
\sigma_{3}\left(\partial_{q^{k}}\left[(E_{+}-E_{-})^{-1}\partial_{q^{\ell}}\Pi_{0}\right]\right)^{OD}
\\
-\frac{\partial^{2}_{p_{k}p_{\ell}}f_{\varepsilon}}{(E_{+}-E_{-})}(\partial_{q^{k}}(E_{+}+E_{-}))
\sigma_{3}(\partial_{q^{\ell}}\Pi_{0})\,.
\end{multline*}
and
\begin{multline*}
\Pi_{2}^{OD}=
-\frac{(\partial_{p_{k}}f_{\varepsilon})(\partial_{p_{\ell}}f_{\varepsilon})}{(E_{+}-E_{-})}
\left(\partial_{q^{k}}\left[(E_{+}-E_{-})^{-1}\partial_{q^{\ell}}\Pi_{0}\right]\right)^{OD}
\\
+\frac{\partial^{2}_{p_{k}p_{\ell}}f_{\varepsilon}}{2(E_{+}-E_{-})^{2}}(\partial_{q^{k}}(E_{+}+E_{-}))
(\partial_{q^{\ell}}\Pi_{0})\,,
\end{multline*}
with
$$
\varepsilon^{-2\delta}\Pi_{2}^{OD}\in
S_{u}\Big(\frac{1}{\langle \sqrt{\frac{\tau'}{\tau''}}q'\rangle\langle \sqrt{\tau'\tau''}p\rangle^{\infty}},g_{\tau};\mathcal{M}_{2}(\cz)\Big)\,.
$$
We have almost proved the
\begin{proposition}
\label{pr.proj2nd}
The pseudodifferential operator
$\hat\Pi(\varepsilon)=\hat\Pi(q,\varepsilon D_{q},\varepsilon)$ given by
\begin{eqnarray*}
&&  \Pi(q,p,\varepsilon)=\Pi_{0}(q,\varepsilon)+\varepsilon
\Pi_{1}(q,p,\varepsilon)+ \varepsilon^{2}\Pi_{2}(q,p,\varepsilon)\\
\text{with}&&
\Pi_{1}(q,p,\varepsilon)=
\Pi_{1}(q,p,\varepsilon)^{OD}=\frac{i
\partial_{p_{k}}f_{\varepsilon}}{(E_{+}-E_{-})}
\sigma_{3}(\partial_{q^{k}}\Pi_{0})
\,,\\
&&\Pi_{2}(q,p,\varepsilon)=
\Pi_{2}(q,p,\varepsilon)^{D}+\Pi_{2}(q,p,\varepsilon)^{OD}\,,\\
&&
\Pi_{2}(q,p,\varepsilon)^{D}
=
-\frac{1}{(E_{+}-E_{-})^{2}}\left[(\partial_{p_{k}}f_{\varepsilon})(\partial_{p_{\ell}}f_{\varepsilon})
\right.\\
&&\hspace{4.5cm}\left.+(E_{+}-E_{})(\partial^{2}_{p_{k}p_{\ell}}f_{\varepsilon})\right]
(\partial_{q^{k}}\Pi_{0})(\partial_{q^{\ell}}\Pi_{0})\,,\\
&&
\Pi_{2}(q,p,\varepsilon)^{OD}
= -\frac{(\partial_{p_{k}}f_{\varepsilon})(\partial_{p_{\ell}}f_{\varepsilon})}{(E_{+}-E_{-})}
\left(\partial_{q^{k}}\left[(E_{+}-E_{-})^{-1}\partial_{q^{\ell}}\Pi_{0}\right]\right)^{OD}
\\
&&\hspace{4.5cm}
+2\frac{\partial^{2}_{p_{k}p_{\ell}}f_{\varepsilon}}{(E_{+}-E_{-})^{3}}(\partial_{q^{k}}(E_{+}+E_{-}))
(\partial_{q^{\ell}}\Pi_{0})\,,\\
\text{and}&&
\varepsilon^{-2\delta}\Pi_{1}, \varepsilon^{-2\delta}\Pi_{2}\in
S_{u}\Big(\frac{1}{\langle \sqrt{\frac{\tau'}{\tau''}}q'\rangle\langle \sqrt{\tau'\tau''}p\rangle^{\infty}}, g_{\tau}; \mathcal{M}_{2}(\cz)\Big)\,,
\end{eqnarray*}
satisfies
$$
\hat \Pi\circ \hat \Pi= \hat \Pi+
\mathcal{O}(\varepsilon^{3+2\delta})
\quad,\quad
\hat \Pi^{*}=\hat \Pi \quad,\quad
\left[\hat H, \hat \Pi\right]= \mathcal{O}(\varepsilon^{3+2\delta})\,,
$$
in $\mathcal{L}(L^{2}(\rz^{d};\cz^{2}))$\,. Moreover the estimates in
$\mathcal{L}(L^{2}(\rz^{d};\cz^{2}))$ of
the remainder terms
 do not depend on the parameter $\tau=(\tau',\tau'')\in (0,1]^{2}$ and
 $\delta>0$, as soon as $\varepsilon\in (0,\varepsilon_{0})$\,.
\end{proposition}
\begindemonstration{}
The above construction gives immediately
$$
\hat \Pi\circ \hat \Pi= \hat \Pi+
\mathcal{O}(\varepsilon^{3})
\quad,\quad
\hat \Pi^{*}=\hat \Pi \quad,\quad
\left[\hat H, \hat \Pi\right]= \mathcal{O}(\varepsilon^{3})\,,
$$
The first improved estimates come from the fact that
$$
\hat \Pi\circ \hat \Pi- \hat \Pi
$$
contains only terms which are Moyal products with  a $\Pi_{1}$ or a
$\Pi_{2}$ factor, with cancellations up to the $\varepsilon^{2}$
coefficient.
Both of them have seminorms of
order $\varepsilon^{2\delta}$\,.\\
For the last one, this is a similar argument after decomposing
$$
\left[\hat H, \hat \Pi\right]
=
\left[\varepsilon^{2\delta}f_{1}(\sqrt{\tau'\tau''}\varepsilon D_{q}), \hat \Pi_{0}\right]
+ \left[\hat H, \varepsilon \hat \Pi_{1}+\varepsilon^{2}\hat\Pi_{2}\right]\,.
$$
\enddemonstration{}
The above result can be improved after considering what
happens at step $n=3$ when
$f_{\varepsilon,\tau}(p)=\varepsilon^{2\delta}f_{1}(\sqrt{\tau'\tau''}p)$
is at most quadratic w.r.t $p$ in
some region.
Before this, let us examine the remainders of order $3$ in
$\varepsilon$.\\
\noindent$\mathbf n=3:$ The remainder term
\begin{eqnarray*}
\varepsilon^{3}G_{3}&=&\Pi^{(2)}\sharp^{\varepsilon}\Pi^{(2)}-\Pi^{(2)}
\\
&=&
\varepsilon\left(\Pi_{0}\sharp^{\varepsilon}
  \Pi_{1}+\Pi_{1}\sharp^{\varepsilon}\Pi_{0}-\Pi_{1}\right)
+\varepsilon^{2}\left(\Pi_{1}\sharp^{\varepsilon}
  \Pi_{1}+\Pi_{0}\sharp^{\varepsilon}
 \Pi_{2}+\Pi_{2}\sharp^{\varepsilon}\Pi_{0}-\Pi_{2}\right)\\
&&+\varepsilon^{3}(\Pi_{1}\sharp^{\varepsilon}\Pi_{2}+\Pi_{2}\sharp^{\varepsilon}\Pi_{1})
+\varepsilon^{4}\Pi_{2}\sharp^{\varepsilon}\Pi_{2}\,.
\end{eqnarray*}
Using the construction of $\Pi_{1}$ and $\Pi_{2}$ and the fact that
$\varepsilon^{-2\delta}\Pi_{1}$ and $\varepsilon^{-2\delta}\Pi_{2}$
belong to $S_{u}\Big(\frac{1}{\langle
  \sqrt{\frac{\tau'}{\tau''}}q'\rangle\langle
  \sqrt{\tau'\tau''}p\rangle^{\infty}}, g_{\tau}; \mathcal{M}_{2}(\cz)\Big)$, the
expansion of the Moyal product \eqref{eq.Wexpan3}-\eqref{eq.Wexpanexp} tells us
\begin{equation}
  \label{eq.G3}
  \Pi^{(2)}\sharp^{\varepsilon}\Pi^{(2)}-\Pi^{(2)}= \varepsilon^{3}\left[A_{G}(\partial_{q}^{2}\Pi_{0},\partial_{p}^{2}\Pi_{1},\varepsilon)
  + B_{G}(\partial_{q}\Pi_{0}, \partial_{p}\Pi_{2},\varepsilon) \right] +\varepsilon^{3+4\delta}R_{G}
\end{equation}
where $R_{G}\in S_{u}\Big(\frac{1}{\langle
  \sqrt{\frac{\tau'}{\tau''}}q'\rangle\langle
  \sqrt{\tau'\tau''}p\rangle^{\infty}}, g_{\tau}; \mathcal{M}_{2}(\cz)\Big)$ and
$A_{G}(.,\varepsilon)$, $B_{G}(.,\varepsilon)$ have an asymptotic expansion in
terms of $\varepsilon$, of which all the terms are bilinear differential expressions
of their arguments.
For the commutator with $\hat{H}$, write
\begin{eqnarray*}
\left[H\sharp^{\varepsilon}\Pi^{(2)}-\Pi^{(2)}\sharp^{\varepsilon}H\right]
&=&\left[f_{\varepsilon}(p)\sharp^{\varepsilon}\Pi_{0}(q)-\Pi_{0}(q)\sharp^{\varepsilon}f_{\varepsilon}(p)\right]
\\
&&+
\left[f_{\varepsilon}(p)\sharp^{\varepsilon}(\varepsilon\Pi_{1}+\varepsilon^{2}\Pi_{2})-(\varepsilon\Pi_{1}+\varepsilon^{2}\Pi_{2})\sharp^{\varepsilon}f_{\varepsilon}(p)\right]
\\
&&
\hspace{-2cm}+
\left[(E_{+}\Pi_{0}+E_{-}(1-\Pi_{0}))\sharp^{\varepsilon}\Pi^{(2)}
-\Pi^{(2)}\sharp^{\varepsilon}(E_{+}\Pi_{0}+E_{-}(1-\Pi_{0}))\right]
\end{eqnarray*}
After eliminating all the terms which are cancelled while constructing
$\Pi_{1}$ and $\Pi_{2}$, the contributions of all three terms of the
right-hand side can be analyzed. The contribution of the third term is
similar to what we got for $G_{3}$:
$$
\varepsilon^{3}\left[A_{H}(\partial_{q}^{2}(E_{\pm}\Pi_{0}), \partial_{p}^{2}\Pi_{1},\varepsilon)
+
B_{H}(\partial_{q}(E_{\pm}\Pi_{0}), \partial_{p}\Pi_{2},\varepsilon)\right]
+\varepsilon^{3+4\delta}R_{H,1}
$$
with $R_{H,1}\in S_{u}\Big(\frac{1}{\langle
  \sqrt{\tau'\tau''}p\rangle^{\infty}}, g_{\tau};
\mathcal{M}_{2}(\cz)\Big)$\,.
With the uniform estimate of $\varepsilon^{-2\delta}\Pi_{1}$,
$\varepsilon^{-2\delta}\Pi_{2}$ and
$\varepsilon^{-2\delta}f_{\varepsilon}$, the contribution of the
second term is estimated as $\varepsilon^{3+4\delta}R_{H,2}$ with
$R_{H,2}\in S_{u}\Big(\frac{1}{\langle\sqrt{\frac{\tau'}{\tau''}}q'\rangle\langle
  \sqrt{\tau'\tau''}p\rangle^{\infty}}, g_{\tau};
\mathcal{M}_{2}(\cz)\Big)$\,. The contribution of the first term is
$\varepsilon^{3}C_{H}(\partial_{p}^{3}f_{\varepsilon}, \partial_{q}^{3}\Pi_{0},\varepsilon)$
and we get
\begin{multline}
  \label{eq.almF3}
  \left[H\sharp^{\varepsilon}\Pi^{(2)}-\Pi^{(2)}\sharp^{\varepsilon}H\right]
=
\varepsilon^{3}\left[
C_{H}(\partial_{p}^{3}f_{\varepsilon}, \partial_{q}^{3}\Pi_{0},\varepsilon)
+
\right.
\\
\left.A_{H}(\partial_{q}^{2}(E_{\pm}\Pi_{0}), \partial_{p}^{2}\Pi_{1},\varepsilon)
+
B_{H}(\partial_{q}(E_{\pm}\Pi_{0}), \partial_{p}\Pi_{2},\varepsilon)
\right]
+\varepsilon^{3+4\delta}R_{H}
\end{multline}
where the expansions of $A_{H},B_{H}$ and $C_{H}$ w.r.t to $\varepsilon$ have
terms which are bilinear differential expression of their arguments,
and the remainder $R_{H}$ belongs to
$S_{u}\Big(\frac{1}{\langle
  \sqrt{\tau'\tau''}p\rangle^{\infty}}, g_{\tau};
\mathcal{M}_{2}(\cz)\Big)$\,.
\begin{proposition}
\label{pr.preci}
  With
  $f_{\varepsilon}(p)=\varepsilon^{2\delta}f_{1}(\sqrt{\tau'\tau''}p)$,
  assume that the third differential $\partial_{p}^{3}f_{1}$ vanishes
  in $\left\{|p|\leq r\right\}$ and fix $r'\in (0,r)$\,. Then the remainders of
  Proposition~\ref{pr.proj2nd} equal
  \begin{eqnarray}
    \label{eq.Piprec}
&&    \hat{\Pi}\circ \hat{\Pi}-\hat{\Pi}=
\varepsilon^{3+2\delta}\hat{R}_{G,1}+
\varepsilon^{3+4\delta}\hat{R}_{G,2}\\
\label{eq.Hprec}
&& \left[\hat{H},\hat{\Pi}\right]= \varepsilon^{3+2\delta}\hat{R}_{H,1}+
\varepsilon^{3+4\delta}\hat{R}_{H,2}
  \end{eqnarray}
where $R_{G,1\text{or}2}$ belong to $OpS_{u}\Big(\frac{1}{\langle\sqrt{\frac{\tau'}{\tau''}}q'\rangle\langle
  \sqrt{\tau'\tau''}p\rangle^{\infty}}, g_{\tau};
\mathcal{M}_{2}(\cz)\Big)$, $R_{H,1\text{or}2}$ belong to $ OpS_{u}\Big(\frac{1}{\sqrt{\tau'\tau''}p\rangle^{\infty}}, g_{\tau};
\mathcal{M}_{2}(\cz)\Big)$ and
$R_{G\text{or}H,1}\equiv 0$ in
$\left\{|\sqrt{\tau'\tau''}p|<r'\right\}$\,. Those estimates are
uniform for $\tau\in (0,1]^{2}$, $\delta\in [0,\delta_{0}]$ and
$\varepsilon\in (0,\varepsilon_{0})$\,.
\end{proposition}
\begindemonstration{}
After noticing that the symbol $\Pi_{1}$ is a linear expression in
$\partial_{p}f_{\varepsilon}$  while the symbol
 $\Pi_{2}$ is the sum of a linear expression of $\partial_{p}^{2}f_{\varepsilon}$
 and quadratic expression in $\partial_{p}f_{\varepsilon}$, the
 identities \eqref{eq.G3} and \eqref{eq.almF3} imply
 \begin{eqnarray*}
    \Pi^{(2)}\sharp^{\varepsilon}\Pi^{(2)}-\Pi^{(2)}
&=& \varepsilon^{3+2\delta+N}R_{N}+ \varepsilon^{3+4\delta}R_{G}\,,\\
\left[H\sharp^{\varepsilon}\Pi^{(2)}-\Pi^{(2)}\sharp^{\varepsilon}H\right]
&=&
\varepsilon^{3+2\delta+ N}R_{N}+ \varepsilon^{3+4\delta}R_{H}\,,
 \end{eqnarray*}
for an arbitrary large $N\in\nz$,
in $\left\{\sqrt{\tau'\tau''}|p|< r\right\}$\,. Choose\footnote{Here
  the estimates become $\delta$-dependent, because a large $\delta$
  requires a large $N$. It is uniformly controlled when $\delta\leq \delta_{0}$.} $N\geq
2\delta$ and take a cut-off function $\chi\in
\mathcal{C}^{\infty}_{0}\left(\left\{|p|<r\right\}\right)$  such that
$\chi\equiv 1$ in a neighborhood $\left\{|p|\leq r'\right\}$\,. Writing
for the symbol
$$
S=
\Pi^{(2)}\sharp^{\varepsilon}\Pi^{(2)}-\Pi^{(2)}\quad\text{or}\quad
S=\left[H\sharp^{\varepsilon}\Pi^{(2)}-\Pi^{(2)}\sharp^{\varepsilon}H\right]\,,$$
$$
S=S\times\chi(\sqrt{\tau'\tau''}p)+ S \times(1-\chi(\sqrt{\tau'\tau''}p))\,,
$$
yields the result.
\enddemonstration{}
\subsection{Unitaries and effective Hamiltonian}
\label{se.projunit}
 We strengthen a little bit the assumptions
 \eqref{eq.hypEPi1}-\eqref{eq.hypEPi2}, with the condition
 \begin{eqnarray}
   \label{eq.hypu0}
&&\Pi_{0}(q,\tau,\varepsilon)=u_{0}(q,\tau,\varepsilon)P_{+}u_{0}(q,\tau,\varepsilon)^{*}
\\
\nonumber
\text{with}
&&
P_{+}=
\begin{pmatrix}
  1&0\\0&0
\end{pmatrix}
\,,\; u_{0}=(u_{0}^{*})^{-1}\in S(1,g_{q,\tau}; \mathcal{M}_{2}(\cz))\,,
\end{eqnarray}
fulfilled in our example. The operator $\hat{u}_{0}$ is nothing but
the local unitary transformation $u_{0}(q,\tau)$, on $L^{2}(\rz^{d};\cz^{2})$\,.\\
With the approximate projection $\hat\Pi=\hat{\Pi}(q,\varepsilon
D_{q},\tau,\varepsilon)$ given in Proposition~\ref{pr.proj2nd},
Proposition~\ref{pr.appproj}
tells us that a true
orthogonal projection $\hat{P}$ can be associated when
$\varepsilon_{0}$ is chosen small enough, by taking
\begin{equation}
  \label{eq.defPeps}
\hat{P}=\frac{1}{2i\pi}\int_{|z-1|=1/2}(z-\hat\Pi)^{-1}~dz\,,
\end{equation}
with
\begin{eqnarray}
\label{eq.hatP}
&&
\hat{P}=P(q,\varepsilon D_{q},\tau,\varepsilon)\quad,\quad
P\in S_{u}(1,g_{\tau};\mathcal{M}_{2}(\cz))\,,\\
  \label{eq.PPi}
&&
\hspace{-3cm}
P(q,p,\tau,\varepsilon)-\Pi(q,p,\tau,\varepsilon)=\varepsilon^{3+2\delta}R_{1}(q,p,\tau,\varepsilon)
+\varepsilon^{3+4\delta}R_{2}(q,p,\tau,\varepsilon)\,,\;\\
\nonumber
\text{with}&&
R_{1},R_{2}\in S_{u}\Big(\frac{1}{\langle \sqrt{\frac{\tau'}{\tau''}}q'\rangle\langle \sqrt{\tau'\tau''}p\rangle^{\infty}},g_{\tau};\mathcal{M}_{2}(\cz)\Big)\,,\\
&&\nonumber
 \hat{P}\circ
\hat{P}=\hat{P}=\hat{P}^{*}\,,
\\
&&
\label{eq.comHP}
\left[\hat{H},\hat{P}\right]=\varepsilon^{3+2\delta}
\hat{C}_{1}(\tau,\varepsilon)+\varepsilon^{3+4\delta}\hat{C}_{2}(\tau,\varepsilon)\,,\\
\nonumber\quad\text{with}&&
C_{1},C_{2}\in S_{u}(1,g_{\tau};\mathcal{M}_{2}(\cz))\,,\\
\text{and}
&&
\label{eq.Pu0P+}
\left\|\hat{P}-\hat{u}_{0}P_{+}\hat{u}_{0}^{*}\right\|_{\mathcal{L}(L^{2})}\leq C\varepsilon\,.
\end{eqnarray}
For a general $f_{1}\in \mathcal{C}^{\infty}_{0}$, $R_{2}$ and $C_{2}$ are included in the main
remainder term. When $f_{1}$ is quadratic in $\left\{|p|<r\right\}$,
then one can assume that $R_{1}$ and $C_{1}$ vanishes in
$\left\{\sqrt{\tau'\tau''}|p|<r'\right\}$ for $r'<r$, according to
Proposition~\ref{pr.preci} and Proposition~\ref{pr.appproj}\\
Instead of constructing unitaries between $\hat{P}(\tau,\varepsilon)$ and $P_{+}$
by the induction  presented in
\cite{PST} and similar to \eqref{eq.defGn1},
\eqref{eq.Pin1D},\eqref{eq.defFn1},\eqref{eq.Pin1OD}, we use like in \cite{MaSo}  Nagy's formula
(\cite{NeSo}, \cite{MaSo}, \cite{PST})
\begin{eqnarray}
\nonumber
  &&P_{2}=WP_{1}W^{*}\quad,\quad W^{*}W=WW^{*}=1\,,\\
\text{with}
&&
\label{eq.nagy}
W=(1-(P_{2}-P_{1})^{2})^{-1/2}\left[P_{2}P_{1}+(1-P_{2})(1-P_{1})\right]\,,\\
\nonumber
\text{when}&&
P_{j}=P_{j}^{2}=P_{j}^{*}\quad \text{for}\, j=1,2\,,\quad
\text{and}\,\|P_{2}-P_{1}\|_{\mathcal{L}}< 1\,,
\end{eqnarray}
easier to handle for direct second order computations in our case.
\begin{proposition}
\label{pr.projunit}
With the definitions \eqref{eq.hypu0} and \eqref{eq.defPeps} after
Proposition~\ref{pr.proj2nd}, there exists a unitary operator
$\hat{U}$ on $L^{2}(\rz;\cz^{2})$ such that
\begin{eqnarray*}
&& \hat P=\hat{U}P_{+}\hat{U}^{*}\,,\\
&& \hat{U}=U(q,\varepsilon D_{q},\tau,\varepsilon)\quad,\quad
U\in S_{u}(1,g_{\tau};\mathcal{M}_{2}(\cz))\,,
\\
\text{with}&&U(q,p,\tau,\varepsilon)=u_{0}(q,\tau,\varepsilon)+\varepsilon u_{1}(q,p,\tau,\varepsilon)+
\varepsilon^{2}u_{2}(q,p,\tau,\varepsilon)\,,
\\
&&
u_{1}(q,p,\varepsilon)=-\frac{i\partial_{p_{k}}f_{\varepsilon}}{(E_{+}-E_{-})}(\partial_{q^{k}}\Pi_{0})
u_{0}\,,\\
\text{and}&& \varepsilon^{-2\delta}u_{1},
\varepsilon^{-2\delta}u_{2}\in S_{u}\Big(
\frac{1}{\langle\sqrt{\frac{\tau'}{\tau''}}\rangle\langle\sqrt{\tau'\tau''}p\rangle^{\infty}},
g_{\tau};\mathcal{M}_{2}(\cz)\Big)\,.
\end{eqnarray*}
Moreover, when $f_{1}$ is a quadratic function in
$\left\{|p|<r\right\}$ and $r'$ is fixed in $ (0,r)$, the term $u_{2}$ can be decomposed into
$$
u_{2}(q,p,\tau,\varepsilon)=
\varepsilon^{2\delta}v_{2}(q,p,\tau,\varepsilon)+
\varepsilon^{4\delta}\widetilde{v_{2}}(q,p,\tau,\varepsilon)
$$
where $v_{2}$ does not depend on $p$ in
$\left\{\sqrt{\tau'\tau''}|p|<r'\right\}$ and $v_{2},
\widetilde{v}_{2}$ belong to the symbol class
$S_{u}\Big(\frac{1}{\langle\sqrt{\frac{\tau'}{\tau''}}\rangle\langle\sqrt{\tau'\tau''}p\rangle^{\infty}},
g_{\tau}; \mathcal{M}_{2}(\cz)\Big)$\,.
\end{proposition}
\begindemonstration{}
The notation $\hat{R}$ will denote a generic remainder term of
the form  $\hat{R}=R(q,\varepsilon D_{q},\tau,\varepsilon)$ with
$R\in S_{u}\Big(\frac{1}{\langle
  \sqrt{\frac{\tau'}{\tau''}}q'\rangle\langle\sqrt{\tau'\tau''}
  p\rangle^{\infty}},g_{\tau};\mathcal{M}_{2}(\cz)\Big)$\,. The notation
$\hat{\underline{R}}$ is used for a symbol $\underline{R}$, like $R$ but which vanishes
around $\left\{\sqrt{\tau'\tau''}|p|<r'\right\}$\,.
We apply Nagy's formula \eqref{eq.nagy} with
$P_{1}=\hat{u}_{0}P_{+}\hat{u}_{0}^{*}=\Pi_{0}(q,\tau,\varepsilon)$
and $P_{2}=\hat{P}(\tau,\varepsilon)$ with
$$
\hat P=\Pi_{0}(q,\tau,\varepsilon)+\varepsilon \Pi_{1}(q,\varepsilon
D_{q},\tau,\varepsilon)+ \varepsilon^{2}\Pi_{2}(q,\varepsilon
D_{q},\tau,\varepsilon)+ \varepsilon^{3+2\delta}\underline{\hat{R}}+\varepsilon^{3+4\delta}\hat{R}\,.
$$
In the expression of $\hat{W}(\varepsilon)$ given by \eqref{eq.nagy},
 the first factor is nothing but
$$
(1-(\hat{P}-\Pi_{0}(q,\tau,\varepsilon))^{2})^{-1/2}=1+
\frac{\varepsilon^{2}}{2}(\Pi_{1})^{2}(q,\varepsilon
D_{q},\tau,\varepsilon)
+
\varepsilon^{3+4\delta}\hat{R}\,,
$$
owing to $P-\Pi_{0}=\mathcal{O}_{S}(\varepsilon^{2\delta})$\,.
In the factor $[P_{2}P_{1}+(1-P_{2})(1-P_{1})]$, the first term equals
\begin{multline*}
\hat{P}\circ \Pi_{0}(q,\varepsilon)=
\Pi_{0}(q,\varepsilon)+\varepsilon\Pi_{1}(q,\tau,\varepsilon
D_{q},\tau,\varepsilon)\circ \Pi_{0}(q,\tau,\varepsilon)
\\
+
\varepsilon^{2}\Pi_{2}(q,\varepsilon
D_{q},\tau,\varepsilon)\circ \Pi_{0}(q,\tau,\varepsilon)
 +\varepsilon^{3+2\delta}\underline{\hat{R}}+\varepsilon^{3+4\delta}\hat{R}\,,
\end{multline*}
while the second term is
\begin{multline*}
  (1-\hat{P})\circ(1-\Pi_{0}(q,\tau,\varepsilon))=(1-\Pi_{0}(q,\tau,\varepsilon))
-\varepsilon\Pi_{1}(q,\varepsilon
D_{q},\tau,\varepsilon)\circ (1-\Pi_{0}(q,\tau,\varepsilon))
\\
-
\varepsilon^{2}\Pi_{2}(q,\varepsilon
D_{q},\tau,\varepsilon)\circ (1-\Pi_{0}(q,\tau\varepsilon))
 +\varepsilon^{3+2\delta}\underline{\hat{R}}+\varepsilon^{3+4\delta}\hat{R}\,.
\end{multline*}
Hence we get
\begin{multline*}
  \hat{W}=
[P_{2}P_{1}+(1-P_{2})(1-P_{1})]=1+\varepsilon\Pi_{1}(q,\varepsilon
D_{q},\tau,\varepsilon)\circ
\sigma_{3}(q,\tau)
\\
+\varepsilon^{2}\Pi_{2}(q,\varepsilon
D_{q},\tau,\varepsilon)\circ \sigma_{3}(q,\tau)
 +\varepsilon^{3+2\delta}\underline{\hat{R}}+\varepsilon^{3+4\delta}\hat{R}\,.
\end{multline*}
The operator $\hat{U}(\varepsilon)$ is given by
$\hat{U}(\varepsilon)=\hat{W}(\varepsilon)\circ\hat{u}_{0}$\,.
The
 semiclassical calculus  recalled in
\eqref{eq.Wexpan3}-\eqref{eq.Wexpanexp} yields the result.
In the decomposition of $u_{2}$, when $f_{1}$ is quadratic around $0$,
the terms which are linear in
$\varepsilon^{2\delta}$ come with the second derivative of
$f_{\varepsilon}$, which does not depend on $p$\,.
\enddemonstration{}

\begin{proposition}
\label{pr.hameff}
Introduce the
notation for $k\in \left\{1,\ldots,d\right\}$
\begin{equation*}
 u_{0}^{*}(\partial_{q^{k}}u_{0})=-i
\begin{pmatrix}
    A_{k}&X_{k}\\
 \overline{X_{k}} & -A_{k}
  \end{pmatrix}\,.
\end{equation*}
If $\hat{U}$ is the unitary operator introduced in
Proposition~\ref{pr.projunit}, the conjugated Hamiltonian
$\hat{U}(\varepsilon)^{*}\hat{H}(\varepsilon)\hat{U}(\varepsilon)$ equals
\begin{equation}
  \label{eq.almdiag}
\hat{U}^{*}\hat{H}\hat{U}
=
\begin{pmatrix}
  \hat{h}_{+}&0\\
0& \hat{h}_{-}
\end{pmatrix}
+ \varepsilon^{3+2\delta}\hat{R}_{1}+\varepsilon^{3+4\delta}\hat{R}_{2}
\end{equation}
where the remainder terms $\hat{R}_{1,2}=R_{1,2}(q,\varepsilon
D_{q},\tau,\varepsilon)$ belong to
$OpS_{u}\left(1,g_{\tau};\mathcal{M}_{2}(\cz)\right)$ and additionally
$R_{1}\equiv 0$ in $\left\{\sqrt{\tau'\tau''}|p|<r'\right\}$ when $f_{1}$
is quadratic in $\left\{|p|<r\right\}$, with $r'<r$\,.\\
The symbol $h_{+}$ and $h_{-}$ are given by
\begin{eqnarray}
\nonumber
h_{+}&=&
f_{\varepsilon}(p)+E_{+}(q)-\varepsilon
(\partial_{p}f_{\varepsilon}).A
+
\frac{\varepsilon^{2}}{2}(\partial^{2}_{p_{k}p_{\ell}}f_{\varepsilon})A_{k}A_{\ell}
\\
\label{eq.h+}
&&\hspace{1cm}
+
\frac{\varepsilon^{2}(\partial^{2}_{p_{k}p_{\ell}}f_{\varepsilon})}
{2}X_{k}\overline{X}_{\ell}
+
\frac{\varepsilon^{2}(\partial_{p_{k}}f_{\varepsilon})(\partial_{p_{\ell}}f_{\varepsilon})}
{E_{+}-E_{-}}X_{k}\overline{X}_{\ell}
\,,
\\
\nonumber
h_{-}&=&
f_{\varepsilon}(p)+E_{-}(q)+\varepsilon
(\partial_{p}f_{\varepsilon}).A
+
\frac{\varepsilon^{2}}{2}(\partial^{2}_{p_{k}p_{\ell}}f_{\varepsilon})A_{k}A_{\ell}
\\
\label{eq.h-}
&&\hspace{1cm}
+\frac{\varepsilon^{2}(\partial^{2}_{p_{k}p_{\ell}}f_{\varepsilon})}{2}
\overline{X_{k}}X_{\ell}
-
\frac{\varepsilon^{2}(\partial_{p_{k}}f_{\varepsilon})(\partial_{p_{\ell}}f_{\varepsilon})}
{E_{+}-E_{-}}\overline{X_{k}}X_{\ell}\,.
\end{eqnarray}
\end{proposition}
\begindemonstration{}
From the semiclassical calculus, we already know that
$\hat{U}^{*}\hat{H}\hat{U}$ is a
semiclassical operator with a symbol in
$S_{u}(\langle \sqrt{\frac{\tau'}{\tau''}}q'\rangle,
g_{\tau};\mathcal{M}_{2}(\cz))$\,. Its off-diagonal part equals
\begin{eqnarray*}
  (1-P_{+})\hat{U}^{*}\hat{H}\hat{U}P_{+}
&+&
P_{+}\hat{U}^{*}\hat{H}\hat{U}(1-P_{+})
\\
&& =
\hat{U}\left[
  \hat{P},\left[\hat{P},\hat{H}\right]\right] \hat{U}\,.
\end{eqnarray*}
The almost diagonal form \eqref{eq.almdiag}
of $\hat{U}^{*}\hat{H}\hat{U}$
is then a consequence of \eqref{eq.comHP}\,.\\
For the second result, it is necessary to compute the diagonal part of the symbol
$U^{*}\sharp^{\varepsilon}H\sharp^{\varepsilon}U$
up to $\mathcal{O}(\varepsilon^{3+2\delta})$ in
$S_{u}(1,g_{\tau};\mathcal{M}_{2})$\,.
Let us compute the diagonal part of
$$
\mathcal{B}:=(u_{0}^{*}+\varepsilon u_{1}^{*}+\varepsilon^{2}u_{2}^{*})\sharp^{\varepsilon}H\sharp^{\varepsilon}(u_{0}+\varepsilon
u_{1}+\varepsilon^{2}u_{2})\,,
$$
or equivalently $(u_{0}\mathcal{B}u_{0}^{*})^{D}$ with
our notations.\\
Since $u_{0}=u_{0}(q,\tau,\varepsilon)$ and
$H=\varepsilon^{2\delta}f_{1}(\sqrt{\tau'\tau''}p)+\mathcal{V}(q,\tau,\varepsilon)$
\,, the first Moyal
product equals according to \eqref{eq.Wexpan3}-\eqref{eq.Wexpanexp},
\begin{multline*}
(u_{0}^{*}+\varepsilon u_{1}^{*}+\varepsilon^{2}u_{2}^{*})\sharp^{\varepsilon}H
=
u_{0}^{*}H+\varepsilon u_{1}^{*}H+\varepsilon^{2}u_{2}^{*}H
\\
+\frac{\varepsilon}{2i}\left\{u_{0}^{*},H\right\}+
\frac{\varepsilon^{2}}{2i}\left\{u_{1}^{*},H\right\}
-\frac{\varepsilon^{2}}{8}(\partial^{2}_{q^{k}q^{\ell}}u_{0}^{*})(\partial^{2}_{p_{k}p_{\ell}}H)
+ \varepsilon^{3+2\delta}R
\\
=u_{0}^{*}H+\varepsilon
\left[u_{1}^{*}H-
\frac{\partial_{p_{k}}f_{\varepsilon}}{2i}\partial_{q^{k}}u_{0}^{*}\right]
+\varepsilon^{2}
\left[u_{2}^{*}H
+ \frac{1}{2i}\left\{u_{1}^{*},H\right\}
\right.
\\
\left.
-\frac{1}{8}(\partial^{2}_{q^{k}q^{\ell}}u_{0}^{*})(\partial^{2}_{p_{k}p_{\ell}}f_{\varepsilon})
\right]+ \varepsilon^{3+2\delta}\underline{R}+\varepsilon^{3+4\delta}R\,,
\end{multline*}
where $R$ and $\underline{R}$ denote  generic element of
$S_{u}(1,g_{\tau};\mathcal{M}_{2}(\cz))$, with the additional property
that $\underline{R}$ vanishes in
$\left\{\sqrt{\tau'\tau''}|p|<r'\right\}$ when $f_{1}$ is quadratic
in $\left\{|p|\leq r\right\}$ with $r'<r$. The reason for the possible
decomposition of the remainder, comes again from the fact that  the
third order  remainder term, proportional to $\varepsilon^{3+2\delta}$,
arises with the third derivative of $f_{\varepsilon}$, the second
derivative w.r.t $p$ of $u_{1}$ and the first derivative w.r.t $p$ of
$u_{2}$\,.\\
In the same way, the complete expression of
$\mathcal{B}$ is given by
\begin{eqnarray*}
&&\hspace{-0.5cm}\mathcal{B}(\varepsilon)=u_{0}^{*}Hu_{0}+\varepsilon\left[u_{0}^{*}H
  u_{1}+u_{1}^{*}Hu_{0}-\frac{(\partial_{p_{k}}f_{\varepsilon})}{2i}(\partial_{q^{k}}u_{0}^{*})u_{0} +
  \frac{(\partial_{p_{k}}f_{\varepsilon})}{2i}u_{0}^{*}(\partial_{q^{k}}u_{0}) \right]
\\
&&\;
+\varepsilon^{2}\left[u_{2}^{*}H u_{0}+
  u_{0}^{*}H u_{2}
+\frac{1}{2i}\left\{u_{0}^{*}H,u_{1}\right\}+
\frac{1}{2i}\left\{u_{1}^{*}H, u_{0}\right\}
+\frac{1}{4}\left\{(\partial_{p_{k}}f_{\varepsilon})(\partial_{q^{k}}u_{0}^{*}),
u_{0}\right\}
\right.\\
&&\quad
+u_{1}^{*}Hu_{1}-\frac{(\partial_{p_{k}}f_{\varepsilon})}{2i}(\partial_{q^{k}}u_{0}^{*})u_{1}
+\frac{1}{2i}\left\{u_{1}^{*},H\right\}u_{0}
\\
&&\quad
\left.
-\frac{1}{8}(\partial^{2}_{q^{k}q^{\ell}}u_{0}^{*})(\partial^{2}_{p_{k}p_{\ell}}f_{\varepsilon}))u_{0}
-\frac{1}{8}u_{0}^{*}(\partial^{2}_{p_{k}p_{\ell}}f_{\varepsilon})(\partial^{2}_{q^{k}q^{\ell}}u_{0})
\right]+\varepsilon^{3+2\delta}\underline{R}+\varepsilon^{3+4\delta}R\,.
\end{eqnarray*}
By recalling that $(u_{1}u_{0}^{*})^{D}=0$ by
Proposition~\ref{pr.projunit}, we get
\begin{multline*}
  (u_{0}\mathcal{B}(\varepsilon)u_{0}^{*})^{D}
=
H+\varepsilon \frac{i}{2}
(\partial_{p_{k}}f_{\varepsilon})\left[u_{0}(\partial_{q^{k}}u_{0}^{*})
-(\partial_{q^{k}}u_{0})u_{0}^{*}
\right]^{D}
\\
+\varepsilon^{2} \mathcal{B}_{2}^{D}+\varepsilon^{3+2\delta}\underline{R}+\varepsilon^{3+4\delta}R\,.
\end{multline*}
where $\mathcal{B}_{2}^{D}$ is made of several terms to be analyzed.
We need the relations
\begin{eqnarray}
\label{eq.deru0u0}
&&\partial_{q^{j}}u_{0}^{*}=-u_{0}^{*}(\partial_{q^{j}}u_{0})u_{0}^{*}\,,\\
\nonumber
&&\partial_{q^{j}q^{j'}}^{2}u_{0}^{*}
=u_{0}^{*}(\partial_{q^{j'}}u_{0})
u_{0}^{*}(\partial_{q^{j}}u_{0})u_{0}^{*}
+
u_{0}^{*}(\partial_{q^{j}}u_{0})
u_{0}^{*}(\partial_{q^{j'}}u_{0})u_{0}^{*}
-u_{0}^{*}(\partial^{2}_{q^{j}q^{j'}}u_{0})u_{0}^{*}
\\
\label{eq.der2u0u0}
&&\qquad=
-(\partial_{q^{j'}}u_{0}^{*})(\partial_{q^{j}}u_{0})u_{0}^{*}
-(\partial_{q^{j}}u_{0}^{*})(\partial_{q^{j'}}u_{0})u_{0}^{*}
-u_{0}^{*}(\partial^{2}_{q^{j}q^{j'}}u_{0})u_{0}^{*}\,,
\,\\
\label{eq.deunitOD}
&&\text{and}\qquad
\left[u_{0}(\partial_{q}u_{0}^{*})\right]^{OD}=- (\partial_{q}\Pi_{0})\sigma_{3}\,,
\end{eqnarray}
coming from $u_{0}^{*}u_{0}=1$ and the differentiation of $u_{0}^{*}\Pi_{0}u_{0}=P_{+}$\,.\\
For example, the first one simplifies the
$\mathcal{O}(\varepsilon)$-term into
$$
 (u_{0}\mathcal{B}(\varepsilon)u_{0}^{*})^{D}
=
H-\varepsilon i
(\partial_{p_{k}}f_{\varepsilon})
\left[(\partial_{q^{k}}u_{0})u_{0}^{*}
\right]^{D}
+\varepsilon^{2} \mathcal{B}_{2}^{D}
+
\varepsilon^{3+2\delta}\underline{R}+\varepsilon^{3+4\delta}R\,.
$$
Many cancellations appear after assembling all the terms in
$\mathcal{B}_{2}^{D}$\,. We need accurate expressions for all of them:
\begin{itemize}
\item By using again $(u_{1}u_{0}^{*})^{D}=0$,
the term $\left[u_{0}u_{2}^{*}H+H
  u_{2}u_{0}^{*}+u_{0}u_{1}^{*}Hu_{1}u_{0}^{*}\right]^{D}$
equals
\begin{multline*}
(f_{\varepsilon}+E_{+})\Pi_{0}(u_{0}u_{2}^{*}+u_{2}u_{0}^{*})\Pi_{0}+
(f_{\varepsilon}+E_{-})(1-\Pi_{0})(u_{0}u_{2}^{*}+u_{2}u_{0}^{*})(1-\Pi_{0})
\\
+
(f_{\varepsilon}+E_{-})\Pi_{0}(u_{0}u_{1}^{*}u_{1}u_{0}^{*})\Pi_{0}+
(f_{\varepsilon}+E_{+})(1-\Pi_{0})(u_{0}u_{1}^{*}u_{1}u_{0}^{*})(1-\Pi_{0})\,.
\end{multline*}
In the relation
$$
u_{0}^{*}u_{2}+u_{2}^{*}u_{0}+u_{1}^{*}u_{1}+\frac{1}{2i}\left\{u_{0}^{*},u_{1}\right\}
+
\frac{1}{2i}\left\{u_{1}^{*},u_{0}\right\}=\varepsilon^{1+2\delta}
\underline{R_{1}}
+\varepsilon^{1+4\delta}R_{1}\,,
$$
the remainder terms satisfy $R_{1},\underline{R}_{1}\in
S_{u}\Big(\frac{1}{\langle \sqrt{\frac{\tau'}{\tau''}}q'\rangle\langle
  \sqrt{\tau'\tau''}p\rangle^{\infty}},g_{\tau};\mathcal{M}_{2}(\cz)\Big)$,
with the same convention as for $R,\underline{R}$\,. This identity is
obtained by  writing that the $\mathcal{O}(\varepsilon^{2})$ remainder
of $U^{*}\sharp^{\varepsilon}U-1$ vanishes and
by noticing that the remainder $\underline{R_{1}}+\varepsilon^{2\delta}R_{1}$
involves second derivatives of $u_{1}$ w.r.t $p$ and first derivatives
of $u_{2}$ w.r.t $p$\,.
We obtain
\begin{multline*}
\left[u_{0}u_{2}^{*}H+H
  u_{2}u_{0}^{*}+u_{0}u_{1}^{*}Hu_{1}u_{0}^{*}\right]^{D}
=(E_{-}-E_{+})(u_{0}u_{1}^{*}u_{1}u_{0}^{*})^{D}\sigma_{3}\\
-\frac{1}{2i}H\left[u_{0}\left\{u_{0}^{*},u_{1}\right\}u_{0}^{*}+u_{0}\left\{u_{1}^{*},u_{0}\right\}u_{0}^{*}\right]^{D}
+\varepsilon^{1+2\delta} \underline{\tilde{R}} + \varepsilon^{1+4\delta}\tilde{R}\,.
\end{multline*}
with $\tilde{R},\underline{\tilde{R}}\in
S_{u}(1,g_{\tau};\mathcal{M}_{2}(\cz))$, again with the same convention.\\
Again with
$(\partial_{p}u_{1})u_{0}^{*}=(\partial_{p}u_{1}u_{0}^{*})^{OD}=\frac{(\partial^{2}_{p,p_{k}}f)}{i(E_{+}-E_{-})}
\partial_{q^{k}}\Pi_{0}$ and
\eqref{eq.deunitOD}, the last factor is
\begin{eqnarray}
  \nonumber
&&\hspace{-1cm}
\left[u_{0}\left\{u_{0}^{*},u_{1}\right\}u_{0}^{*}
+u_{0}\left\{u_{1}^{*},u_{0}\right\}u_{0}^{*}\right]^{D}
=
-(u_{0}\partial_{q^{\ell}}u_{0}^{*})^{OD}
\frac{(\partial^{2}_{p_{\ell}p_{k}}f_{\varepsilon})}{i(E_{+}-E_{-})}(\partial_{q^{k}}\Pi_{0})
\\
\nonumber
&&\hspace{4cm}
-
\frac{(\partial^{2}_{p_{\ell}p_{k}}f_{\varepsilon})}{i(E_{+}-E_{-})}(\partial_{q^{k}}\Pi_{0})((\partial_{q^{\ell}}u_{0})u_{0}^{*})^{OD}
\\
\label{eq.poissdiag}
&&\hspace{3cm}=-2\frac{(\partial^{2}_{p_{\ell}p_{k}}f_{\varepsilon})}{i(E_{+}-E_{-})}(\partial_{q^{k}}\Pi_{0})(\partial_{q^{\ell}}\Pi_{0})\sigma_{3}\,.
\end{eqnarray}
With
$u_{1}u_{0}^{*}=\frac{(\partial_{p_{k}}f_{\varepsilon})}{i(E_{+}-E_{-})}(\partial_{q^{k}}\Pi_{0})$\,,
we have proved
\begin{multline}
  \label{eq.21}
\left[u_{0}u_{2}^{*}H+H
  u_{2}u_{0}^{*}+u_{0}u_{1}^{*}Hu_{1}u_{0}^{*}\right]^{D}
=-\frac{(\partial_{p_{k}}f_{\varepsilon})(\partial_{p_{\ell}}f_{\varepsilon})}{(E_{+}-E_{-})}
(\partial_{q^{k}}\Pi_{0})(\partial_{q^{\ell}}\Pi_{0})\sigma_{3}
\\
-H\frac{(\partial^{2}_{p_{\ell}p_{k}}f_{\varepsilon})}{E_{+}-E_{-}}
(\partial_{q^{k}}\Pi_{0})(\partial_{q^{\ell}}\Pi_{0})\sigma_{3}
+ \varepsilon^{1+2\delta} \underline{\tilde{R}_{1}}+\varepsilon^{1+4\delta}\tilde{R}_{1}\,.
\end{multline}
\item  The term
  $$
\left[\frac{1}{2i}u_{0}\left\{u_{0}^{*}f_{\varepsilon},u_{1}\right\}u_{0}^{*}+
    \frac{1}{2i}u_{0}\left\{u_{1}^{*}f_{\varepsilon},u_{0}\right\}u_{0}^{*}-\frac{(\partial_{p_{k}}f_{\varepsilon})}{2i}
u_{0}(\partial_{q^{k}}u_{0}^{*})u_{1}u_{0}^{*}\right]^{D}
$$
equals
  \begin{eqnarray*}
    &&
\frac{1}{2i}f_{\varepsilon}\left[u_{0}\left\{u_{0}^{*},u_{1}\right\}u_{0}^{*}+u_{0}\left\{u_{1}^{*},u_{0}\right\}u_{0}^{*}\right]^{D}
\\
&&
+\frac{(\partial_{p_{k}}f_{\varepsilon})}{2i}
\left[(\partial_{q^{k}}u_{1})u_{0}^{*}+u_{0}u_{1}^{*}(\partial_{q^{k}}u_{0})u_{0}^{*}
-u_{0}(\partial_{q^{k}}u_{0}^{*})u_{1}u_{0}^{*}\right]^{D}\,.
  \end{eqnarray*}
The diagonal part of $(\partial_{q^{k}}u_{1})u_{0}^{*}$ is
$$
\left[(\partial_{q^{k}}u_{1})u_{0}^{*}\right]^{D}=\frac{\partial_{p_{\ell}}f_{\varepsilon}}{i(E_{+}-E_{-})}
(\partial^{2}_{q^{k}q^{\ell}}\Pi_{0})^{D}\,.
$$
But differentiating the relation
$(\partial_{q^{k}}\Pi_{0})\Pi_{0}+\Pi_{0}(\partial_{q^{\ell}}\Pi_{0})=\partial_{q^{k}}\Pi_{0}$
w.r.t $q^{\ell}$
leads to
\begin{equation}
  \label{eq.dersecdiag}
(\partial^{2}_{q^{k}q^{\ell}}\Pi_{0})^{D}=-(\partial_{q^{k}}\Pi_{0}\partial_{q^{\ell}}\Pi_{0}
+ \partial_{q^{\ell}}\Pi_{0}\partial_{q^{k}}\Pi_{0})\sigma_{3}\,.
\end{equation}
With
\eqref{eq.deunitOD} and $(u_{1}u_{0}^{*})=(u_{1}u_{0}^{*})^{OD}=\frac{\partial_{p_{\ell}}f_{\varepsilon}}{i(E_{+}-E-)}\partial_{q^{\ell}}\Pi_{0}$, we obtain
\begin{eqnarray*}
&&\hspace{-1cm}\left[u_{0}u_{1}^{*}(\partial_{q^{k}}u_{0})u_{0}^{*}
-u_{0}(\partial_{q^{k}}u_{0}^{*})u_{1}u_{0}^{*}\right]^{D}
=
[-u_{0}u_{1}^{*}\sigma_{3}(\partial_{q^{k}}\Pi_{0})+(\partial_{q^{k}}\Pi_{0})\sigma_{3}u_{1}u_{0}^{*}]
\\
&&
\hspace{2cm}=
-\frac{(\partial_{p_{\ell}}f_{\varepsilon})}{i(E_{+}-E_{-})}
\left[(\partial_{q^{\ell}}\Pi_{0})(\partial_{q^{k}}\Pi_{0})
+(\partial_{q^{k}}\Pi_{0})(\partial_{q^{\ell}}\Pi_{0})\right]\sigma_{3}\,.
\end{eqnarray*}
We have found
\begin{multline}
\label{eq.22}
\left[\frac{1}{2i}u_{0}\left\{u_{0}^{*}f_{\varepsilon},u_{1}\right\}u_{0}^{*}+
    \frac{1}{2i}u_{0}\left\{u_{1}^{*}f_{\varepsilon},u_{0}\right\}u_{0}^{*}-\frac{(\partial_{p_{k}}f_{\varepsilon})}{2i}
u_{0}(\partial_{q^{k}}u_{0}^{*})u_{1}u_{0}^{*}\right]^{D}\\
=
\frac{f_{\varepsilon}(\partial^{2}_{p_{k}p_{\ell}}f_{\varepsilon})}{E_{+}-E_{-}}
(\partial_{q^{k}}\Pi_{0})(\partial_{q^{\ell}}\Pi_{0})\sigma_{3}
+
2\frac{(\partial_{p_{k}}f_{\varepsilon})(\partial_{p_{\ell}}f_{\varepsilon})}{E_{+}-E_{-}}
(\partial_{q^{k}}\Pi_{0})(\partial_{q^{\ell}}\Pi_{0})\sigma_{3}\,.
\end{multline}
\item The diagonal part of
  $\frac{1}{2i}u_{0}\left[\left\{u_{0}^{*}\mathcal{V},u_{1}\right\}
+\left\{u_{1}^{*}\mathcal{V},u_{0}\right\}\right]u_{0}^{*}$ equals
$$
\frac{1}{2i}\left[-u_{0}(\partial_{q^{k}}u_{0}^{*})\mathcal{V}(\partial_{p_{k}}u_{1})u_{0}^{*}
-(\partial_{q^{k}}\mathcal{V})(\partial_{p_{k}}u_{1})u_{0}^{*} + u_{0}(\partial_{p_{k}}u_{1}^{*})\mathcal{V}(\partial_{q^{k}}u_{0})u_{0}^{*}\right]^{D}\,.
$$
By using \eqref{eq.deunitOD} with
\begin{eqnarray*}
&&(\partial_{p_{k}}u_{1})u_{0}^{*}=\left[(\partial_{p_{k}}u_{1})u_{0}^{*}\right]^{OD}
=\frac{\partial^{2}_{p_{k}p_{\ell}}f_{\varepsilon}}{i(E_{+}-E_{-})}(\partial_{q^{\ell}}\Pi_{0})\,,
\\
&&
(\partial_{q}\mathcal{V})^{OD}=(E_{+}-E_{-})(\partial_{q}\Pi_{0})\,,
\\
\text{and}&&
E_{-}\Pi_{0}+E_{+}(1-\Pi_{0})=\mathcal{V}+(E_{-}-E_{+})\sigma_{3}\,,
\end{eqnarray*}
it becomes
\begin{multline*}
\frac{\partial^{2}_{p_{k}p_{\ell}}f_{\varepsilon}}{2}(\partial_{q^{k}}\Pi_{0})(\partial_{q^{\ell}}\Pi_{0})
+
\frac{1}{2i}\mathcal{V}\left[u_{0}\left\{u_{0}^{*},u_{1}\right\}u_{0}^{*}
+u_{0}\left\{u_{1}^{*},u_{0}\right\}u_{0}^{*}\right]^{D}
\\
+\frac{1}{2i}
(E_{-}-E_{+})\left[u_{0}\left\{u_{0}^{*},u_{1}\right\}u_{0}^{*}
+u_{0}
\left\{u_{1}^{*},u_{0}\right\}u_{0}^{*}\right]^{D}\sigma_{3}
\,.
\end{multline*}
The relation \eqref{eq.poissdiag} yields
\begin{multline}
  \label{eq.23}
\frac{1}{2i}\left[u_{0}\left\{u_{0}^{*}\mathcal{V},u_{1}\right\}u_{0}^{*}
+u_{0}\left\{u_{1}^{*}\mathcal{V},u_{0}\right\}u_{0}^{*}\right]
=
-\frac{\partial^{2}_{p_{k}p_{\ell}}f_{\varepsilon}}{2}(\partial_{q^{k}}\Pi_{0})(\partial_{q^{\ell}}\Pi_{0})
\\
+\mathcal{V}\frac{\partial^{2}_{p_{k}p_{\ell}}f_{\varepsilon}}{E_{+}-E_{-}}(\partial_{q^{k}}\Pi_{0})(\partial_{q^{\ell}}\Pi_{0})\,.
\end{multline}
\item The term
  $\frac{1}{2i}\left[u_{0}\left\{u_{1}^{*},H\right\}\right]^{D}$
  is the sum of two terms
$$
\frac{1}{2i}\left[u_{0}\left\{u_{1}^{*},f_{\varepsilon}\right\}\right]^{D}
+
\frac{1}{2i}\left[u_{0}\left\{u_{1}^{*},\mathcal{V}\right\}\right]^{D}\,.
$$
Since
$u_{0}\partial_{p_{\ell}}u_{1}^{*}=i\frac{\partial^{2}_{p_{k}p_{\ell}}f_{\varepsilon}}{(E_{+}-E_{-})}(\partial_{q^{k}}\Pi_{0})$
is off-diagonal while
$$
(\partial_{q^{\ell}}\mathcal{V})^{OD}=(E_{+}-E_{-})(\partial_{q^{\ell}}\Pi_{0})\,,
$$
the second term equals
$$
\frac{\partial^{2}_{p_{k}p_{\ell}}f_{\varepsilon}}{2}
(\partial_{q^{k}}\Pi_{0})(\partial_{q^{\ell}}\Pi_{0})\,.
$$
The first term  a priori contains more terms because $u_{0}^{*}$ has to be differentiated:
\begin{multline*}
\frac{1}{2}
\Big[
-(\partial_{p_{k}}f_{\varepsilon})(\partial_{p_{\ell}}f_{\varepsilon})\partial_{q^{\ell}}\left(\frac{1}{E_{+}-E_{-}}\right)(\partial_{q^{k}}\Pi_{0})
-
\frac{(\partial_{p_{k}}f_{\varepsilon})(\partial_{p_{\ell}}f_{\varepsilon})}{E_{+}-E_{-}}
(\partial^{2}_{q^{\ell}q^{k}}\Pi_{0})
\\
-u_{0}(\partial_{q^{\ell}}u_{0}^{*}) \frac{(\partial_{p_{k}}f_{\varepsilon})(\partial_{p_{\ell}}f_{\varepsilon})}{E_{+}-E_{-}}
(\partial_{q^{k}}\Pi_{0})\Big]^{D}\,.
\end{multline*}
The first part is off-diagonal and vanishes after taking the diagonal
part.
By using again \eqref{eq.deunitOD} and \eqref{eq.dersecdiag} and $\partial_{q}\Pi_{0}=(\partial_{q}\Pi_{0})^{OD}$, we
obtain
$$
\frac{(\partial_{p_{k}}f_{\varepsilon})
(\partial_{p_{\ell}}f_{\varepsilon})}{E_{+}-E_{-}}(\partial_{q^{k}}\Pi_{0})
(\partial_{q^{\ell}}\Pi_{0})\sigma_{3}
+
\frac{(\partial_{p_{k}}f_{\varepsilon})
(\partial_{p_{\ell}}f_{\varepsilon})}{E_{+}-E_{-}}(\partial_{q^{k}}\Pi_{0})\sigma_{3}
(\partial_{q^{\ell}}\Pi_{0})=0\,.
$$
Hence we have proved
\begin{equation}
  \label{eq.23bis}
  \frac{1}{2i}\left[u_{0}\left\{u_{1}^{*},H\right\}\right]^{D}
=\frac{\partial^{2}_{p_{k}p_{\ell}}f_{\varepsilon}}{2}
(\partial_{q^{k}}\Pi_{0})(\partial_{q^{\ell}}\Pi_{0})\,.
\end{equation}
\item The last term is
  \begin{multline*}
T=\frac{1}{8}\Big[2u_{0}\left\{(\partial_{p_{k}}f_{\varepsilon})(\partial_{q^{k}}u_{0}^{*}),
    u_{0}\right\}u_{0}^{*}\\
-u_{0}(\partial^{2}_{q^{k}q^{\ell}}u_{0}^{*})(\partial^{2}_{p_{k}p_{\ell}}f_{\varepsilon})
-
(\partial^{2}_{p_{k}p_{\ell}}f_{\varepsilon})(\partial^{2}_{q^{k}q^{\ell}}u_{0})u_{0}^{*}\Big]^{D}\,.
\end{multline*}
Forgetting the $\frac{1}{8}$ factor, the first part equals
$2(\partial^{2}_{p_{k}p_{\ell}}f_{\varepsilon})u_{0}(\partial_{q^{k}}u_{0}^{*})
(\partial_{q^{\ell}}u_{0})u_{0}^{*}$\,. By using \eqref{eq.der2u0u0}
in the second part gives
$$
8T=4(\partial_{p_{k}p_{\ell}}f_{\varepsilon})
\left[u_{0}(\partial_{q^{k}}u_{0}^{*})(\partial_{q^{\ell}}u_{0})u_{0}^{*}\right]^{D}\,.
$$
Owing to \eqref{eq.diagprod} and \eqref{eq.deunitOD}, this gives
$$
2T=
(\partial^{2}_{p_{k}p_{\ell}}f_{\varepsilon})
[u_{0}(\partial_{q^{k}}u_{0}^{*})]^{D}[(\partial_{q^{\ell}}u_{0})u_{0}^{*}]^{D}
+
(\partial^{2}_{p_{k}p_{\ell}}f_{\varepsilon})(\partial_{q^{k}}\Pi_{0})(\partial_{q^{\ell}}\Pi_{0})\,.
$$
With $u_{0}(\partial_{q^{k}}u_{0}^{*})=-(\partial_{q^{k}}u_{0})u_{0}^{*}$, the last term equals
\begin{equation}
  \label{eq.24}
  T=
  -\frac{(\partial^{2}_{p_{k}p_{\ell}}f_{\varepsilon})}{2}[(\partial_{q^{k}}u_{0})u_{0}^{*}]^{D}
[(\partial_{q^{\ell}}u_{0})u_{0}^{*}]^{D}
+\frac{(\partial^{2}_{p_{k}p_{\ell}}f_{\varepsilon})}{2}
(\partial_{q^{k}}\Pi_{0})(\partial_{q^{\ell}}\Pi_{0})\,.
\end{equation}
\end{itemize}
By summing \eqref{eq.21},\eqref{eq.22},\eqref{eq.23},\eqref{eq.23bis},\eqref{eq.24}, we
obtain
\begin{multline*}
B_{2}^{D}=  -\frac{(\partial^{2}_{p_{k}p_{\ell}}f_{\varepsilon})}{2}[(\partial_{q^{k}}u_{0})u_{0}^{*}]^{D}
[(\partial_{q^{\ell}}u_{0})u_{0}^{*}]^{D}
+ \frac{(\partial_{p_{k}}f_{\varepsilon})(\partial_{p_{\ell}}f_{\varepsilon})}{E_{+}-E_{-}}
(\partial_{q^{k}}\Pi_{0})(\partial_{q^{\ell}}\Pi_{0})\sigma_{3}
\\
+\frac{\partial^{2}_{p_{k}p_{\ell}}f_{\varepsilon}}{2}
(\partial_{q^{k}}\Pi_{0})(\partial_{q^{\ell}}\Pi_{0})
+\varepsilon^{1+2\delta}\underline{\tilde{R}_{1}}+ \varepsilon^{1+4\delta}\tilde{R}_{1}
\,.
\end{multline*}
Hence the diagonal symbol that we seek, is
\begin{multline*}
\begin{pmatrix}
  h_{+} &0\\
0&h_{-}
\end{pmatrix}
=
\begin{pmatrix}
  f_{\varepsilon}(p,\tau)+E_{+}(q,\tau) & 0\\
0& f_{\varepsilon}(p,\tau)+E_{-}(q,\tau)
\end{pmatrix}
\\-
i\varepsilon
(\partial_{p_{k}}f_{\varepsilon})u_{0}^{*}\left[(\partial_{q^{k}}u_{0})u_{0}^{*}\right]^{D}u_{0}
-\varepsilon^{2}\frac{\partial^{2}_{p_{k}p_{\ell}}f_{\varepsilon}}{2}u_{0}^{*}
\left[(\partial_{q^{k}}u_{0})u_{0}^{*}\right]^{D}\left[(\partial_{q^{k}}u_{0})u_{0}^{*}\right]^{D}u_{0}
\\
+\varepsilon^{2}
\frac{(\partial_{p_{k}}f_{\varepsilon})(\partial_{p_{\ell}}f_{\varepsilon})}{E_{+}-E_{-}}
u_{0}^{*}(\partial_{q^{k}}\Pi_{0})(\partial_{q^{\ell}}\Pi_{0})\sigma_{3}u_{0}
+
\varepsilon^{2}\frac{\partial^{2}_{p_{k}p_{\ell}}f_{\varepsilon}}{2}u_{0}^{*}(\partial_{q^{k}}\Pi_{0})(\partial_{q^{\ell}}\Pi_{0})u_{0}
\,.
\end{multline*}
The symbol
$u_{0}^{*}\left[(\partial_{q^{k}}u_{0})u_{0}^{*}\right]^{D}u_{0}$ equals
\begin{multline*}
  u_{0}^{*}\Pi_{0}(\partial_{q^{k}}u_{0})u_{0}^{*}\Pi_{0}u_{0}+
 u_{0}^{*}(1-\Pi_{0})(\partial_{q^{k}}u_{0})u_{0}^{*}(1-\Pi_{0})u_{0}
\\
=
P_{+}u_{0}^{*}(\partial_{q^{k}}u_{0})P_{+}
+
(1-P_{+})u_{0}^{*}(\partial_{q^{k}}u_{0})(1-P_{+})
=-i
\begin{pmatrix}
  A_{k}&0\\
0& -A_{k}
\end{pmatrix}\,.
\end{multline*}
For the last term we deduce from $\Pi_{0}=u_{0}P_{+}u_{0}^{*}$ and
$(\partial_{q}u_{0})u_{0}^{*}+ u_{0}(\partial_{q}u_{0})=0$,
\begin{eqnarray*}
&&
P_{+}\left[u_{0}^{*}(\partial_{q^{k}}\Pi_{0})u_{0}\right]
\left[u_{0}^{*}(\partial_{q^{\ell}}\Pi_{0})u_{0}\right]P_{+}
\\
&&\qquad
=
P_{+}\left[u_{0}^{*}(\partial_{q^{k}}u_{0})P_{+}+
  P_{+}(\partial_{q^{k}}u_{0}^{*})u_{0}\right]
\left[u_{0}^{*}(\partial_{q^{\ell}}u_{0})P_{+}+
  P_{+}(\partial_{q^{\ell}}u_{0}^{*})u_{0}\right]
P_{+}
\\
&&\qquad
=
P_{+}u_{0}^{*}(\partial_{q^{k}}u_{0})P_{+}u_{0}^{*}(\partial_{q^{\ell}}u_{0})P_{+}
+
P_{+}u_{0}^{*}(\partial_{q^{k}}u_{0})P_{+}(\partial_{q^{\ell}}u_{0}^{*})u_{0}P_{+}
\\
&&\qquad\qquad
+
P_{+}(\partial_{q^{k}}u_{0}^{*})u_{0}u_{0}^{*}(\partial_{q^{\ell}}u_{0})P_{+}
+
P_{+}(\partial_{q^{k}}u_{0}^{*})u_{0}P_{+}(\partial_{q^{\ell}}u_{0}^{*})u_{0}P_{+}
\\
&&\qquad
=
-P_{+}u_{0}^{*}(\partial_{q^{k}}u_{0})(1-P_{+})u_{0}^{*}(\partial_{q^{\ell}}u_{0})P_{+}
= X_{k}\overline{X}_{\ell}\,.
\end{eqnarray*}
Taking the  bracket with $(1-P_{+})$ is even simpler and gives
\begin{multline*}
(1-P_{+})\left[u_{0}^{*}(\partial_{q^{k}}\Pi_{0})u_{0}\right]
\left[u_{0}^{*}(\partial_{q^{\ell}}\Pi_{0})u_{0}\right](1-P_{+})
\\
=
-(1-P_{+})u_{0}^{*}(\partial_{q^{k}}u_{0})P_{+}u_{0}^{*}(\partial_{q^{\ell}}u_{0})(1-P_{+})
=
\overline{X_{k}}X_{\ell}\,.
\end{multline*}
This ends the proof.
\enddemonstration{}

\subsection{Discussion about the adiabatic approximation of the Born-Oppenheimer Hamiltonian}
\label{se.discuss}
The Theorem~\ref{th.BornOp} is a direct application of
Proposition~\ref{pr.hameff} by taking
$$
f_{\varepsilon}(p)= \varepsilon^{2\delta}\tau'\tau''|p|^{2}\gamma(\tau'\tau''|p|^{2})\,.
$$
The operators
\begin{eqnarray*}
h_{\pm, BO}(q,\varepsilon D_{q},\varepsilon)&=&
\varepsilon^{2\delta}\tau'\tau''\left[\sum_{k=1}^{d}(p_{k}
  \mp \varepsilon A_{k})^{2}\right]^{Weyl}+\varepsilon^{2+2\delta}\tau'\tau''\sum_{k=1}^{d}|X_{k}|^{2}\\
&=&
\varepsilon^{2+2\delta}\tau'\tau''\left[
-\Delta
+ |A|^{2}
\mp \left[(\frac{1}{i}\nabla).A
+ A.(\frac{1}{i}\nabla)\right]
+|X|^{2}\right]
\\
&=&
\varepsilon^{2+2\delta}\tau'\tau''
\left[
|\frac{1}{i}\nabla \mp A|^{2}+|X|^{2}
\right]
\,,
\end{eqnarray*}
is nothing but
the usual adiabatic effective Hamiltonians which can be found in the
physics literature, including the Born-Huang potential
$|X|^{2}=\sum_{k=1}^{d}|X_{k}|^{2}$\,.
We refer to \cite{PST} and \cite{PST2} for a discussion of the various
presentations of the calculations and additional references.\\
Even in the region $\left\{\sqrt{\tau'\tau''}|p|\leq
  r_{\gamma}\right\}$ with $p$ quantized into $\varepsilon D_{q}$,
this approximation makes sense, only for
$\delta>0$, because of the additional term
$$
\mp
2\frac{\varepsilon^{2+4\delta}\tau'\tau''}{E_{+}-E_{-}}\bar{X_{k}}X_{\ell}
(\varepsilon \partial_{q_{k}})(\varepsilon \partial_{q_{\ell}})
$$
coming from the last terms of \eqref{eq.h+} and \eqref{eq.h-}, that we
have included in the remainder. It is
not surprising (see \cite{Sor}, \cite{MaSo}) that the degree of the
differential operators increases with the degree in $\varepsilon$ in
the adiabatic expansion of Schr{\"o}dinger type Hamiltonians. The argument
of physicists says that this effective Hamiltonian is used
for  relatively small frequencies (or momentum) so that $X_{k}X_{\ell}p^{k}p^{\ell}$ is negligible
w.r.t $|p-\varepsilon A|^{2}+\varepsilon^{2}|X|^{2}$\,. The introduction of
the additional factor $\varepsilon^{2\delta}$ with $\delta >0$
provides
 a mathematically accurate and rather flexible implementation of this
 approximation.

\section{Adaptation of the adiabatic asymptotics to the full nonlinear
minimization problem}
\label{se.minipb}
In this section, we adapt our rather general adiabatic result to our
 nonlinear problem. In a first step, we give an explicit
 form of Theorem~\ref{th.BornOp} in our specific framework. Those
results are effective when applied with wave functions localized in
the frequency variable, $\psi=\chi(\sqrt{\tau_{x}\tau_{y}}\varepsilon
D_{q})\psi$ for some compactly supported $\chi$\,. It could suffice if
we considered  minimizing the energy among such well prepared quantum
states. We can do better by using a partition of unity in the
frequency variable, which will be combined, in the end, with the a priori
estimates coming from the complete and reduced minimization problems.
Finally an estimate of the effect of the unitary transform $\hat{U}$
on the nonlinear term is provided.
\subsection{Adiabatic approximation for the explicit Schr{\"o}dinger
  Hamiltonian}
\label{se.explred}
Let us specify the result of Theorem~\ref{th.BornOp} by going back to the coordinates
$q=(q',q'')=(x,y)$ and $\tau=(\tau',\tau'')=(\tau_{x},\tau_{y})$\,.
 Provided that  $V_{\varepsilon,\tau}(x,y)$
 fulfills the proper  assumptions, the unitary transform introduced in
 Theorem~\ref{th.BornOp}  transforms
 $H_{Lin}$ given by (\ref{Hlin}) into the Born-Oppenheimer Hamiltonian
with a good accuracy
in the low frequency region, that is when applied to wave
 functions $\psi$ such that
 $\psi=\chi(\sqrt{\tau_{x}\tau_{y}}\varepsilon D_{q})\psi$ for some
 compactly supported $\chi$\,.

The operator $H_{Lin}$ is the $\varepsilon$-quantization of the
symbol
\begin{eqnarray*}
&&\varepsilon^{2\delta}\tau_{x}\tau_{y}|p|^{2} + u_{0}(q,\tau)
\begin{pmatrix}
  E_{+}(q,\tau,\varepsilon)& 0\\
0&E_{-}(q,\tau,\varepsilon)
\end{pmatrix}u_{0}(q,\tau)^{*}\,,\\
\text{with}&&E_{\pm}(q,\tau,\varepsilon)= V_{\varepsilon,\tau}(x,y)\pm
\Omega(\sqrt{\frac{\tau_{x}}{\tau_{y}}}x)=V_{\varepsilon,\tau}(x,y)\pm
\sqrt{1+\frac{\tau_{x}}{\tau_{y}}x^{2}}\,.
\end{eqnarray*}
The operator $u_{0}(q,\tau)$ equals
\begin{eqnarray*}
&&u_{0}(q,\tau)=
\begin{pmatrix}
  C& Se^{i\varphi}\\
 Se^{-i\varphi}& -C
\end{pmatrix}=u_{0}(q,\tau)^{*}\\
\text{with}&&C=\cos(\frac{\theta}{2})\,,\quad S=\sin(\frac{\theta}{2})\,.
\\
\text{and}&&\cot(\theta)=\sqrt{\frac{\tau_{x}}{\tau_{y}}}x\quad,\quad
\varphi=\sqrt{\frac{\tau_{y}}{\tau_{x}}}y\,.
\end{eqnarray*}
 \begin{proposition}
\label{pr.adiabexpl1}
   Assume $\delta \in (0,\delta_{0}]$, $(\tau_{x},\tau_{y})\in (0,1]^{2}$
and assume that $V_{\varepsilon,\tau}$ belongs to $S_{u}(\langle
   \sqrt{\frac{\tau_{x}}{\tau_{y}}}x\rangle,\frac{\frac{\tau_{x}}{\tau_{y}}dx^{2}}{\langle
     \sqrt{\frac{\tau_{x}}{\tau_{y}}}
x\rangle^{2}}+ \frac{\tau_{y}}{\tau_{x}}dy^{2})$ then the
matricial
   potential
$$
\mathcal{V}(q,\tau,\varepsilon)=V_{\varepsilon,\tau}(x,y)+\Omega( \sqrt{\frac{\tau_{x}}{\tau_{y}}}x)
   \begin{pmatrix}
     \cos(\theta)&e^{i\varphi}\sin (\theta)\\
e^{-i\varphi}\sin(\theta)& -\cos(\theta)
   \end{pmatrix}
$$
fulfills the assumption of Theorem~\ref{th.BornOp}. \\
Choose the function $\gamma\in \mathcal{C}^{\infty}_{0}(\rz)$ such that  $\gamma\equiv 1$ in a
neighborhood of $[-r_{\gamma}^{2},r_{\gamma}^{2}]$,
as in Theorem~\ref{th.BornOp}
and consider a cut-off function $\chi\in
\mathcal{C}_{0}^{\infty}((-r_{\gamma}^{2},r_{\gamma}^{2}))$. When
$\hat{U}=U(q,\varepsilon D_{q},\tau,\varepsilon)\in
OpS_{u}(1,g_{\tau}; \mathcal{M}_{2}(\cz))$
 is given in Theorem~\ref{th.BornOp}, the identities
\begin{eqnarray}
\nonumber
&&\hat{U}^{*}H_{Lin}\hat{U}\chi(\tau_{x}\tau_{y}|\varepsilon D_{q}|^{2})
-
\varepsilon^{2+4\delta}R_{1}(\tau,\varepsilon)=\\
\label{eq.U*HUchi}
&&
\hspace{-1.5cm}
\varepsilon^{2+2\delta}\tau_{x}\tau_{y}
{\small\begin{pmatrix}
 e^{+i\sqrt{\frac{\tau_{y}}{\tau_{x}}}\frac{y}{2}}\hat{H}_{+}e^{-i\sqrt{\frac{\tau_{y}}{\tau_{x}}}\frac{y}{2}}&0
\\
0&e^{-i\sqrt{\frac{\tau_{y}}{\tau_{x}}}\frac{y}{2}}\hat{H}_{-}e^{i\sqrt{\frac{\tau_{y}}{\tau_{x}}}\frac{y}{2}}
\end{pmatrix}
}\chi(\tau_{x}\tau_{y}|\varepsilon D_{q}|^{2})
\,,\\
\nonumber
 &&H_{Lin}\chi(\tau_{x}\tau_{y}|\varepsilon D_{q}|^{2})-\varepsilon^{2+4\delta}R_{2}(\tau,\varepsilon)=
\\
\label{eq.UHU*chi}
&&\hspace{-1.5cm}
 \varepsilon^{2+2\delta}\tau_{x}\tau_{y} \hat{U}
{\small\begin{pmatrix}
 e^{+i\sqrt{\frac{\tau_{y}}{\tau_{x}}}\frac{y}{2}}\hat{H}_{+}e^{-i\sqrt{\frac{\tau_{y}}{\tau_{x}}}\frac{y}{2}}&0\\
0&e^{-i\sqrt{\frac{\tau_{y}}{\tau_{x}}}\frac{y}{2}}\hat{H}_{-}e^{i\sqrt{\frac{\tau_{y}}{\tau_{x}}}\frac{y}{2}}
\end{pmatrix}
}\hat{U}^{*}\chi(\tau_{x}\tau_{y}|\varepsilon D_{q}|^{2})
\,,
\end{eqnarray}
hold
with
\begin{eqnarray*}
  \hat{H}_{+}&=& -\partial_{x}^{2}-
\left(\partial_{y}+ i\frac{x}{2\sqrt{1+\frac{\tau_{x}}{\tau_{y}}x^{2}}}\right)^{2}
\\
&&+ \frac{1}{\varepsilon^{2+2\delta}\tau_{x}\tau_{y}}\left[V_{\varepsilon,\tau}(x,y)+\sqrt{1+\frac{\tau_{x}}{\tau_{y}}x^{2}}\right]+W_{\tau}(x,y)
\,,\\
 \hat{H}_{-}&=& -\partial_{x}^{2}-
\left(\partial_{y}-
i\frac{x}{2\sqrt{1+\frac{\tau_{x}}{\tau_{y}}x^{2}}}\right)^{2}
\\
&&+ \frac{1}{\varepsilon^{2+2\delta}\tau_{x}\tau_{y}}\left[V_{\varepsilon,\tau}(x,y)-\sqrt{1+\frac{\tau_{x}}{\tau_{y}}x^{2}}\right]
+W_{\tau}(x,y)\,,
\\
W_{\tau}(x,y)&=&
\frac{\tau_{x}}{\tau_{y}(1+\frac{\tau_{x}}{\tau_{y}}x^{2})^{2}}+ \frac{\tau_{y}}{\tau_{x}(1+\frac{\tau_{x}}{\tau_{y}}x^{2})}\,,
\end{eqnarray*}
and the estimates
$$
\|R_{1}(\tau,\varepsilon)\|_{\mathcal{L}(L^{2})}+\|R_{2}(\tau,\varepsilon)\|_{\mathcal{L}(L^{2})}\leq
C
$$
which are uniform w.r.t $\delta\in (0,\delta_{0}]$, $\tau\in
(0,1]^{2}$ and $\varepsilon\in
(0,\varepsilon_{0}]$\,.
 \end{proposition}
\begindemonstration{}
Following the approach of Section~\ref{se.BO}, the Hamiltonian
$H_{Lin}$ is decomposed into
\begin{eqnarray*}
  &&H_{Lin}= H(q,\varepsilon D_{q},\tau,\varepsilon)+
  R_{\gamma}(\varepsilon D_{q},\tau,\varepsilon)\\
\text{with}
&&H(q,p,\varepsilon)=\varepsilon^{2\delta}\tau_{x}\tau_{y}|p|^{2}\gamma(\tau_{x}\tau_{y}|p|^{2})
+\Big[u_{0}
\begin{pmatrix}
  E_{+}& 0\\
0&E_{-}
\end{pmatrix}u_{0}^{*}\Big](q,\tau)\,,\\
\text{and}&&R_{\gamma}(p,\tau,\varepsilon)=\varepsilon^{2\delta}\tau_{x}\tau_{y}|p|^{2}(1-\gamma(\tau_{x}\tau_{y}|p|^{2}))\,.
\end{eqnarray*}
The Hamiltonian $H(q,\varepsilon D_{q},\tau,\varepsilon)$ fulfills the
assumptions of Theorem~\ref{th.BornOp}, with the metric
$g_{\tau}=\frac{\frac{\tau_{x}}{\tau_{y}}dx^{2}}{\langle\sqrt{\frac{\tau_{x}}{\tau_{y}}}
  x\rangle^{2}}+
\frac{\tau_{y}}{\tau_{x}}dy^{2}+\frac{\tau_{x}\tau_{y}dp^{2}}{\langle
  \sqrt{\tau_{x}\tau_{y}}p\rangle^{2}}$,   because we assumed
 $$
V_{\varepsilon,\tau}\in S_{u}(\langle
\sqrt{\frac{\tau_{x}}{\tau_{y}}}x\rangle, \frac{\frac{\tau_{x}}{\tau_{y}}dx^{2}}{\langle\sqrt{\frac{\tau_{x}}{\tau_{y}}} x\rangle^{2}}+\frac{\tau_{y}}{\tau_{x}}dy^{2}
)
$$
while the gap $E_{+}(q,\tau,\varepsilon)-E_{-}(q,\tau,\varepsilon)$ equals
$2\sqrt{1+\frac{\tau_{x}}{\tau_{y}}x^{2}}$\,.
Actually the estimates
$$
|\partial_{x}^{\alpha}\partial_{y}^{\beta}u_{0}(x,y)|\leq
C_{\alpha,\beta}\left\langle \sqrt{\frac{\tau_{x}}{\tau_{y}}}
x\right\rangle^{-|\alpha|}\left(\frac{\tau_{x}}{\tau_{y}}\right)^{\frac{|\alpha|-|\beta|}{2}}
$$
are due to
$\partial_{x}\theta=-\frac{\sqrt{\tau_{x}}}{\sqrt{\tau_{y}}(1+\frac{\tau_{x}}{\tau_{y}}x^{2})}$ and
$\partial_{y}e^{i\varphi}=i\sqrt{\frac{\tau_{y}}{\tau_{x}}}e^{i\varphi}$\,.
Moreover the explicit computation with $\alpha+ \beta=1$ leads to
\begin{eqnarray}
 \label{eq.Dxu0} &&
  \begin{pmatrix}
    A_{x}& X_{x}\\
\overline{X_{x}}& -A_{x}
  \end{pmatrix}= iu_{0}^{*}(\partial_{x}u_{0})=
    \frac{i\partial_{x}\theta}{2}
    \begin{pmatrix}
      0 & e^{i\varphi}\\
     -e^{-i\varphi}&0
    \end{pmatrix}
\\
\label{eq.Dyu0}
\text{and}
&&
 \begin{pmatrix}
    A_{y}& X_{y}\\
\overline{X_{y}}& -A_{y}
  \end{pmatrix}= iu_{0}^{*}(\partial_{y}u_{0})=
\sqrt{\frac{\tau_{y}}{\tau_{x}}} S
  \begin{pmatrix}
  S  & -Ce^{i\varphi}\\
-Ce^{-i\varphi} & -S
  \end{pmatrix}\,.
\end{eqnarray}
Hence the effective Hamiltonians $h_{B0,\pm}(q,\varepsilon
D_{q},\varepsilon)$, when restricted to the region
$\left\{\sqrt{\tau_{x}\tau_{y}}|p|\leq r_{\gamma}\right\}$, are given by the symbols
\begin{multline*}
  h_{BO,\pm}=\varepsilon^{2\delta}\tau_{x}\tau_{y}\left[p_{x}^{2} +
  (p_{y}\mp \varepsilon
\sqrt{\frac{\tau_{y}}{\tau_{x}}}
\sin^{2}(\frac{\theta}{2}))^{2}
+ \frac{1}{4}(|\partial_{x}\theta|^{2}+
\frac{\tau_{y}}{\tau_{x}}\sin^{2}(\theta))
\right]\\
+ V_{\varepsilon,\tau}(x,y)\pm \sqrt{1+\frac{\tau_{x}}{\tau_{y}}x^{2}}\,.
\end{multline*}
With
$\cos(\theta)=\frac{\sqrt{\tau_{x}}x}{\sqrt{\tau_{y}}\sqrt{1+\frac{\tau_{x}}{\tau_{y}}x^{2}}}$
and
$\sin(\theta)=\frac{1}{\sqrt{1+\frac{\tau_{x}}{\tau_{y}}x^{2}}}$,
the Schr{\"o}dinger Hamiltonian corresponding to the RHS is
$\varepsilon^{2+2\delta}\tau_{x}\tau_{y}e^{\pm
  i\frac{\sqrt{\tau_{y}}}{2\sqrt{\tau_{x}}}y}\hat{H}_{\pm}e^{\mp i\frac{\sqrt{\tau_{y}}}{2\sqrt{\tau_{x}}}y}$
with
\begin{multline*}
\hat{H}_{\pm}=-\partial_{x}^{2}-
\left(\partial_{y}\pm i\frac{x}{2\sqrt{1+\frac{\tau_{x}}{\tau_{y}}x^{2}}}\right)^{2}
+
\frac{\tau_{x}}{\tau_{y}(1+\frac{\tau_{x}}{\tau_{y}}x^{2})^{2}}
+ \frac{\tau_{y}}{\tau_{x}(1+\frac{\tau_{x}}{\tau_{y}}x^{2})}
\\
+\frac{1}{\tau_{x}\tau_{y}\varepsilon^{2+2\delta}}
\left[V_{\varepsilon,\tau}(x,y)\pm
  \sqrt{1+\frac{\tau_{x}}{\tau_{y}}x^{2}}\right]\,.
\end{multline*}
The remainder term in Theorem~\ref{th.BornOp} is
$$
\varepsilon^{2+4\delta}R_{1}(q,\varepsilon
D_{q},\tau,\varepsilon)+ \varepsilon^{3+2\delta} R_{2}(q,\varepsilon
D_{q},\tau,\varepsilon)\,,
$$
with $R_{1,2}\in S_{u}(1,g_{\tau};\mathcal{M}_{2}(\cz))$ and where
$R_{2}$ vanishes in a neighborhood of
$\left\{\sqrt{\tau_{x}\tau_{y}}|p|\leq r_{\gamma}\right\}$\,.
The first term provides the expected
$\mathcal{O}(\varepsilon^{2+4\delta})$
estimate in $L^{2}(\rz^{2};\cz^{2})$\,.
\\
It remains to check the effect of truncations.  All the factors,
including the left terms
$$
R_{\gamma}(p, \tau,\varepsilon)\,,\quad
\varepsilon^{2\delta}\tau_{x}\tau_{y}\left[p_{x}^{2} +
  (p_{y}\mp \varepsilon
\sin^{2}(\frac{\theta}{2}))^{2}
\right](1-\gamma(\tau_{x}\tau_{y}|p|^{2}))\,,
$$
 belong
to $OpS_{u}(\langle
\sqrt{\tau_{x}\tau_{y}}p\rangle^{2},g_{\tau};\mathcal{M}_{2}(\cz))$\,.
For any $a,b$ belonging to the class $OpS_{u}(\langle
\sqrt{\tau_{x}\tau_{y}}p\rangle^{2},g_{\tau};\mathcal{M}_{2}(\cz))$,
where $b$ vanishes in a neighborhood of
$\left\{\sqrt{\tau_{x}\tau_{y}}|p|\leq r_{\gamma}\right\}$,  and two cut-off functions
$\chi_{1},\chi_{2}\in \mathcal{C}^{\infty}_{0}((-r_{\gamma}^{2},r_{\gamma}^{2}))$ such that
$\chi_{1}\prec \chi_{2}$ (see Definition~\ref{de.ordercutoff}), the
pseudo-differential calculus says
\begin{eqnarray*}
&&(1-\chi_{2}(\tau_{x}\tau_{y}|p|^{2}))\sharp^{\varepsilon}
a\sharp^{\varepsilon} \chi_{1}(\tau_{x}\tau_{y}|p|^{2})\in \mathcal{N}_{u,g_{\tau}}\,,
\\
&& b\sharp^{\varepsilon}\chi_{1}(\tau_{x}\tau_{y}|p|^{2})\in
\mathcal{N}_{u,g_{\tau}}\,,
\end{eqnarray*}
with uniform estimates of all the seminorms w.r.t $\tau\in (0,1]^{2}$
and $\delta\in (0,\delta_{0}]$\,.
Applying this with $\chi_{1}=\chi$ and various $\chi_{2}$ such that
$\chi_{1}\prec \chi_{2}\prec \gamma$, implies that the
remainder terms due to truncations are
$\mathcal{O}(\varepsilon^{N})$ elements of $\mathcal{L}(L^{2})$ for
any $N\in \nz$\,, uniformly w.r.t $\tau\in (0,1]^{2}$ and $\delta\in
(0,\delta_{0}]$\,.
Fixing $N\geq 2+4\delta_{0}$ ends the proof of \eqref{eq.U*HUchi}.\\
For \eqref{eq.UHU*chi} use \eqref{eq.U*HUchi} with a cut-off function
$\chi_{1}$ such that $\chi\prec \chi_{1}$ and conjugate with $\hat{U}$~:
\begin{multline*}
H_{Lin}\hat{U}\chi_{1}(\tau_{x}\tau_{y}|\varepsilon D_{q}|^{2})\hat{U}^{*}
=
\\
 \varepsilon^{2+2\delta}\tau_{x}\tau_{y}
\hat{U}
\begin{pmatrix}
 e^{i\frac{\sqrt{\tau_{y}}}{2\sqrt{\tau_{x}}}y}\hat{H}_{+}e^{-i\frac{\sqrt{\tau_{y}}}{2\sqrt{\tau_{x}}}y}&0\\
0&e^{-i\frac{\sqrt{\tau_{y}}}{2\sqrt{\tau_{x}}}y}\hat{H}_{-}e^{i\frac{\sqrt{\tau_{y}}}{2\sqrt{\tau_{x}}}y}
\end{pmatrix}\chi_{1}(\tau_{x}\tau_{y}|\varepsilon
D_{q}|^{2})\hat{U}^{*}\\
+
\varepsilon^{2+4\delta}R_{1}'(\varepsilon)\,.
\end{multline*}
Right-composing  with $\chi(\tau_{x}\tau_{y}|\varepsilon D_{q}|^{2})$
and noticing that
$$\chi_{1}(\tau_{x}\tau_{y}|\varepsilon D_{q}|^{2})\hat{U}^{*}\chi(\tau_{x}\tau_{y}|\varepsilon D_{q}|^{2})-
\hat{U}^{*}\chi(\tau_{x}\tau_{y}|\varepsilon
D_{q}|^{2})=R(q,\varepsilon D_{q},\tau,\varepsilon)\,,
$$ with $R\in \mathcal{N}_{u,g_{\tau}}$ lead to
\eqref{eq.UHU*chi} like above.
\enddemonstration{}

\subsection{Linear energy estimates for non truncated states}
\label{se.linenest}
\begin{proposition}
\label{pr.adiabexpl2}
Assume $\delta\in (0,\delta_{0}]$, $\tau=(\tau_{x},\tau_{y})\in
(0,1]^{2}$ and that $V_{\varepsilon,\tau}$
 belongs to the parametric symbol class $S_{u}(\langle
   \sqrt{\frac{\tau_{x}}{\tau_{y}}}x\rangle,\frac{\frac{\tau_{x}}{\tau_{y}}dx^{2}}{\langle
     \sqrt{\frac{\tau_{x}}{\tau_{y}}}x\rangle^{2}}+
   \frac{\tau_{y}}{\tau_{x}}dy^{2})$\,. Set $\hat{\chi}=
   \chi(\tau_{x}\tau_{y}|\varepsilon D_{q}|^{2})$ for $\chi\in
   \mathcal{C}^{\infty}_{0}((-r_{\gamma}^{2},r_{\gamma}^{2}))$\,.
 When
$\hat{U}$ is the unitary semiclassical operator
$U(q,\varepsilon D_{q},\tau,\varepsilon)\in OpS_{u}(1,g_{\tau}; \mathcal{M}_{2}(\cz))$,
  given in Theorem~\ref{th.BornOp} and parametrized by a truncation
 in $\left\{\sqrt{\tau_{x}\tau_{y}}|p|< r_{\gamma}\right\}$,  then for any $\chi\in
 \mathcal{C}^{\infty}_{0}((-r_{\gamma}^{2},r_{\gamma}^{2}))$ the estimates
\begin{eqnarray}
\label{eq.U*HU}
 &&
\left| \left\langle\psi\,,\,\big[\hat{U}^{*}H_{Lin}\hat{U}
-
\varepsilon^{2+2\delta}\tau_{x}\tau_{y}
H_{BO}
\big]
\psi\right\rangle\right|
\leq
C\varepsilon^{1+2\delta}\left[\varepsilon^{1+2\delta}\|\psi\|_{L^{2}}^{2}
\right.\\
\nonumber
&&\hspace{1cm}
+\|(1-\hat{\chi})\psi\|_{L^{2}}^{2}
\left.
+
\|(1-\hat{\chi})\psi\|_{L^{2}}\times
\|\sqrt{\tau_{x}\tau_{y}}|\varepsilon D_{q}|(1-\hat{\chi})\psi\|_{L^{2}}
\right]\,,
\\
\label{eq.UHU*}
&&
 \left| \left\langle\psi\,,\,\big[H_{Lin}
-
 \varepsilon^{2+2\delta}\tau_{x}\tau_{y}
\hat{U}
H_{BO}
\hat{U}^{*}\big]
\psi\right\rangle\right|
\leq
C\varepsilon^{1+2\delta}\left[\varepsilon^{1+2\delta}\|\psi\|_{L^{2}}^{2}
\right.
\\
\nonumber
&&\hspace{1cm}\left.+\|(1-\hat{\chi})\psi\|_{L^{2}}^{2} +
\|(1-\hat{\chi})\psi\|_{L^{2}}\times
\|\sqrt{\tau_{x}\tau_{y}}|\varepsilon D_{q}|(1-\hat{\chi})\psi\|_{L^{2}}
\right]\,,
\end{eqnarray}
hold uniformly w.r.t $\tau\in (0,1]^{2}$ and $\delta\in
(0,\delta_{0}]$, for all $\psi\in L^{2}(\rz^{2}; \cz^{2})$,
with
\begin{eqnarray*}
H_{BO}&=&
\begin{pmatrix}
 e^{i\frac{\sqrt{\tau_{y}}}{2\sqrt{\tau_{x}}}y}\hat{H}_{+}e^{-i\frac{\sqrt{\tau_{y}}}{2\sqrt{\tau_{x}}}y}&0\\
0&e^{-i\frac{\sqrt{\tau_{y}}}{2\sqrt{\tau_{x}}}y}\hat{H}_{-}
e^{i\frac{\sqrt{\tau_{y}}}{2\sqrt{\tau_{x}}}y}
\end{pmatrix}\\
  \hat{H}_{\pm}&=& -\partial_{x}^{2}-
\left(\partial_{y}\pm  i\frac{x}{2\sqrt{1+\frac{\tau_{x}}{\tau_{y}}x^{2}}}\right)^{2}
+ \varepsilon^{-2-2\delta}\left[V_{\varepsilon,\tau}(x,y)\pm \sqrt{1+x^{2}}\right]
\\
&&\hspace{6cm}+
W_{\tau}(x,y)\,,
\\
W_{\tau}(x,y)&=&
\frac{\tau_{x}}{\tau_{y}(1+\frac{\tau_{x}}{\tau_{y}}x^{2})^{2}}+ \frac{\tau_{y}}{\tau_{x}(1+\frac{\tau_{x}}{\tau_{y}}x^{2})}
\,.
\end{eqnarray*}
 \end{proposition}
\begindemonstration
  Set
$$
\mathcal{D}=\hat{U}^{*}H_{Lin}\hat{U}
-\varepsilon^{2+2\delta}\tau_{x}\tau_{y}
\begin{pmatrix}
 e^{i\frac{\sqrt{\tau_{y}}}{2\sqrt{\tau_{x}}}y}\hat{H}_{+}e^{-i\frac{\sqrt{\tau_{y}}}{2\sqrt{\tau_{x}}}y}&0\\
0&e^{-i\frac{\sqrt{\tau_{y}}}{2\sqrt{\tau_{x}}}y}\hat{H}_{-}
e^{i\frac{\sqrt{\tau_{y}}}{2\sqrt{\tau_{x}}}y}
\end{pmatrix}
$$
and bound the terms $\langle \hat{\chi}\psi\,,\,
\mathcal{D} \psi\rangle$ and $\langle (1-\hat{\chi}
\psi\,,\, \mathcal{D} \hat{\chi}\psi\rangle$
by $C\varepsilon^{2+4\delta}\|\psi\|_{L^{2}}^{2}$ with the help of
 Proposition~\ref{pr.adiabexpl1}.
The remaining term is
$$
\langle (1-\hat{\chi})\psi\,,\, \mathcal{D}
(1-\hat{\chi})\psi\rangle\,,
\quad\text{with}\quad \hat{\chi}=\chi(\tau_{x}\tau_{y}|\varepsilon
D_{q}|^{2})\,.
$$
The operator $\mathcal{D}$ can be decomposed according to
$\mathcal{D}=\varepsilon^{2+2\delta}\tau_{x}\tau_{y}\mathcal{D}_{kin}+
\mathcal{D}_{pot}$ with
\begin{eqnarray}
\label{eq.Dkin}
 \hspace{-0.8cm}
\mathcal{D}_{kin}&=&\hat{U}^{*}|D_{q}|^{2}\hat{U}-
  \begin{pmatrix}
    |D_{q}-A|^{2}+|X|^{2}&0\\
 0&|D_{q}+A|^{2}+|X|^{2}
  \end{pmatrix}\,,
\\
\label{eq.Dpot}
\hspace{-0.8cm}
\mathcal{D}_{pot}
&=&
\hat{U}^{*}\hat{U}_{0}
\begin{pmatrix}
  E_{+}&0\\
0&E_{-}
\end{pmatrix}
\hat{U}_{0}^{*}\hat{U}
-V_{\varepsilon,\tau}(x,y)-
\sqrt{1+\frac{\tau_{x}}{\tau_{y}}x^{2}}
\begin{pmatrix}
  +1&0\\
0&-1
\end{pmatrix}\,.
\end{eqnarray}
The normalization in \eqref{eq.Dpot} allows to use directly the
semiclassical calculus if one remembers that
$\hat{U}_{0}^{*}\hat{U}=\Id +\varepsilon \hat{R}$
and $\hat{U}^{*}\hat{U}_{0}=\Id +\varepsilon \hat{R}'$ with
$\varepsilon^{-2\delta}R, \varepsilon^{-2\delta}R' \in
S_{u}(\frac{1}{\langle
  \sqrt{\frac{\tau_{x}}{\tau_{y}}}x\rangle},g_{\tau};\mathcal{M}_{2}(\cz))$
and
$E_{\pm}(q,\tau,\varepsilon)=V_{\varepsilon,\tau}(x,y)\pm\sqrt{1+\frac{\tau_{x}}{\tau_{y}}x^{2}}$\,.
We obtain
$$
\left|\left\langle (1-\hat{\chi})\psi\,,\,
    \mathcal{D}_{pot} (1-\hat{\chi}))\psi\right\rangle\right|
\leq C \varepsilon^{1+2\delta}\|(1-\hat{\chi})\psi\|_{L^{2}}^{2}\,,
$$
 which corresponds to the second term of our right-hand sides.\\
The kinetic energy term \eqref{eq.Dkin} is decomposed into
$$
\varepsilon^{2+2\delta}\tau_{x}\tau_{y}\mathcal{D}_{kin}=\varepsilon^{2+2\delta}\tau_{x}\tau_{y}
\mathcal{D}_{kin}^{0}+
\mathcal{D}_{kin}^{1}\,,
$$
with
\begin{eqnarray}
  \label{eq.Dkin1}
  \mathcal{D}_{kin}^{1}&=& \varepsilon^{2\delta}\hat{U}^{*}\tau_{x}\tau_{y}|\varepsilon
  D_{q}|^{2}\hat{U}-
 \varepsilon^{2\delta}\hat{U}_{0}^{*}\tau_{x}\tau_{y}|\varepsilon
 D_{q}|^{2}\hat{U}_{0}\,,\\
\label{eq.Dkin0}
  \mathcal{D}_{kin}^{0}&=& \hat{U}_{0}^{*}
  |D_{q}|^{2}\hat{U}_{0}-
  \begin{pmatrix}
    |D_{q}-A|^{2}+|X|^{2}&0\\
 0&|D_{q}+A|^{2}+|X|^{2}
  \end{pmatrix}\,.
\end{eqnarray}
Writing \eqref{eq.Dkin1} in the form
$$
\mathcal{D}_{kin}^{1}=
\varepsilon^{2\delta}(\hat{U}^{*}-\hat{U}_{0}^{*})\tau_{x}\tau_{y}|\varepsilon
D_{q}|^{2}\hat{U}+
\varepsilon^{2\delta}\hat{U}_{0}^{*}\tau_{x}\tau_{y}|\varepsilon
D_{q}|^{2}(\hat{U}-\hat{U}_{0})\,,
$$
while $\tau_{x}\tau_{y}|\varepsilon D_{q}|^{2}\in OpS_{u}(\langle
\sqrt{\tau_{x}\tau_{y}}p\rangle^{2},g_{\tau}; \mathcal{M}_{2}(\cz))$,  $\hat{U},\hat{U}_{0}\in
OpS_{u}(1,g_{\tau};\mathcal{M}_{2}(\cz))$ and
$\varepsilon^{-1-2\delta}(\hat{U}-\hat{U}_{0})\in
OpS_{u}\Big(\frac{1}{\langle \sqrt{\frac{\tau_{x}}{\tau_{y}}}x\rangle\langle
  \sqrt{\tau_{x}\tau_{y}}p\rangle^{\infty}},g_{\tau};\mathcal{M}_{2}(\cz)\Big)$
\,, leads to
$$
\left|\left\langle (1-\hat{\chi})\psi\,,\,
   \mathcal{D}_{kin}^{1}(1-\hat{\chi})
\psi\right\rangle\right|
\leq C\varepsilon^{1+4\delta}\|(1-\hat{\chi})\psi\|_{L^{2}}^{2}\,,
$$
which is even smaller than the $\mathcal{D}_{pot}$ upper bound.\\
In \eqref{eq.Dkin0}, the first term can be computed via
$$
\langle\Phi\,,\,
U_{0}^{*}|D_{q}|^{2}U_{0}(q)\Psi\rangle=\sum_{q_{j}\in \left\{x,y\right\}}
\langle
D_{q_{j}}(u_{0}(q,\tau)\Phi)\,,\, D_{q_{j}}(u_{0}(q,\tau)\Psi)\rangle\,,
$$
with $D_{q_{j}}(u_{0}f)=u_{0}(D_{q_{j}}-iu_{0}^{*}\partial_{q_{j}}u_{0})f$ and equals
\begin{eqnarray*}
\hat{U}_{0}^{*}|D_{q}|^{2}\hat{U}_{0}&=&
|D_{q}-iu_{0}^{*}\partial_{q}u_{0}|^{2}
\\
&=&
\sum_{q_{j}\in \{x,y\}}\Big[D_{q_{j}}^{2}-(iu_{0}^{*}(q)\partial_{q_{j}}u_{0}(q))^{2})D_{q_{j}}
-D_{q_{j}}(iu_{0}^{*}(q)\partial_{q_{j}}u_{0}(q))^{2})
\\
&&
\hspace{6cm}+(iu_{0}^{*}(q)\partial_{q_{j}}u_{0}(q))^{2}\Big]\,.
\end{eqnarray*}
Meanwhile expanding the entries of the second term in \eqref{eq.Dkin0}
gives
$$
|D_{q}\mp A(q)|^{2}= \sum_{q_{j}\in \left\{x,y\right\}}\left[D_{q_{j}}^{2}\mp
A_{j}(q)D_{q_{j}}\mp D_{q_{j}}A_{j}(q)
+ A_{j}(q)^{2}\right]\,.
$$
By using the expressions \eqref{eq.Dxu0} and
\eqref{eq.Dyu0} for $A, X$ and $iu_{0}^{*}\partial_{q}u_{0}$, we
obtain
\begin{eqnarray*}
&&\mathcal{D}_{kin}^{0}
=\frac{1}{2}
\begin{pmatrix}
  0 & R_{-}\\
R_{+}  &0
\end{pmatrix}\\
\text{with}&&
R_{\pm}=\pm i(D_{x}(\partial_{x}\theta)+(\partial_{x}\theta)
  D_{x})e^{\mp i\varphi}
-\sqrt{\frac{\tau_{y}}{\tau_{x}}}\sin(\theta)(D_{y}e^{\mp
  i\varphi}+e^{\mp i\varphi}D_{y})
\\
\text{and}
&&
(\partial_{x}\theta)=-\frac{\sqrt{\tau_{x}}}{\sqrt{\tau_{y}}(1+\frac{\tau_{x}}{\tau_{y}}x^{2})}
\quad,\quad
\sin(\theta)=\frac{1}{\sqrt{1+\frac{\tau_{x}}{\tau_{y}}x^{2}}}\,.
\end{eqnarray*}
Therefore, we obtain
\begin{equation*}
\left|\left\langle
    (1-\hat{\chi})\psi\,,\,\mathcal{D}_{kin}^{0}(1-\hat{\chi})
\psi\right\rangle\right|
\leq
4\max(\sqrt{\frac{\tau_{x}}{\tau_{y}}}, \sqrt{\frac{\tau_{y}}{\tau_{x}}})\||D_{q}|(1-\hat{\chi})\psi\|_{L^{2}}
\|(1-\hat{\chi})\psi\|_{L^{2}}
\,.
\end{equation*}
and, owing to $\tau_{x},\tau_{y}\in (0,1]$,
\begin{multline*}
\left|\left\langle
    (1-\hat{\chi})\psi\,,\,\varepsilon^{2+2\delta}\tau_{x}\tau_{y}\mathcal{D}_{kin}^{0}(1-\hat{\chi})
\psi\right\rangle\right|
\\
\leq
4\varepsilon^{1+2\delta}\|\sqrt{\tau_{x}\tau_{y}}|\varepsilon D_{q}|(1-\hat{\chi})\psi\|_{L^{2}}
\|(1-\hat{\chi})\psi\|_{L^{2}}\,.
\end{multline*}
This ends the proof of \eqref{eq.U*HU}.\\
For \eqref{eq.UHU*} it suffices to replace $\psi$ in \eqref{eq.U*HU}
by $\hat{U}^{*}\psi$ with a kinetic energy cut-off function $\chi_{1}$ such that
$\chi\prec \chi_{1}$ and then to use $(1-\hat\chi_{1})\hat{U}^{*}\hat\chi\in
Op\mathcal{N}_{u,g_{\tau}}$, with uniform seminorm estimates w.r.t
$\tau\in (0,1]^{2}$ and $\delta\in (0,\delta_{0}]$\,. The
$L^{2}$-norm of the corresponding additional error term is
$\mathcal{O}(\varepsilon^{N})$, for any $N$, and one fixes $N\geq 2+4\delta_{0}$\,.
\enddemonstration

 \subsection{Control of the nonlinear term}
 \label{se.nonlin}
In this subsection, we estimate the effect of the operator
$\hat{U}=U(q,\varepsilon D_{q},\tau,\varepsilon)$ belonging to
$OpS_{u}(1,g_{\tau}; \mathcal{M}_{2}(\cz))$ on the
nonlinear term $\int_{\rz^{2}}|\psi|^{4}~dxdy$\,.
\begin{proposition}
  \label{pr.L4err}
Let $\hat{U}$ be the unitary operator introduced in Theorem~\ref{th.BornOp}.
The inequalities
\begin{eqnarray}
\label{eq.L4U}
&&  \int |\psi(x,y)|^{4}~dxdy \geq (1-C\varepsilon^{1+2\delta})
\int |(\hat{U}\psi)(x,y)|^{4}~dxdy\,,\\
\label{eq.L4U*}
\text{and}&&
\int |\psi(x,y)|^{4}~dxdy \geq (1-C\varepsilon^{1+2\delta})
\int |(\hat{U}^{*}\psi)(x,y)|^{4}~dxdy
\end{eqnarray}
hold for any
$\psi\in L^{4}(\rz^{2};\cz^{2})$\,.
\end{proposition}
\begindemonstration
  For $\psi_{1}$ and $\psi_{2}$ belonging to $L^{4}(\rz^{2};\cz^{2})$,
  the local relations
  \begin{eqnarray*}
    |\psi_{1}|^{4}(q)&=&\left(|\psi_{2}(q)|^{2}+2\Real \langle
      \psi_{2}(q)\,,\,
      (\psi_{1}-\psi_{2})(q)\rangle+|\psi_{1}(q)|^{2}\right)^{2}\\
&=&
|\psi_{2}(q)|^{4}+2|\psi_{2}|^{2}|\psi_{1}-\psi_{2}|^{2}
+ 4\left(\Real\langle \psi_{2}\,,\,
  \psi_{1}-\psi_{2}\rangle+\frac{1}{2}|\psi_{1}-\psi_{2}|^{2}
\right)^{2}
\\
&&\hspace{6cm}
 +4|\psi_{2}|^{2}\Real\langle
\psi_{2}\,,\,(\psi_{1}-\psi_{2})\rangle
\\
&\geq&
|\psi_{2}(q)|^{4}-4|\psi_{2}(q)|^{3}|\psi_{1}(q)-\psi_{2}(q)|\,,
  \end{eqnarray*}
is integrated w.r.t $q=(x,y)\in \rz^{2}$, with H{\"o}lder inequality, into
$$
\|\psi_{1}\|_{L^{4}}^{4}=\int |\psi_{1}|^{4}~dxdy \geq \|\psi_{2}\|_{L^{4}}^{4} -
4\|\psi_{2}\|_{L^{4}}^{3}\|\psi_{1} - \psi_{2}\|_{L^{4}}\,.
$$
With $\psi_{1}=\hat{U}_{0}\psi=u_{0}(q)\psi$,
$|\psi(q)|^{2}=|\psi_{1}(q)|^{2}$ for all $q\in \rz^{2}$, and $\psi_{2}=
\hat{U}\psi=\psi_{1}+(\hat{U}-\hat{U}_{0})\psi$, we obtain
$$
\|\psi\|_{L^{4}}^{4}=\|\psi_{1}\|_{L^{4}}^{4}\geq
\|\hat{U}\psi\|_{L^{4}}^{4}-4\|\hat{U}\psi\|_{L^{4}}^{3}\|(\hat{U}-\hat{U_{0}})\psi\|_{L^{4}}\,.
$$
The operator $\hat{U}-\hat{U_{0}}$ equals
$\varepsilon^{1+2\delta}r(q,\varepsilon D_{q},\tau,\varepsilon)$ with
 $r$ belonging to the class $S_{u}\Big(\frac{1}{\langle \sqrt{\frac{\tau_{x}}{\tau_{y}}}x\rangle\langle
  \sqrt{\tau_{x}\tau_{y}}p\rangle^{\infty}}, g_{\tau}; \mathcal{M}_{2}(\cz)\Big)$, where we
recall $g_{\tau}=\frac{\frac{\tau_{x}}{\tau_{y}}dx^{2}}{\langle
  \sqrt{\frac{\tau_{x}}{\tau_{y}}}x\rangle^{2}}
+\frac{\tau_{y}}{\tau_{x}}dy^{2}+\frac{\tau_{x}\tau_{y}dp^{2}}{\langle\sqrt{\tau_{x}\tau_{y}}
   p\rangle^{2}}$\,. After introducing the isometric transform on
$L^{4}(\rz^{2};\cz^{2})$
$$
(T_{\varepsilon,\tau}\varphi)(x,y)=
\varepsilon\sqrt{\tau_{x}\tau_{y}}\varphi(\varepsilon \sqrt{\tau_{x}\tau_{y}} x,
\varepsilon \sqrt{\tau_{x}\tau_{y}} y)\,,
$$
the difference $\hat{U}-\hat{U}_{0}$ becomes
$$
\hat{U}-\hat{U}_{0}=\varepsilon^{1+\delta}T_{\varepsilon,\tau}^{-1}r_{1}(\varepsilon \tau_{x}x,
\varepsilon \tau_{y}y, D_{x}, D_{y},\varepsilon,\tau)T_{\varepsilon,\tau}
$$
with $r_{1}$ uniformly bounded in
$S(1, \frac{d(\tau_{x}x)^{2}}{\langle
  \tau_{x}x\rangle^{2}}+d(\tau_{y}y)^{2}+\frac{dp^{2}}{\langle
  p\rangle^{2}};\mathcal{M}_{2}(\cz))$\,. A fortiori, the symbol $r_{1}(\varepsilon \tau_{x}x,\varepsilon
\tau_{y}y,\xi,\eta;\varepsilon,\tau)$ is uniformly bounded
in
$S(1, dq^{2}+\frac{dp^{2}}{\langle p\rangle^{2}})$ and the
Lemma~\ref{le.Taylor} below provides the uniform bound
$$
\|\hat{U}-\hat{U}_{0}\|_{\mathcal{L}(L^{4})}\leq C_{0}\varepsilon^{1+\delta}\,.
$$
We have proved
$$
\|\psi\|_{L^{4}}^{4}\geq \|\hat{U}\psi\|_{L^{4}}^{4}-C\varepsilon^{1+2\delta}\|\hat{U}\psi\|_{L^{4}}^{3}\|\psi\|_{L^{4}}\,.
$$
which implies \eqref{eq.L4U}. The second inequality \eqref{eq.L4U*} is
proved similarly with $\hat{U}^{*}=\hat{U}_{0}^{*}+ (\hat{U}^{*}-\hat{U}_{0}^{*})$\,.
\enddemonstration
The result below is a particular case of the general $L^{p}$ bound,
$1<p<\infty$, for pseudodifferential operator in
$OpS(1,dq^{2}+\frac{dp^{2}}{\langle p\rangle^{2}};
\mathcal{L}(\mathcal{H}_{1};\mathcal{H}_{2}))$, $\mathcal{H}_{i}$
Hilbert spaces, stated in \cite{Tay}-Proposition~5.7  and relying on
Calderon-Zygmund analysis of singular integral operators.
\begin{lemme}
 \label{le.Taylor}
For any $p\in(1,+\infty)$, there exists a seminorm $\mathbf{n}$ on
$S(1,dq^{2}+\frac{dp^{2}}{\langle p\rangle^{2}};
\mathcal{M}_{2}(\cz))$ such that
$$
\forall a\in S(1, dq^{2}+\frac{dp^{2}}{\langle p\rangle^{2}};
\mathcal{M}_{2}(\cz)),\quad
\|a(q,D_{q})\|_{\mathcal{L}(L^{p})}
\leq \mathbf{n}(a)\,.
$$
\end{lemme}
\begindemonstration
  The Proposition~5.7 of \cite{Tay} says that for any $a\in S(1,
  dq^{2}+\frac{dp^{2}}{\langle p\rangle}^{2}; \mathcal{M}_{2}(\cz))$,
  the operator $a(q,D_{q})$ is bounded on
  $L^{p}(\rz^{2};\cz^{2})$\,. It is not difficult to follow the
  control of the constants in the previous pages of \cite{Tay} in
  order to check that $\|a(q,D_{q})\|_{\mathcal{L}(L^{p})}$ is
  estimated by a seminorm of $a$\,.
More efficiently, a linear mapping from
a Fr{\'e}chet space into a Banach space is continuous as soon as it is
bounded on bounded sets. Apply this argument  with the result of
\cite{Tay} to
$$
S(1,dq^{2}+\frac{dp^{2}}{\langle
  p\rangle^{2}};\mathcal{M}_{2}(\cz))\ni a\mapsto a(q,D_{q})\in
\mathcal{L}(L^{p}(\rz^{2};\cz^{2}))
\,.
$$
\enddemonstration

\section{Reduced minimization problems}
\label{se.explnlmin-}
In this section, we assume that the potential $V_{\varepsilon,\tau}$ satisfies
\eqref{eq.defVepsIntro}-\eqref{eq.defvIntro}.
 After the first
paragraph of this section and in the rest of the paper, we
focus on the case $\tau_{y}=1$, $\tau_{x}\to 0$  (and $\varepsilon\to
0$)\,.
Two reduced problems have to be considered: 1) the one obtained as
$\varepsilon\to 0$ and $\tau_{x}$ is fixed; 2) the one derived from
the previous one as $\tau_{x}\to 0$ and which is parametrized only by
$(G,\ell_{V})$\,.
The linear part of this latter reduced problem
is a purely quadratic Schr{\"o}dinger Hamiltonian (with a constant
magnetic field), from which many a priori information can be obtained.
This section is divided into three parts. First we specify the potential
$V_{\varepsilon,\tau}$ and check our main assumptions for the general
theory.
Then we review some properties of the reduced  Gross-Pitaevskii problem
parametrized by $(G,\ell_{V})$\,. Finally we make the comparison with
the reduced Gross-Pitaevskii problem
parametrized by $(G,\ell_{V},\tau_{x})$ as $\tau_{x}\to 0$ and deduce
properties which will be necessary for the study of the complete
minimization problem.

\subsection{Reduced minimization problems}
\label{se.uppbd}

\begin{lemme}
\label{le.Vepstau}
The potential $V_{\varepsilon,\tau}$ defined by
\eqref{eq.defVepsIntro}-\eqref{eq.defvIntro} belongs to the class
$S_{u}\left(\langle\sqrt{\frac{\tau_{x}}{\tau_{y}}}x\rangle,
 g_{q,\tau}\right)$ with the metric
$g_{q,\tau}= \frac{\tau_{x}dx^{2}}{\tau_{y}\langle 1+\frac{\tau_{x}}{\tau_{y}}x^{2}
\rangle}+\frac{\tau_{y}}{\tau_{x}}dy^{2}$\,.
\end{lemme}
\begindemonstration{}
After the change of variable
$(x',y')=(\sqrt{\frac{\tau_{x}}{\tau_{y}}}x,\sqrt{\frac{\tau_{y}}{\tau_{x}}}y)$,
it is equivalent to check
$$
\frac{\varepsilon^{2+2\delta}}{\ell_{V}^{2}}v(\tau_{y}x,\tau_{x}y)
+
\sqrt{1+x^{2}}-\varepsilon^{2+2\delta}
\left[
\frac{\tau_{x}^{2}}{(1+x^{2})^{2}}+
\frac{\tau_{y}^{2}}{1+x^{2}}
\right]\;\in S_{u}(\langle x\rangle, \frac{dx^{2}}{\langle
  x\rangle^{2}}+dy^{2})\,.
$$
It is done if
$v(\tau_{y}x,\tau_{x}y) \in S_{u}(\langle x\rangle, \frac{dx^{2}}{\langle
  x\rangle^{2}}+dy^{2})$\,.
We know $v\in S(1,\frac{dx^{2}+dy^{2}}{1+x^{2}+y^{2}})$\,.
Hence for all $(\alpha,\beta)\in \nz^{2}$ there
exists $C_{\alpha,\beta}>0$ such that
\begin{eqnarray*}
 \forall \tau\in (0,1]^{2}, \forall x,y\in\rz^{2},\quad
|\partial_{x}^{\alpha}\partial_{y}^{\beta}\left(v(\tau_{y}x,\tau_{x}y)\right)|
&&\leq
C_{\alpha,\beta}\frac{\tau_{x}^{\beta}\tau_{y}^{\alpha}}{(1+\tau_{y}^{2}x^{2}+\tau_{x}^{2}y^{2})^{-\frac{\alpha+\beta}{2}}}\\
&&\hspace{-2cm}\leq
C_{\alpha,\beta}\frac{1}{(\frac{1}{\tau_{y}^{2}}+x^{2})^{\alpha/2}}\leq
C_{\alpha,\beta}\langle x\rangle^{1-|\alpha|}\,,
\end{eqnarray*}
which is what we seek.
\enddemonstration{}
If the error terms of Proposition~\ref{pr.adiabexpl2} and
Proposition~\ref{pr.L4err}
 are assumed to be negligible, the energy
$\mathcal{E}_{\varepsilon}(\psi)$ of a state $\psi=\hat{U}^{*}
\begin{pmatrix}
  0\\a_{-}
\end{pmatrix}
$  is  close to
\begin{eqnarray}
\nonumber
&&\varepsilon^{2+2\delta}\tau_{x}\tau_{y}\langle a_{-}\,, \hat{H}_{-}
a_{-}\rangle + \frac{G_{\varepsilon,\tau}}{2}\int |a_{-}|^{4}~dxdy
=
\varepsilon^{2+2\delta}\tau_{x}\tau_{y}\mathcal{E}_{\tau}(a_{-})\,,\\
\text{with}
&&
\label{eq.defER}
\mathcal{E}_{\tau}(a_{-})= \langle a_{-}\,, \hat{H}_{-}
a_{-}\rangle + \frac{G}{2}\int |a_{-}|^{4}~dxdy\,,\\
\nonumber
&&
\hat{H}_{-}=
-\partial_{x}^{2}-
\left(\partial_{y}- i\frac{x}{2\sqrt{1+\frac{\tau_{x}}{\tau_{y}}x^{2}}}\right)^{2}
+ \frac{1}{\ell_{V}^{2}\tau_{x}\tau_{y}}v(\sqrt{\tau_{x}\tau_{y}}x,\sqrt{\tau_{x}\tau_{y}}y)\,.
\end{eqnarray}
with the potential
 $v$ chosen from (\ref{eq.defvIntro}).\\
When $\frac{\tau_{x}}{\tau_{y}}$ and $\tau_{x}\tau_{y}$ are small, in
particular in the regime $\tau_{x}\ll 1$ and $\tau_{y}=1$ that we
shall consider,
this energy is well approximated by
\begin{equation}
\label{eq.defEH}
\mathcal{E}_{H}(\varphi)=
\langle
\varphi\,,\,
\left[-\partial_{x}^{2}-(\partial_{y}-\frac{ix}{2})^{2}
  + \frac{x^{2}+y^{2}}{\ell_{V}^{2}}\right]\varphi\rangle
+\frac{G}{2}\int|\varphi|^{4}\,,
\end{equation}
as this will be checked and specified in the next paragraph.
Although more general asymptotics could be considered, we concentrate
from now on the regime $\tau_{x}\ll 1$, $\tau_{y}=1$\,. The parameters
$\ell_{V}$ and $G$ are assumed to be fixed as $\tau_{x}\to 0$\,.

In order to prove that the ground states of ${\mathcal E}_{H}$ and ${\mathcal E}_{\varepsilon}$
are close, we need good estimates on the energy ${\mathcal E}_{H}$.
\subsection{Properties of the harmonic approximation}
\label{se.harmapp}

The energy functional ${\mathcal
  E}_{H}$  does not any more depend on $\tau_{x}$ and is
parametrized only by $(G,\ell_{V})$\,. Let us start with its properties.
 We introduce the spaces ${\mathcal H}_{1}$ and ${\mathcal H}_{2}$ which are
given by
\begin{equation}
  \label{eq.defcalHs}
{\mathcal H}_{s}=\left\{u\in L^{2}(\rz^{2})\,
 \sum_{|\alpha|+|\beta|\leq s} \|q^{\alpha}D_{q}^{\beta}u\|_{L^{2}}
 <+\infty\right\}\,,\quad s=1,2, (q=(x,y))
\end{equation}
endowed with the norm $\|u\|_{{\mathcal
    H}_{s}}^{2}=\sum_{|\alpha|+|\beta|\leq s}
\|q^{\alpha}D_{q}^{\beta}u\|_{L^{2}}^{2}$\,.
For a compact set $K$ of ${\mathcal H}_{s}$ and for $u\in {\mathcal H}_{s}$,
the distance $d_{s}(u,K)$ follows the usual definition $\min_{v\in K}\|u-v\|_{{\mathcal H}_{s}}$\,.
The self-ajoint operator associated with the linear part of ${\mathcal
  E}_{H}$ is denoted by
$$
H_{\ell_{V}}= -\partial_{x}^{2}-(\partial_{y}-\frac{ix}{2})^{2}
  + \frac{x^{2}+y^{2}}{\ell_{V}^{2}}\,, \quad (0<\ell_{V}<+\infty)\,.
$$
Its domain is ${\mathcal H}_{2}$ while its form domain is ${\mathcal
  H}_{1}$\,. Note also the compact embeddings ${\mathcal
  H}_{2}\subset\subset {\mathcal H}_{1}\subset\subset L^{2}\cap L^{4}$\,.
Following the general scheme presented in \cite{HiPr,Sjo},
its spectrum equals
$$
\sigma(H_{\ell_{V}})=\left\{(1+2n_{+})r_{+}+(1+2n_{-})r_{-}\,,\quad
  (n_{+},n_{-})\in \nz^{2}\right\}
$$
with $r_{\pm}=\frac{1}{2\sqrt{2}}\sqrt{1+\frac{8}{\ell_{V}^{2}}\pm\sqrt{1+\frac{4}{\ell_{V}^{2}}}}$\,.
\begin{proposition}
\label{pr.harmapp}
The functional ${\mathcal E}_{H}$ admits minima on $\left\{u\in {\mathcal
    H}_{1}, \|u\|_{L^{2}}=1\right\}$\,, with a minimum value ${\mathcal
  E}_{H,min}$ satisfying
$$
{\mathcal E}_{H,min}\geq r_{+}+r_{-}\geq \frac{\sqrt{2}}{2}\,.
$$
The set of minimizers $\mathrm{Argmin~}{\mathcal E}_{H}$ is a bounded subset of
${\mathcal H}_{2}$ and therefore a compact subset of $\left\{u\in {\mathcal
    H}_{1}, \|u\|_{L^{2}}=1\right\}$\,. Moreover for any $\varphi\in
\mathrm{Argmin~}{\mathcal E}_{H}$, $\varphi$ is an eigenvector of $H_{0}+
G|\varphi|^{2}$\,.
Finally there exist two constants $C=C_{\ell_{V},G}>0$ and
$\nu=\nu_{\ell_{V},G}\in (0,1/2]$ such that the conditions $u\in {\cal
  H}_{1}$,
$\|u\|_{L^{2}}=1$ and ${\mathcal E}_{H}(u)\leq {\mathcal E}_{H,min}+1$, imply
\begin{equation}
  \label{eq.lojas}
d_{{\mathcal H}_{1}}(u,\mathrm{Argmin}~{\mathcal E}_{H})\leq C({\mathcal E}_{H}(u)-{\mathcal E}_{H,min})^{\nu}\,.
\end{equation}
\end{proposition}
\begindemonstration{}
On ${\mathcal H}_{1}$ and ${\mathcal H}_{2}$, the scalar products
\begin{eqnarray*}
&&\langle u\,,\, v\rangle_{1,\ell_{V}}=\langle u\,,\,
H_{\ell_{V}}v\rangle_{L^{2}}
\\
&&
\langle u\,,\, v\rangle_{2,\ell_{V}}=\langle H_{\ell_{V}}u\,,\,
H_{\ell_{V}}v\rangle_{L^{2}}
\end{eqnarray*}
provide norms $\|u\|_{k,\ell_{V}}$, $k=1,2$, respectively  equivalent to $\|u\|_{{\mathcal
    H}_{k}}$\,. In this proof, all the ``uniform'' estimates are
actually parametrized by $(G,\ell_{V})$\,.
The  nonlinearity $\frac{G}{2}\int|u|^{4}(x)~dx$ as well as the
constraint $\|u\|_{L^{2}}=1$ are continuous functions on $L^{2}\cap
L^{4}$ while the quadratic part of ${\mathcal E}_{H}(u)$ is simply
$\|u\|_{1,\ell_{V}}^{2}$ with ${\mathcal E}_{H}(u)\geq
\|u\|_{1,\ell_{V}}^{2}\geq r_{+}+r_{-}\geq \frac{\sqrt{2}}{2}$\,.
The compact embedding ${\mathcal H}_{1}\subset\subset L^{2}\cap L^{4}$
thus implies that the infimum $\inf_{u\in {\mathcal H}_{1}\,,\,
  \|u\|_{L^{2}}=1}{\mathcal E}_{H}(u)$ is achieved.\\
A minimizer $\varphi\in \mathrm{Argmin}~{\mathcal E}_{H}$ solves in a
distributional sense the Euler-Lagrange equation
$$
H_{\ell_{V}}\varphi +G|\varphi|^{2}u=\lambda_{\varphi}\varphi
$$
where $\lambda_{\varphi}$ is the Lagrange multiplier associated with the
constraint $\|\varphi\|_{L^{2}}=1$\,. By taking the scalar product with $\varphi$,
one obtains the  bounds for $\lambda_{\varphi}$:
$$
{\mathcal E}_{H,min}\leq \lambda_{\varphi}\leq 2{\mathcal E}_{H,min}\,.
$$
Since ${\mathcal H}_{1}$ is also (compactly) embedded in $L^{6}(\rz^{2})$,
the equation
$$
H_{\ell_{V}}\varphi = -G|\varphi|^{2}\varphi+\lambda_{\varphi}\varphi
$$
ensures that $\|\varphi\|_{2,\ell_{V}}$ is uniformly bounded on
$\mathrm{Argmin}~{\mathcal E}_{H}$\,.
Therefore $\mathrm{Argmin}~{\mathcal E}_{H}$ is a bounded subset of ${\mathcal
  H}_{2}$ and a compact subset of ${\mathcal H}_{1}$\,.
In an ${\mathcal H}_{1}$-neighborhood of $\varphi\in \mathrm{Argmin}~{\mathcal
  E}_{H}$ ($\|\varphi\|_{L^{2}}=1$), the $L^{2}$-sphere $\left\{u\in
  {\mathcal H}_{1},\, \|u\|_{L^{2}}=1\right\}$ can be parametrized by
$$
u=(1-\|v\|_{L^{2}}^{2})\varphi + v\,,\quad \langle \varphi\,,\, v\rangle_{L^{2}}=0\,.
$$
 Notice also that the potential
$G|\varphi|^{2}$ is a relatively compact perturbation of
$H_{\ell_{V}}$, so that $H_{\ell_{V}}+ G|\varphi|^{2}$ is a
self-adjoint operator in $L^{2}(\rz^{2})$ with  domain ${\mathcal H}_{2}$
and with a compact resolvent.
With $\langle \varphi\,,\,
v\rangle_{1,\ell_{V}}=-G\int_{\rz^{2}}|\varphi|^{2}\overline{\varphi}v~dx$
for $\varphi$ is
an eigenvector of $H_{\ell_{V}}+G|\varphi|^{2}$ and $v\perp \varphi$,
the energy ${\mathcal E}_{H}(u)$ becomes
\begin{eqnarray*}
  {\mathcal E}_{H}((1-\|v\|_{L^{2}}^{2})\varphi+v)&=&\|v\|_{1,\ell_{V}}^{2} +(1-\|v\|_{L^{2}})\langle
  \varphi\,,\,\varphi\rangle_{1,\ell_{V}}
\\
&&\hspace{-3cm}
-2G(1-\|v\|_{L^{2}}^{2})\Real\int_{\rz^{2}}|\varphi|^{2}\overline{\varphi}v~dx
+\frac{G}{2}\int_{\rz^{2}}|(1-\|\varphi\|_{L^{2}}^{2})\varphi+v|^{4}~dx\\
&=&\|v\|_{1,\ell_{V}}^{2}+ F_{\varphi}(v)\,,
\end{eqnarray*}
where $v$ lies in the closed subset ${\mathcal H}_{1,\varphi}=\left\{v\in {\mathcal H}_{1}, \langle
  \varphi\,, v\rangle_{L^{2}}=0\right\}$ of ${\mathcal H}_{1}$ and
$F_{\varphi}(v)$ is the composition of the compact embedding ${\mathcal
  H}_{1}\to L^{2}\cap L^{4}$ with a real analytic, real-valued,
functional on $L^{2}\cap L^{4}$\,. Hence on ${\mathcal H}_{1,\varphi}$ endowed
with the scalar product $\langle~,~\rangle_{1,\ell_{V}}$, the Hessian
of ${\mathcal E}_{H}((1-\|v\|_{L^{2}}^{2})\varphi+v)$ equals $\Id+
D^{2}F_{\ell_{V}}(0)$, with $D^{2}F_{\ell_{V}}(0)$ compact (and
self-adjoint)\,.
We can apply the Lojasiewicz-Simon inequality which says that there
exist two constants $C_{\varphi}>0$, $\nu_{\varphi}\in (0,1/2]$, such that
$$
\|v\|_{1,\ell_{V}}\leq C_{\varphi}\left({\mathcal E}_{H}((1-\|v\|_{L^{2}}^{2})\varphi+v)-{\mathcal E}_{H,min}\right)^{\nu_{\varphi}}\,.
$$
Since the set $\mathrm{Argmin}~{\mathcal E}_{H}$ is a compact subset of
${\mathcal H}_{1}$, it can be covered by a finite number of neighborhoods
of $\varphi_{i}\in \mathrm{Argmin}~{\mathcal E}_{H}$, $1\leq i\leq N$,  where a
Lojasiewicz-Simon inequality holds. Take
$$
\nu_{\ell_{V},G}=\min_{1\leq i\leq
  N}\nu_{\varphi_{i}}\quad\text{and}\quad
C_{\ell_{V},G}=2\max_{1\leq i\leq N}C_{\varphi_{i}}\,.
$$
\enddemonstration{}
\begin{remarque}
\label{re.loja}
  The Lojasiewicz inequality is a classical result of real algebraic
  geometry (see a.e. \cite{Loj,BCR}) proved by Lojasiewicz
after Tarski-Seidenberg Theorem. It is usually written as
  $|\nabla f(x)|\leq C |f(x)|^{\nu}$ with $\nu\in (0,1]$ for a real
  analytic function of $x$ lying around $x_{0}$ with
  $f(x_{0})=0$\,. The variational form is a variant of it. It was
  extended to the infinite dimensional case with applications to PDE's
 by L.~Simon in \cite{Sim}.
We refer the reader also to \cite{Chi,HaJe,Hua} and
\cite{BDLM}
for recent
texts and references concerned with the infinite dimensional case or
 the extension with $o$-minimal structures.
\item The nonlinear Euler-Lagrange equation is usually studied after
  linearization via the Liapunov-Schmidt process.
  Here using some coordinate representation of the
  constraint submanifold, especially when it is a sphere for a simple
  norm, allows to use directly the standard result for the
  minimization of real analytic functionals.

\item When the minimization problem is non degenerate at every
  $\varphi\in \mathrm{Argmin}~{\mathcal E}_{H}$, i.e. in the present
  case when the kernel of $\Id+D^{2}F_{\varphi}(0)$  is restricted to
  $\left\{0\right\}$, the compact set $\mathrm{Argmin}~{\mathcal E}_{H}$ is
  made of a finite number of point. When $\ell_{V}$ is fixed so that
  $r_{+}$ and $r_{-}$ are rationally independent, the spectrum of
  $H_{\ell_{V}}$ is made of simple eigenvalues and when $G$ is small enough,
  $G < G_{\ell_{V}}$, the non degeneracy assumption is satisfied via a perturbation
  argument from the case $G=0$\,. For large $G$, we can only say that
  the set of $(G,\ell_{V})\in (0,+\infty)^{2}$ such that all the
  minima are non degenerate, $\nu_{\ell_{V},G}=1/2$, is a subanalytic
  subset of $\rz^{2}$\,. In our case with a linear part
  $H_{\ell_{V}}$ which is a complex operator with no rotational
  symmetry, no standard methods like in \cite{AJR} allow to reduce the
  minimization problem to some radial nonlinear ODE.

 From the
  information given by the Lowest-Landau-Level reduction, when $G$ and
  $\ell_{V}$ are large, the supposed hexagonal symmetry, after
  removing some trivial rotational invariance,  of the
  problem (see \cite{ABN,Nie}) suggests that there are
  presumably several minimizers.
\end{remarque}

A change of variable $\varphi(x,y)e^{-ixy/4}=\alpha u(\alpha x,\alpha y)$ with $\alpha^2=1/(\ell_{V}\sqrt{G})$
 leads to \begin{equation}{\mathcal E}_H(\varphi)=\frac 1{\ell_{V}\sqrt{G}}\tilde { \mathcal E}_H(u)=\frac 1{\ell_{V}\sqrt{G}}\int
 |(\nabla - \frac i 2 \ell_{V}\sqrt{G} e_{z}\times \vec{r})
 u|^2+G(r^2|u|^2+\frac 12 |u|^4)\,,
\end{equation}
with the notations $\vec{r}=
\begin{pmatrix}
  x\\y
\end{pmatrix}
$ $r=|\vec{r}|$\,.
 This implies that $\ell_{V}\sqrt{G}$ is equivalent to a rotation value. We have the following results from
  the literature\begin{itemize}\item when $\ell_{V}\sqrt{G}$ is small and $G$ is large, the minimizer is unique
   up to rotation and vortex-free \cite{AJR}: namely $u(x,y)=f(r)e^{ic}$ for some real number $c$, where $f$ does not vanish.
   If $\ell_{V}\sqrt{G}=0$, this is an adaptation of a result of \cite{BrOs}. When $\ell_{V}\sqrt{G}$ is non zero, this
    requires refined estimates for the jacobian.
   \item when $\ell_{V}\sqrt{G}$ is large, then vortices are expected in the system and this can be
   analyzed in details in the LLL regime (lowest Landau level) if additionally $\sqrt G/\ell_{V}$ is small \cite{AB,ABN}.
   More precisely, if $\sqrt G/\ell_{V}$ is small, then
   \begin{equation} \inf \tilde {\mathcal E}_H -\frac 12- \inf E_{LLL}=o \left (\frac { \sqrt G}{\ell_{V}}\right )\end{equation}
   where \begin{equation}E_{LLL}(u)=\int G(r^2|u|^2+\frac 12 |u|^4)\end{equation} for functions $u$ such that $u(x,y)e^{\ell_{V}\sqrt{G}r^2/4}$
    is a holomorphic function of $x+iy$. This space is called the LLL. If $u$ is a ground state of $\tilde {\mathcal E}_H$
     and $w$ its projection onto the LLL, then $|u-w|$ tends to 0 in $H^1$ and $C^{0,\alpha}$ as $\sqrt G/\ell_{V}$
      tends to 0. If additionally, $\ell_{V}\sqrt{G}$ is large, then one can estimate $\inf E_{LLL}$ \cite{ABN} thanks to test
      functions with vortices and $\inf E_{LLL}=O\left (\frac { \sqrt G}{\ell_{V}}\right )$.\item if $\ell_{V}\sqrt{G}$ is large,
       and $\sqrt G/\ell_{V}$ is large, then this is a Thomas Fermi regime where the energy can be estimated as well \cite{Aft}
        and is of order $\sqrt G/\ell_{V}$.
\end{itemize}

We complete the previous result with another comparison statement
which will be useful in the sequel.
\begin{proposition}
\label{pr.harmH2}
There exists $C=C_{\ell_{V},G}>0$ such that when
$u\in \mathcal{H}_{1}$ satisfy $\mathcal{E}_{H}(u)\leq
\mathcal{E}_{H,min}+1$, $\|u\|_{L^{2}}=1$, and solves
$$
H_{\ell_{V},G}u+G|u|^{2}u= \lambda_{u} u + r
$$
with $\lambda_{u}\in \rz$ and $r\in L^{2}$, then
\begin{itemize}
\item $u\in \mathcal{H}_{2}$\,;
\item there exists $u_{0}\in \mathrm{Argmin}~\mathcal{E}_{H}$, with
  Lagrange multiplier $\lambda_{u_{0}}$,  such that
$$
|\lambda_{u}-\lambda_{u_{0}}|+\|u-u_{0}\|_{\mathcal{H}_{2}}\leq
C\left(\|r\|_{L^{2}}+ \left(\mathcal{E}_{H}(u)-\mathcal{E}_{H,min}\right)^{\nu}\right)\,,
$$
where $\nu=\nu_{\ell_{V},G}\in (0,\frac{1}{2}]$ is the exponent given
in Proposition~\ref{pr.harmapp}.
\end{itemize}
\end{proposition}
\begindemonstration{}
Since $\mathrm{Argmin}~\mathcal{E}_{H}$ is compact,
Proposition~\ref{pr.harmapp} already  provides $u_{0}\in
\mathrm{Argmin}~\mathcal{E}_{H}$ such that
$$
\|u-u_{0}\|_{\mathcal{H}_{1}}\leq
C(\mathcal{E}_{H}(u)-\mathcal{E}_{H,min})^{\nu}\,.
$$
Taking the difference of the equation for $u$ and the Euler-Lagrange
equation for $u_{0}$, we obtain
$$
H_{\ell_{V},G}(u-u_{0})=\left(\lambda_{u}-\lambda_{u_{0}}\right)u_{0}+
\lambda_{u}(u-u_{0})+G(|u_{0}|^{2}u_{0}-|u|^{2}u)+r\,.
$$
Taking the scalar product with $u_{0}$, with
$$
\||u_{0}|^{2}u_{0}-|u|^{2}u\|_{L^{2}}\leq
C\left(\mathcal{E}_{H}(u_{0})+\mathcal{E}_{H}(u)\right)
\|u-u_{0}\|_{\mathcal{H}_{1}}\leq C'\left(\mathcal{E}_{H}(u)-\mathcal{E}_{H,min}\right)^{\nu}
$$
implies
$$
|\lambda_{u}-\lambda_{u_{0}}|\leq C'' \left(\|r\|_{L^{2}}+ \left(\mathcal{E}_{H}(u)-\mathcal{E}_{H,min}\right)^{\nu}\right)\,.
$$
Using the ellipticity of $H_{\ell_{V},G}$ and the equivalence of the
norms $\|\varphi\|_{\mathcal{H}_{2}}$ and $\|H_{\ell_{V},G}\varphi\|_{L^{2}}$ ends
the proof.
\enddemonstration{}
\subsection{Comparison of the two reduced minimization problems}
\label{se.compred}

In the regime $\tau_{y}=1$ and $\tau_{x}\to 0$, while $\ell_{V}>0$ and
$G>0$ are fixed, we compare the two minimization problems for the energies ${\mathcal
  E}_{\tau}$ and ${\mathcal E}_{H}$ defined in
\eqref{eq.defER}-\eqref{eq.defEH}.
We start with the next Lemma which is a simple application of the so
called IMS localization formula (see a.e. \cite{CFKS}). We shall use
the functional spaces ${\mathcal H}_{s}$ defined by \eqref{eq.defcalHs}
associated with $\mathcal{E}_{H}$ as well as the standard Sobolev
spaces $H^{s}(\rz^{2})$ associated with $\mathcal{E}_{\tau}$, with
$s=1,2$ and $\mathcal{H}_{s}\subset H^{s}(\rz^{2})$\,.
\begin{lemme}
\label{le.compred}
Let $\chi_{1}, \chi_{2} \in \mathcal{C}_{b}^{\infty}(\rz^{2})$ satisfy
$\chi_{1}^{2}+\chi_{2}^{2}=1$, $\supp \chi_{1}\subset
\left\{x^{2}+y^{2}< 1\right\}$ and take $\alpha\in (0,\frac{1}{2}]$\,. Then the
following identity
\begin{multline}
  \label{eq.IMSH}
\mathcal{E}_{H}(u)=
\mathcal{E}_{H}(\chi_{1}(\tau_{x}^{\alpha}.)u)+\mathcal{E}_{H}(\chi_{2}(\tau_{x}^{\alpha}.)u)
-
\tau_{x}^{2\alpha}\sum_{j=1}^{2}\int_{\rz^{2}}|(\nabla\chi_{j})(\tau_{x}^{\alpha}.)|^{2}|u|^{2}
\\
+G\int_{\rz^{2}}(\chi_{1}^{2}\chi_{2}^{2})(\tau_{x}^{\alpha}.)|u|^{4}\,,
\end{multline}
holds for all $u\in \mathcal{H}_{1}$,
with the same formula for $\mathcal{E}_{\tau}(u)$ when $u\in
H^{1}(\rz^{2})$\,. Moreover, $\mathcal{E}_{\tau}$ and $\mathcal{E}_{H}$ satisfy
\begin{multline}
  \label{eq.IMSR}
\mathcal{E}_{\tau}(u)=\mathcal{E}_{H}(\chi_{1}(\tau_{x}^{\alpha}.)u)+\mathcal{E}_{\tau}(\chi_{2}(\tau_{x}^{\alpha}.)u)
- \tau_{x}^{2\alpha}\sum_{j=1}^{2}\int_{\rz^{2}}|(\nabla
\chi_{j})(\tau_{x}^{\alpha}.)|^{2}|u|^{2}\\
+G\int_{\rz^{2}}(\chi_{1}^{2}\chi_{2}^{2})(\tau_{x}^{\alpha}.)|u|^{4}+ R(u)\,,
\end{multline}
for all $u\in H^{1}(\rz^{2})$ with
$$
|R(u)|\leq\frac{1}{4}
\left(\mathcal{E}_{\tau}(\chi_{1}(\tau_{x}^{\alpha}.)u)^{1/2}+
\mathcal{E}_{H}(\chi_{1}(\tau_{x}^{\alpha}.)u)^{1/2}\right)\|u\|_{L^{2}}
\tau_{x}^{1-3\alpha}\,.
$$
\end{lemme}
\begindemonstration{}
The first identity is a direct application of the IMS localization
formula (see a.e. \cite{CFKS}) which comes from the identity
$$
P\chi^{2}P-\chi P^{2}
\chi=\left[P,\chi\right]^{2}-\frac{1}{2}\left[\chi^{2},P\right]P
-\frac{1}{2}P\left[P,\chi^{2}\right]
$$
when $P$ is a  differential operator of order $\leq 1$ and $\chi$ is a
$\mathcal{C}^{\infty}$ function. Simply combine it with the identity
$$
|u|^{4}=
|\chi_{1}(\tau_{x}^{\alpha}.)u|^{4}+|\chi_{2}(\tau_{x}^{\alpha}.)u|^{4}+
2\chi_{1}^{2}\chi_{2}^{2}(\tau_{x}^{\alpha}.)|u|^{4}\,.
$$
Using the same argument for $\mathcal{E}_{\tau}$ provides the same
identity after replacing $\mathcal{E}_{H}$ with $\mathcal{E}_{\tau}$, and
it suffices to compare
$\mathcal{E}_{\tau}(\chi_{1}(\tau_{x}^{\alpha}.)u)$ with $\mathcal{E}_{H}
(\chi_{1}(\tau_{x}^{\alpha}.)u)$\,. The definition \eqref{eq.defvIntro} of
the potential $v$ and the condition $\alpha\leq \frac{1}{2}$
imply
$$
v(\tau_{x}^{1/2}.)= \tau_{x}(x^{2}+y^{2})\,\quad\text{on}~\supp \chi_{1}(\tau_{x}^{\alpha}.)\,.
$$
Therefore, we obtain, by setting $u_{\tau}=
\chi_{1}(\tau_{x}^{\alpha}.)u$\,,
\begin{eqnarray*}
&&  |\mathcal{E}_{\tau}(\chi_{1}(\tau_{x}^{\alpha}.)u)-
  \mathcal{E}_{H}(\chi_{1}(\tau_{x}^{\alpha}.)u)|
=
\\
&&\qquad
\left|
\int_{\rz^{2}}|(\partial_{y}-i\frac{x}{2\sqrt{1+\tau_{x}x^{2}}})u_{\tau}|^{2}
- |(\partial_{y}-i\frac{x}{2})u_{\tau}|^{2}\right|\\
&&\leq
\left(
\|(\partial_{y}-i\frac{x}{2\sqrt{1+\tau_{x}x^{2}}})u_{\tau}\|_{L^{2}}+
\|(\partial_{y}-i\frac{x}{2})u_{\tau}\|_{L^{2}}
\right)
\\
&&\hspace{4cm}\times
\|\frac{\tau_{x}x^{3}/2}{1+\sqrt{1+\tau_{x}x^{2}}}\chi_{1}(\tau_{x}^{\alpha}.)u\|_{L^{2}}
\\
&&
\leq\frac{1}{4}
\left(\mathcal{E}_{\tau}(\chi_{1}(\tau_{x}^{\alpha}.)u)^{1/2}+
\mathcal{E}_{H}(\chi_{1}(\tau_{x}^{\alpha}.)u)^{1/2}\right)\|u\|_{L^{2}}
\tau_{x}^{1-3\alpha}\,.
\end{eqnarray*}
\enddemonstration{}
\begin{proposition}
\label{pr.redHR}
For any given $(\ell_{V},G)\in (0,+\infty)^{2}$, there exists
$\tau_{\ell_{V},G}>0$ such that the following properties hold when
$\tau_{x}\leq \tau_{\ell_{V},G}$\,.\\
\begin{itemize}
\item The minimization problem
$$
\inf_{u\in H^{1}(\rz^{^{2}})\,,\, \|u\|_{L^{2}}=1}\mathcal{E}_{\tau}(u)
$$
admits a solution $u\in H^{1}(\rz^{2})$\,.
\item A solution $u\in H^{1}(\rz^{2})$ to the above minimization
  problem,
solves an Euler-Lagrange equation
$$
\left[-\partial_{x}^{2}-(\partial_{y}-\frac{i}{2}\frac{x}{\sqrt{1+\tau_{x}x^{2}}})^{2}
  + \frac{v(\tau_{x}^{1/2}.)}{\ell_{V}^{2}\tau_{x}}+G|u|^{2}\right]u=\lambda_{u}u
$$
with $0\leq \lambda_{u}\leq 2\mathcal{E}_{\tau,min}$ and belongs to $H^{2}(\rz^{2})$\,.
\item Moreover the minimum value
$\mathcal{E}_{\tau,min}=\min_{u\in H^{1}(\rz^{2})\,,\,
  \|u\|_{L^{2}}=1}\mathcal{E}_{\tau}(u)$ satisfies the estimate
$$
|\mathcal{E}_{\tau,min}-\mathcal{E}_{H,min}|\leq
C_{\ell_{V},G}\tau_{x}^{2/3}\,.
$$
\item For $u\in \mathrm{Argmin}~\mathcal{E}_{\tau}$ and any pairs $\chi=(\chi_{1},\chi_{2})$ in $
  \mathcal{C}^{\infty}_{b}(\rz^{2})^{2}$ such that
  $\chi_{1}^{2}+\chi_{2}^{2}=1$ with $\supp
  \chi_{1}\subset\left\{x^{2}+y^{2}<1\right\}$ and $\chi_{1}\equiv 1$
  in $\left\{x^{2}+y^{2}\leq 1/2\right\}$, the functions
  $\chi_{j}(\tau_{x}^{1/9}.)u$, $j=1,2$, satisfy
  \begin{eqnarray}
    \label{eq.estimu2}
&&    \|\chi_{2}(\tau_{x}^{1/9}.)u\|^{2}_{L^{2}}\leq
    C_{\chi,\ell_{V},G}\tau_{x}^{2/3}\\
&&
\label{eq.estimEnu1}
\mathcal{E}_{\tau}(\chi_{1}(\tau_{x}^{1/9}.)u)\leq \mathcal{E}_{\tau,min}+ C_{\chi,\ell_{V},G}\tau_{x}^{2/3}
\\
\label{eq.estimH2u1}
&& d_{\mathcal{H}_{2}}(\chi_{1}(\tau_{x}^{1/9}.)u,
\mathrm{Argmin}~\mathcal{E}_{H})\leq
C_{\chi,\ell_{V},G}\tau_{x}^{2\nu_{\ell_{V,G}}/3}\,,\\
&&
\nu_{\ell_{V},G}\in (0,\frac{1}{2}]\,.
  \end{eqnarray}
\end{itemize}
A constant $C_{a,b,c}$ is a constant which is fixed once $(a,b,c)$
are given.
\end{proposition}
\begindemonstration{}
Fix $\ell_{V}$ and $G$. We drop the indices $\ell_{V},G$ in the constants. The
exponent $\alpha$ will be fixed to the value $\frac{1}{9}$ within the proof.\\
\noindent\underline{First step, upper bound for
  $\inf\{\mathcal{E}_{\tau}(u)\,,\; u\in H^{1}(\rz^{2})\,,\,
    \|u\|_{L^{2}}=1\}$:}\\
Let $\chi=(\chi_{1},\chi_{2})$
and $\tilde{\chi}=(\tilde{\chi}_{1},\tilde{\chi_{2}})$ be two pairs as
in our statement such that $\tilde{\chi}_{1}\prec\chi_{1}$ according to
Definition~\ref{de.ordercutoff}.
Take $u_{0}\in \mathrm{Argmin}~\mathcal{E}_{H}\subset
\mathcal{H}_{1}$\,. According to Proposition~\ref{pr.harmapp}, it
belongs to a bounded set
of $\mathcal{H}_{2}$ so that $\||q|^{2}u_{0}\|_{L^{2}}$, with
$q=(x,y)$,  is uniformly
bounded. Hence,
$0\not\in\supp\nabla\tilde{\chi}_{j}\cup \supp
\tilde{\chi}_{1}\tilde{\chi}_{2}$ implies
$$
\int_{\rz^{2}}|\nabla\tilde{\chi}_{j}(\tau_{x}^{\alpha}.)|^{2}|u_{0}|^{2}=
\mathcal{O}(\tau_{x}^{4\alpha})\quad\text{while}\quad
\int_{\rz^{2}}\tilde{\chi}_{1}^{2}\tilde{\chi}_{2}^{2}(\tau_{x}^{\alpha}.)|u_{0}|^{4}
\geq 0\,.
$$
Lemma~\ref{le.compred} above with the pair $\tilde{\chi}$
and $\alpha\in (0,\frac{1}{2}]$ gives:
$$
\mathcal{E}_{H, min}=\mathcal{E}_{H}(u_{0})
\geq \mathcal{E}_{H}(\tilde{\chi}_{1}(\tau_{x}^{\alpha}. )u_{0})
+
\mathcal{E}_{H}(\tilde{\chi}_{2 }(\tau_{x}^{\alpha}. )u_{0})- C\tau_{x}^{6\alpha}\,.
$$
On $\supp \tilde{\chi}_{2}(\tau_{x}^{\alpha}.)$, the potential
$\frac{v(\tau_{x}^{1/2}.)}{\ell_{v}^{2}\tau_{x}}$ is bounded from below by
$\frac{1}{C'\tau_{x}^{2\alpha}}$\,. Thus we get
$$
\mathcal{E}_{H,min}\left(\|\tilde{\chi_{1}}u_{0}\|_{L^{2}}^{2}+\|\tilde{\chi}_{2}u_{0}\|_{L^{2}}^{2}\right)
\geq \mathcal{E}_{H,min}\|\tilde{\chi_{1}}u_{0}\|_{L^{2}}^{2}
+ \frac{1}{C'\tau_{x}^{2\alpha}}\|\tilde{\chi}_{2}u_{0}\|_{L^{2}}^{2} -C\tau_{x}^{6\alpha}\,,
$$
and finally
$$
\|\tilde{\chi}_{2}u_{0}\|_{L^{2}}^{2}\leq C''\tau_{x}^{8\alpha}\,,\quad
\|\tilde{\chi}_{1}u_{0}\|_{L^{2}}^{2}=1+\mathcal{O}(\tau_{x}^{8\alpha})\,,
$$
as soon as
$\tau_{x}<(C'\mathcal{E}_{H,min})^{-1/2\alpha}$\,.\\

The function $u_{1}=
\|\tilde{\chi}_{1}u_{0}\|_{L^{2}}^{-1}\tilde{\chi}_{1}u_{0}$ is
normalized with
$$
\mathcal{E}_{H,min}\leq \mathcal{E}_{H}(u_{1})\leq
\mathcal{E}_{H,min}+
\mathcal{O}(\tau_{x}^{6\alpha})\,,
$$
and $\chi_{1}(\tau_{x}^{\alpha}.)u_{1}=u_{1}$,
$\chi_{2}(\tau_{x}^{\alpha}.)u_{1}=0$\,.
Applying the second formula of Lemma~\ref{le.compred} with, now, the
pair $\chi$, leads to
\begin{eqnarray*}
&&  \mathcal{E}_{\tau}(u_{1})= \mathcal{E}_{H}(u_{1})+ R(u_{1})
= \mathcal{E}_{H,min}+ \mathcal{O}(\tau_{x}^{6\alpha})+ R(u_{1})\\
\text{with}
&& R(u_{1})\leq \frac{1}{4}(\mathcal{E}_{\tau}(u_{1})^{1/2}+
(\mathcal{E}_{H,min}+\mathcal{O}(\tau_{x}^{6\alpha}))^{1/2})\tau_{x}^{1-3\alpha}\,.
\end{eqnarray*}
With the estimate $\frac{\sqrt{2}}{2}\leq \mathcal{E}_{H,\min}\leq C$,
we deduce
$$
\mathcal{E}_{\tau}(u_{1})= \mathcal{E}_{H,min}+
\mathcal{O}(\tau_{x}^{6\alpha}+ \tau_{x}^{1-3\alpha})\,.
$$
It's time to fix $\alpha$ to the value $\frac{1}{9}$ so that
$\tau_{x}^{6\alpha}=\tau_{x}^{1-3\alpha}
=\tau_{x}^{2/3}$
 and
$$
\inf_{u\in H^{1}(\rz^{2})}\mathcal{E}_{\tau}(u)\leq
\mathcal{E}_{\tau}(u_{1})\leq \mathcal{E}_{H,min}+ \kappa \tau_{x}^{2/3}\,.
$$
\noindent\underline{Second step - Existence of a minimizer:} Once the function
$u_{1}\in H^{1}(\rz^{2})$ has been constructed as above, consider
$\tau_{x}<\tau_{0}$ with $\mathcal{E}_{H,min}+\kappa
\tau_{0}^{2/3}\leq \frac{1}{\ell_{V}^{2}\tau_{0}}$\,.
The functional
$$
\mathcal{E}_{\tau}(u)-\frac{1}{\ell_{V}^{2}\tau_{x}}\|u\|_{L^{2}}^{2}\,,
$$
is the sum of a convex strongly continuous functional (and therefore
weakly continuous) on
$H^{1}(\rz^{2})$ and a negative functional
$$
\langle u, \frac{1}{\ell_{V}^{2}\tau_{x}}\left[v(\tau_{x}^{1/2}.)-1\right]_{-}u\rangle\,.
$$
Due to the compact support of $v-1$, it is also continuous w.r.t the
weak topology on $H^{1}(\rz^{2})$\,.  Out a  minimizing sequence
$(u_{n})_{n\in \nz^{*}}$,
extract a weakly converging subsequence in $H^{1}(\rz^{2})$\,. The
weak limit, $u_{\infty}$,  satisfies
$$
\mathcal{E}_{\tau}(u_{\infty})-
\frac{1}{\ell_{V}^{2}\tau_{x}}\|u_{\infty}\|_{L^{2}}^{2}=\lim_{k\to\infty}\mathcal{E}_{\tau}(u_{n_{k}})
-
\frac{1}{\ell_{V}^{2}\tau_{x}}\|u_{n_{k}}\|_{L^{2}}^{2}\,.
$$
with $\|u_{\infty}\|_{L^{2}}\leq 1$\,.
The same convergence holds also for the energy
$\mathcal{E}_{\tau}(u)-\frac{1}{2\ell_{V}^{2}\tau_{x}}\|u\|_{L^{2}}^{2}$, so
that actually $\|u_{\infty}\|_{L^{2}}=\lim_{k\to\infty}\|u_{n_{k}}\|_{L^{2}}=1$
and $u_{\infty}$ realizes the minimum of $\mathcal{E}_{\tau}(u)$ under
the constraint $\|u\|_{L^{2}}=1$\,.\\
The Euler-Lagrange equation can thus be written, with the stated
straightforward consequences.\\
\noindent\underline{Third step - a priori estimate for minimizers of
  $\mathcal{E}_{\tau}$:}\\
Let $u\in H^{1}(\rz^{2})$ satisfy $\|u\|_{L^{2}}=1$ and
$\mathcal{E}_{\tau}(u)=\mathcal{E}_{\tau,min}\leq \mathcal{E}_{H}+\kappa \tau_{x}^{2/3}$\,.
Take two pairs $\tilde{\chi}=(\tilde{\chi}_{1},\tilde{\chi}_{2})$ and
$\chi'=(\chi'_{1},\chi'_{2})$, like
in our statement, and such that $\chi_{1}'\prec \tilde{\chi}_{1}$. The
identities \eqref{eq.IMSR} for $\mathcal{E}_{\tau}$ and  \eqref{eq.IMSR} in
Lemma~\ref{le.compred} provide
\begin{eqnarray*}
&&\mathcal{E}_{\tau,min}\geq \mathcal{E}_{\tau}(\tilde{\chi}_{1}(\tau_{x}^{1/9}.)u)-C_{\tilde{\chi}}\tau_{x}^{2/9}
\\
&&
\mathcal{E}_{\tau,min}
\geq \mathcal{E}_{H}(\tilde{\chi}_{1}(\tau_{x}^{1/9}.)u)-C_{\tilde{\chi}}\tau_{x}^{2/9}
\\
&&\qquad-\frac{1}{4}\left[\mathcal{E}_{\tau}(\tilde{\chi}_{1}(\tau_{x}^{1/9}.)u)^{1/2}+
\mathcal{E}_{H}(\tilde{\chi}_{1}(\tau_{x}^{1/9}.)u)^{1/2}\right]\tau_{x}^{2/3}\,.
\end{eqnarray*}
The first line says
$$
\mathcal{E}_{\tau}(\tilde{\chi}_{1}(\tau_{x}^{1/9}.)u)
\leq \mathcal{E}_{\tau,min}+C_{\tilde{\chi}}\tau_{x}^{2/9}\leq
\mathcal{E}_{H,min}+ \kappa \tau_{x}^{2/3}+
C_{\tilde{\chi}}\tau_{x}^{2/9}\leq C_{\tilde{\chi}}'\,,
$$
which combined with the second line provides the uniform estimate
$$
\mathcal{E}_{\tau}(\tilde{\chi}_{1}(\tau_{x}^{1/9}.)u)\leq C''_{\tilde{\chi}}\,.
$$
Therefore $\tilde{\chi}_{1}(\tau_{x}^{1/9}.)u$ is uniformly bounded in
$\mathcal{H}_{1}$ with respect to $\tau_{x}$\,.
Consider now the Euler-Lagrange equation
$$
\left[-\partial_{x}^{2}-(\partial_{y}-\frac{i}{2}\frac{x}{\sqrt{1+\tau_{x}x^{2}}})^{2}
  + \frac{v(\tau_{x}^{1/2}.)}{\ell_{V}^{2}\tau_{x}}+G|u|^{2}\right]u=\lambda_{u}u
$$
and write its local version for $u_{1}'=\chi'_{1}(\tau_{x}^{2/9}.)u$ in the
form
\begin{eqnarray}
\nonumber
  &&
\left[-\partial_{x}^{2}-(\partial_{y}-\frac{ix}{2})^{2}
  + \frac{x^{2}+y^{2}}{\ell_{V}^{2}}\right]u'_{1}
=\lambda_{u}u_{1}'-G|\tilde{u}_{1}|^{2}u'_{1}\\
&&
\label{eq.HHu1}\hspace{6cm}+f_{u}^{1}+f^{2}_{u}\,,
\\
\nonumber
\text{with}
&&f^{1}_{u}=-ix(\frac{1}{\sqrt{1+\tau_{x}x^{2}}}-1)\chi'(\tau_{x}^{1/9}.)\partial_{y}\tilde{u}_{1}
- \frac{x^{2}}{4}\frac{\tau_{x}x^{2}}{1+\tau_{x}x^{2}}u'_{1}\,,\\
\nonumber
\text{and}
&&f_{u}^{2}=
-2\tau_{x}^{1/9}(\nabla \chi_{1}')(\tau_{x}^{1/9}.).\nabla
\tilde{u}_{1}-\tau_{x}^{2/9}(\Delta\chi_{1}')(\tau_{x}^{1/9}.)\tilde{u}_{1}\,,
\end{eqnarray}
after setting
$\tilde{u}_{1}=\tilde{\chi}_{1}(\tau_{x}^{2/9}.)u$\,. Both functions,
$\tilde{u}_{1}$ and therefore $u'_{1}=\chi'_{1}(\tau_{x}^{2/9}.)\tilde{u}_{1}$
are uniformly estimated in $\mathcal{H}_{1}$ and therefore in
$H^{1}(\rz^{2})$\,.
From the embedding $H^{1}(\rz^{2})\subset L^{6}(\rz^{2})$, the term
$G|\tilde{u}_{1}|^{2}u'_{1}$ is uniformly bounded in $L^{2}(\rz^{2})$\,.
For the term $f_{u}^{1}$,  the support condition $\supp
\chi'_{1}(\tau_{x}^{1/9}.)\subset \left\{|x|\leq\tau_{x}^{-1/9}\right\}$ imply
$$
\|f^{1}_{u}\|_{L^{2}}\leq C_{\chi'}(\tau_{x}^{1-1/3}+
\tau_{x}^{1-4/9})\leq C''_{\chi'}\,,
$$
while the estimate
$$
\|f_{u}^{2}\|\leq C_{\chi'}\tau_{x}^{1/9}\leq C_{\chi'}''
$$
is straightforward.
Hence the right-hand side of \eqref{eq.HHu1} is uniformly bounded in
$L^{2}(\rz^{2})$ and we have proved
\begin{equation}
  \label{eq.estim1H2}
\|\chi'_{1}(\tau_{x}^{1/9}.)u\|_{\mathcal{H}_{2}}\leq C^{3}_{\chi'}
\end{equation}
for any good pair of cut-offs $\chi'=(\chi_{1}',\chi_{2}')$\,.
\\
\noindent\underline{Fourth step- accurate comparison of minimal energies:}\\
We already know $\mathcal{E}_{\tau,min}\leq
\mathcal{E}_{H,min}+\kappa\tau_{x}^{2/3}$ and we want to check the
reverse inequality. Consider a minimizer $u$ of $\mathcal{E}_{\tau}$  and
take two pairs of cut-off $\chi=(\chi_{1},\chi_{2})$ and
$\chi'=(\chi_{1}',\chi_{2}')$, such that $\chi_{1}\prec \chi_{1}'$\,.
The identity \eqref{eq.IMSH} for $\mathcal{E}_{\tau}$ and
\eqref{eq.IMSR} of Lemma~\ref{le.compred} used with $\chi$
 imply
\begin{eqnarray}
\nonumber
&&
\mathcal{E}_{\tau,min}\geq
 \mathcal{E}_{\tau}(\chi_{1}(\tau_{x}^{1/9}.)u)+
 \mathcal{E}_{\tau}(\chi_{2}(\tau_{x}^{1/9}.)u)
\\
\label{eq.ERmin1}
&&
\hspace{3cm}
-\tau_{x}^{2/9}\sum_{j=1}^{2}\int_{\rz^{2}}|(\nabla\chi_{j})(\tau_{x}^{1/9}.)|^{2}|u|^{2}\,,
\\
\nonumber
&&
\mathcal{E}_{\tau,min}
\geq \mathcal{E}_{H}(\chi_{1}(\tau_{x}^{1/9}.)u)
\\
&&\label{eq.ERmin2}
\hspace{3cm}-\tau_{x}^{2/9}\sum_{j=1}^{2}\int_{\rz^{2}}|(\nabla\chi_{j})(\tau_{x}^{1/9}.)|^{2}|u|^{2}+R(u)\,.
\end{eqnarray}
After setting
$u_{1}'=\chi_{1}'(\tau_{x}^{1/9}.)u$, we get the bound
$$
\int_{\rz^{2}}|(\nabla \chi_{j})(\tau_{x}^{1/9}.)|^{2}|u|^{2}
=
\int_{\rz^{2}}\frac{|(\nabla
  \chi_{j})(\tau_{x}^{1/9}.)|^{2}}{|q|^{4}}||q|^{2}u'_{1}|^{2}\leq C_{\chi}\tau_{x}^{4/9}\,,
$$
while we already know from the third step the bounds
$\mathcal{E}_{\bullet}(\chi_{1}(\tau_{x}^{1/9}.)u)\leq C_{\chi}$ when
$\bullet$ stands for $\tau$ or $H$, which implies $|R(u)|\leq
C_{\chi}\tau_{x}^{2/3}$\,.
From \eqref{eq.ERmin1}, we deduce, as we did in the first step with
the energy $\mathcal{E}_{H}$,
$$
\|\chi_{2}(\tau_{x}^{1/9}.)u\|_{L^{2}}^{2}\leq
C_{\chi}\tau_{x}^{2/3}\quad
,\quad \|\chi_{1}(\tau_{x}^{1/9}.)u\|_{L^{2}}=1+\mathcal{O}(\tau_{x}^{2/3})\,,
$$
while the second line implies
$$
\mathcal{E}_{\tau,min}\geq
\left(1-C_{\chi}\tau_{x}^{2/3}\right)\mathcal{E}_{H,min}-C_{\chi}'\tau_{x}^{2/3}
\geq \mathcal{E}_{H,min}-C''_{\chi}\tau_{x}^{2/3}\,.
$$
\noindent\underline{Fifth step- accurate comparison of minimizers:}\\
The function
$u_{1}=\|\chi_{1}(\tau_{x}^{1/9}.)u\|_{L^{2}}^{-1}\chi_{1}(\tau_{x}^{1/9}.)u$,
satisfies
$$
\mathcal{E}_{H}(u_{1})\leq \mathcal{E}_{H,min}+ C_{\chi}\tau_{x}^{2/3}
$$
while $\chi_{1}(\tau_{x}^{1/9}.)u$ solves the equation \eqref{eq.HHu1} for some pair
$\tilde{\chi}=(\tilde{\chi}_{1},\tilde{\chi}_{2})$ such that
$\chi_{1}\prec \tilde{\chi}_{1}$ after replacing $(\chi', u_{1}')$
with $(\chi,\chi_{1}(\tau_{x}^{1/9}.)u)$\,. After normalization by
setting $\tilde{u}_{1}=\|\chi(\tau_{x}^{1/9}.)u\|_{L^{2}}^{-1}\chi_{1}(\tau_{x}^{1/9}.)u$ it becomes
\begin{eqnarray}
\label{eq.HHu2}
  &&
H_{\ell_{V},G}u_{1}+G|u_{1}|^{2}u_{1}
=\lambda_{u}u_{1}+f_{u}^{1}+f^{2}_{u}+f^{3}_{u}\,,
\\
\nonumber
&&u_{1}=\|\chi(\tau_{x}^{1/9}.)u\|_{L^{2}}^{-1}\chi_{1}(\tau_{x}^{1/9}.)u\quad,\quad
\tilde{u}_{1}=\|\chi(\tau_{x}^{1/9}.)u\|_{L^{2}}^{-1}\tilde\chi_{1}(\tau_{x}^{1/9}.)u\,,\\
\nonumber
\text{with}
&&f^{1}_{u}=-ix(\frac{1}{\sqrt{1+\tau_{x}x^{2}}}-1)\chi(\tau_{x}^{1/9}.)
\partial_{y}\tilde{u}_{1}
- \frac{x^{2}}{4}\frac{\tau_{x}x^{2}}{1+\tau_{x}x^{2}}u_{1}\,,
\\
\nonumber
&&f_{u}^{2}=
-2\tau_{x}^{1/9}(\nabla \chi_{1})(\tau_{x}^{1/9}.).\nabla
\tilde{u}_{1}-\tau_{x}^{2/9}(\Delta\chi_{1})(\tau_{x}^{1/9}.)\tilde{u}_{1}\,,
\\
\nonumber
\text{and}
&&
\hspace{-0.5cm}f_{u}^{3}=G\|\chi_{1}(\tau_{x}^{1/9}.)u\|_{L^{2}}^{2}(|\tilde{u}_{1}|^{2}-|u_{1}|^{2})u_{1}+
G(1-\|\chi_{1}(\tau_{x}^{1/9}.)u\|_{L^{2}}^{2})|u_{1}|^{2}u_{1}\,.
\end{eqnarray}

 The estimate $\eqref{eq.estim1H2}$,
for any new good pair
$\chi'=(\chi_{1}',\chi_{2}')$ such that
$\chi_{1}\prec\tilde{\chi}_{1}\prec\chi'_{1}$,
implies that the terms $f^{1}_{u}$ and $f_{u}^{2}$ of the right-hand side of \eqref{eq.HHu2}
 have an
$L^{2}$-norm
of order $\tau_{x}^{2/3}$\,. For the third term the estimate
\eqref{eq.estim1H2} also implies that
$$
(|\tilde{u}_{1}|^{2}-|u_{1}|^{2})u_{1}=\left(1-\chi_{1}^{2}(\tau_{x}^{1/9}.)\right)|\tilde{u_{1}}|^{2}u_{1}\,,
$$
has an $L^{2}$-norm of order $\mathcal{O}(\tau_{x}^{2/3})$ (use the
$L^{\infty}$ bound for $|\tilde{u}_{1}|^{2}$ with
$\||q|^{2}u_{1}\|_{L^{2}}\leq C_{\chi}$).
We conclude by applying Proposition~\ref{pr.harmH2}\,.
\enddemonstration{}
We end this section with a comparison property similar to
Proposition~\ref{pr.harmH2}.
\begin{proposition}
\label{pr.appH2}
Let $\ell_{V},G$ be fixed positive numbers and take $u\in
H^{1}(\rz^{2})$ such that
$$
\mathcal{E}_{\tau}(u)\leq \mathcal{E}_{\tau,min}+C_{\ell_{V,G}}\tau_{x}^{2/3}
$$
 and which solves
$$
-\partial_{x}^{2}u-(\partial_{y}-i\frac{x}{\sqrt{1+\tau_{x}x^{2}}})^{2}+
\frac{v(\tau_{x}^{1/2}.)}{\ell_{V}^{2}\tau_{x}}u+G|u|^{2}u=\lambda_{u}u+ r_{u}
$$
with $\lambda_{u}\in \rz$ and $\|r_{u}\|_{L^{2}}\leq 1$\,.
Then for any  pair
$\chi=(\chi_{1},\chi_{2})\in \mathcal{C}^{\infty}_{b}(\rz^{2})$ so
that $\chi_{1}^{2}+\chi_{2}^{2}=1$ with $\supp
\chi_{1}\subset\left\{x^{2}+y^{2}<1\right\}$ and $\chi_{1}\equiv 1$
  in $\left\{x^{2}+y^{2}\leq 1/2\right\}$\,, there exists
$\tau_{\chi,\ell_{V},G}$ and $C_{\chi,\ell_{V}}$ such that
  \begin{eqnarray}
    \label{eq.estimu2a}
&&    \|\chi_{2}(\tau_{x}^{1/9}.)u\|^{2}_{L^{2}}\leq
    C_{\chi,\ell_{V},G}\tau_{x}^{2/3}\\
&&
\label{eq.estimEnu1a}
|\mathcal{E}_{\tau}(\chi_{1}(\tau_{x}^{1/9}.)u)-\mathcal{E}_{H,min}|\leq  C_{\chi,\ell_{V},G}\tau_{x}^{2/3}
\\
\label{eq.estimH2u1a}
&& d_{\mathcal{H}_{2}}(\chi_{1}(\tau_{x}^{1/9}.)u,
\mathrm{Argmin}~\mathcal{E}_{H})\leq
C_{\chi,\ell_{V},G}(\tau_{x}^{2\nu_{\ell_{V,G}}/3}+\|r_{u}\|_{L^{2}})\,,
  \end{eqnarray}
when $\tau_{x}< \tau_{\chi,\ell_{V},G}$ and
where $\nu_{\ell_{V},G}\in (0,\frac{1}{2}]$ is the exponent given in Proposition~\ref{pr.harmapp}.
\end{proposition}
\begindemonstration{}
The analysis follows essentially the same line as the study of the
minimizers of $\mathcal{E}_{\tau}$ in the proof of
Proposition~\ref{pr.redHR}.
By taking two pairs $\chi'=(\chi'_{1},\chi'_{2})$ and
$\tilde{\chi}=\tilde(\tilde{\chi}_{1},\tilde{\chi}_{2})$ such that
$\chi'_{1}\prec \tilde{\chi_{1}}$, we obtain successively like in the Third~Step
in the proof of Proposition~\ref{pr.redHR}~:
\begin{itemize}
\item $\mathcal{E}_{\tau}(\tilde{\chi}_{1}(\tau_{x}^{1/9}.)u)+
\mathcal{E}_{\tau}(\tilde{\chi}_{2}(\tau_{x}^{1/9}.)u) \leq
C_{\tilde\chi}$\,;
\item $\left[-\partial_{x}^{2}-(\partial_{y}-\frac{ix}{2})^{2}
  + \frac{x^{2}+y^{2}}{\ell_{V}^{2}}\right]u'_{1}
=\lambda_{u}u_{1}'-G|\tilde{u}_{1}|^{2}u'_{1}+f_{u}^{1}+f^{2}_{u}+r_{u}$
where $u_{1}'$,$\tilde{u}_{1}$, $f_{u}^{1,2}$ have the same
expressions as in \eqref{eq.HHu1};
\item $\|\chi'_{1}(\tau_{x}^{1/9}.)u\|_{\mathcal{H}_{2}}\leq C_{\chi'}$
 owing to $\|r_{u}\|_{L^{2}}\leq 1$\,.
\end{itemize}
From the last estimate, the refined comparison of energies like in
the Fourth~Step gives for a pair $\chi=(\chi_{1},\chi_{2})$ such that
$\chi_{1}\prec \chi'_{1}$:
\begin{eqnarray*}
  \|\chi_{2}(\tau_{x}^{1/9}.)u\|_{L^{2}}^{2}\leq C_{\chi}\tau_{x}^{2/3}\quad
  \mathcal{E}_{H}(u_{1})\leq \mathcal{E}_{H,min}+C_{\chi}\tau_{x}^{2/3}\,,
\end{eqnarray*}
with
$u_{1}=\|\chi_{1}(\tau_{x}^{1/9}.)u\|^{-1}\chi_{1}(\tau_{x}^{1/9}.)u$
and $C_{\chi}\tau_{x}^{2/3}\leq 1$ for $\tau_{x}\leq \tau_{\chi}$\,.
The equation~\eqref{eq.HHu2} is replaced by
$$
H_{\ell_{V},G}u_{1}+G|u_{1}|^{2}u_{1}=\lambda_{u}u_{1}+f_{u}^{1}+f_{u}^{2}+f_{u}^{3}+\|\chi_{1}(\tau_{x}^{1/9}.)u\|^{-1}r_{u}
$$
without changing the expressions of $f_{u}^{1,2,3}$\,. Again the
estimate $\|\chi'_{1}(\tau_{x}^{1/9}.)u\|_{\mathcal{H}_{2}}\leq C_{\chi'}$ is used with various cut-offs
$\chi_{1}'$, in order to get
$\|f^{1}_{u}+f^{2}_{u}+f^{3}_{u}\|_{L^{2}}=\mathcal{O}(\tau_{x}^{2/3})$\,. We
conclude with the help of Proposition~\ref{pr.harmH2} applied to
$u_{1}$\,.
\enddemonstration{}

\section{Analysis of the complete minimization problem}
\label{se.compmini}

We consider the complete minimization problem for the energy
$$
\mathcal{E}_{\varepsilon}(\psi)=\langle
\psi,H_{Lin}\psi\rangle+\frac{G_{\varepsilon,\tau}}{2}\int|\psi^{4}|
= \varepsilon^{2+2\delta}\tau_{x}\left[\langle
\psi,\varepsilon^{-2-2\delta}\tau_{x}^{-1}
H_{Lin}\psi\rangle+\frac{G}{2}\int|\psi^{4}|\right]
$$
and compare its solutions to the minimization of the reduced energies
$\mathcal{E}_{\tau}$ and
$\mathcal{E}_{H}$, introduced in the previous sections. We work with
$\tau_{y}=1$, $\tau_{x}\to 0$, $\varepsilon\to 0$, while $\ell_{V},G$
and $\delta\in (0,\delta_{0}]$ are fixed.
The analysis follows the same lines as the proof of Proposition~\ref{pr.redHR}.

\subsection{Upper bound for
  $\inf\left\{\mathcal{E}_{\varepsilon}(\psi), \|\psi\|_{L^{2}}=1\right\}$}
\label{se.uppcomplete}

The potential $V_{\varepsilon}$ is chosen according
to \eqref{eq.defVepsIntro}-\eqref{eq.defvIntro} while $\ell_{V}$ and $G$ are
fixed. The parameter $\tau_{x}$ is assumed to be smaller than
$\tau_{\ell_{V},G}$ so that the minimal energy
$\mathcal{E}_{\tau,min}(\tau_{x})$
of $\mathcal{E}_{\tau}$ is achieved (see Proposition~\ref{pr.redHR})  and
$|\mathcal{E}_{\tau,min}(\tau_{x})-\mathcal{E}_{H,min}|\leq
C_{\ell_{V},G}\tau_{x}^{2/3}$\,.
Moreover Proposition~\ref{pr.redHR} also says that by truncating an
element of $\mathrm{Argmin}~\mathcal{E}_{\tau}$, one can find
$a_{-}\in \mathcal{H}_{2}$ such that
\begin{equation}
  \label{eq.choixtest}
\|a_{-}\|_{L^{2}}=1\,,\quad |\mathcal{E}_{\tau}(a_{-})-\mathcal{E}_{H,min}|\leq
C_{\ell_{V},G}\tau_{x}^{2/3}
\quad\text{and}\quad
\|a_{-}\|_{\mathcal{H}_{2}}\leq C_{\ell_{V},G}\,.
\end{equation}
\begin{proposition}
\label{pr.bornesup}
  Under the above assumptions, take $\psi=
  \hat{U}\begin{pmatrix}
    0\\e^{-i\frac{y}{2\sqrt{\tau_{x}}}}a_{-}
  \end{pmatrix}$ where $a_{-}$ satisfies \eqref{eq.choixtest} and
  $\hat{U}=U(q,\varepsilon D_{q},\tau,\varepsilon)$ is the unitary
  operator introduced in Theorem~\ref{th.BornOp}.
The estimate
\begin{equation}
\label{bornesup}
\left|\mathcal{E}_{\varepsilon}(\psi)- \varepsilon^{2+2\delta}\tau_{x}\mathcal{E}_{\tau}(a_{-})
\right|\leq C_{\ell_{V},G}\varepsilon^{2+4\delta}\,.
\end{equation}
hold uniformly w.r.t $\tau_{x}\in (0,\tau_{\ell_{V},G}]$ and
$\delta\in (0,\delta_{0}]$\,.
\end{proposition}
\begindemonstration
Let us compare first the linear part by estimating
\begin{multline*}
\left|\langle \psi\,,
H_{Lin}\psi\rangle-\varepsilon^{2+2\delta}\tau_{x}\langle a_{-}\,,
\left[-\partial_{x}^{2}-(\partial_{y}-i
  \frac{x}{2\sqrt{1+\tau_{x}x^{2}}})^{2}+\frac{v(\sqrt{\tau_{x}}.)}{\ell_{V}^{2}\tau_{x}}\right]a_{-}\rangle\right|
\\
=\left|
\langle \psi\,,
(H_{Lin}-\varepsilon^{2+2\delta}\tau_{x}U^{*}H_{BO}U )\psi\rangle
\right|
\end{multline*}
By Proposition~\ref{pr.adiabexpl2}, it suffices to estimate
$$
\|(1-\hat{\chi})\psi\|_{L^{2}}\quad\text{and}\quad
\|\sqrt{\tau_{x}}|\varepsilon D_{q}|(1-\hat{\chi})\psi\|_{L^{2}}
$$
for some given cut-off function $\chi\in
\mathcal{C}^{\infty}_{0}(\rz)$ with
$\hat{\chi}=\chi(\tau_{x}\varepsilon^{2}D_{q}^{2})$\,.
Notice
$$
e^{i\frac{y}{2\sqrt{\tau_{x}}}}(\sqrt{\tau_{x}}\varepsilon D_{q})e^{-i\frac{y}{2\sqrt{\tau_{x}}}}
=
\begin{pmatrix}
  \sqrt{\tau_{x}}\varepsilon D_{x}\\
   \sqrt{\tau_{x}}\varepsilon D_{y}-\frac{\varepsilon}{2}
\end{pmatrix}
$$
Hence by using the functional calculus of $\sqrt{\varepsilon}D_{q}$,
we can say that there exists a cut-off $\chi'\in
\mathcal{C}^{\infty}_{0}(-r_{\gamma}^{2},r_{\gamma}^{2})$, with
$\chi'\equiv 1$ around $0$ and $\chi'\prec\chi$ such that
\begin{eqnarray*}
  &&
e^{i\frac{y}{2\sqrt{\tau_{x}}}}(1-\hat{\chi})^{2}e^{-i\frac{y}{2\sqrt{\tau_{x}}}}\leq
(1-\hat{\chi'})^{2}\,,
\\
\text{and}
&&
e^{i\frac{y}{2\sqrt{\tau_{x}}}}(\tau_{x}\varepsilon^{2}|D_{q}|^{2})(1-\hat{\chi})^{2}e^{-i\frac{y}{2\sqrt{\tau_{x}}}}
\leq 2(\tau_{x}\varepsilon^{2}|D_{q}|^{2})(1-\hat{\chi'})^{2}+ 2(1-\hat{\chi'})^{2}
\end{eqnarray*}
as soon as $\varepsilon\leq \varepsilon_{0}\leq 1$, for a convenient choice
of $\varepsilon_{0}$ and $r_{\gamma}$\,.
By using $\|a\|_{H^{2}(\rz^{2})}\leq C_{\ell_{V},G}$, we deduce
$$
\|(1-\hat{\chi})\psi\|_{L^{2}}^{2}\leq
\|(1-\hat{\chi'})a_{-}\|_{L^{2}}^{2}\leq
C_{\ell_{V},G}\tau_{x}^{2}\varepsilon^{4}\,,
$$
and
\begin{eqnarray*}
\|(\sqrt{\tau_{x}}\varepsilon|D_{q}|)(1-\hat{\chi})\psi\|_{L^{2}}^{2}
&\leq &
2\|(\sqrt{\tau_{x}}\varepsilon|D_{q}|)(1-\hat{\chi'})a_{-}\|_{L^{2}}^{2}
+2\|(1-\hat{\chi'})a_{-}\|_{L^{2}}^{2}\\
&\leq& C_{\ell_{V},G}\tau_{x}\varepsilon^{2}\,.
\end{eqnarray*}
By Proposition~\ref{pr.adiabexpl2}, we obtain
\begin{multline}
  \label{eq.erreurenerlin}
  \left|
\langle \psi\,,
(H_{Lin}-\varepsilon^{2+2\delta}\tau_{x}U^{*}H_{BO}U )\psi\rangle
\right|\\
\leq C_{\ell_{V},G}\left[\varepsilon^{2+4\delta}+
  \varepsilon^{5+2\delta}\tau_{x}^{2}+
  \varepsilon^{4+2\delta}\tau_{x}^{3/2}\right]
\leq C_{\ell_{V},G}'\varepsilon^{2+4\delta}\,.
\end{multline}
For the nonlinear part of the energy,
Proposition~\ref{pr.L4err} gives
$$
(1-C\varepsilon^{1+2\delta})\int_{\rz^{2}}|a_{-}|^{4}\leq
\int_{\rz^{2}}|\psi|^{4}\leq (1+C\varepsilon^{1+2\delta})\int_{\rz^{2}}|a_{-}|^{4}
$$
and the bound,
$\int_{\rz^{2}}|e^{-i\frac{y}{2\sqrt{\tau_{x}}}}a_{-}|^{4}=\int_{\rz^{2}}|a_{-}|^{4}\leq
C_{\ell_{V},G}$, leads to
$$
\left|\frac{G\tau_{x}\varepsilon^{2+2\delta}}{2}
\int_{\rz^{2}}|\psi|^{4}-\varepsilon^{2+2\delta}\frac{G\tau_{x}}{2}\int_{\rz^{2}}|a_{-}|^{4}\right|
\leq
C_{\ell_{V},G}\varepsilon^{2+2\delta}\tau_{x}\times \varepsilon^{1+2\delta}\,,
$$
which is smaller than the error term for the linear part.
This ends the proof of \eqref{bornesup}\,.
\enddemonstration
\begin{remarque}
\label{re.epstau}
\begin{itemize}
\item The energy
  $\mathcal{E}_{\tau}(a_{-})=\mathcal{E}_{H,min}+\mathcal{O}(\tau_{x}^{2/3})$\,. Therefore
  the error given by \eqref{bornesup} is relevant, as compared with
  the energy scale of
  $\varepsilon^{2+2\delta}\tau_{x}\mathcal{E}_{\tau,min}$,  when
  \begin{equation}
    \label{eq.compepstaurel1}
    \varepsilon^{2\delta}\leq c_{\ell_{V},G}\tau_{x}\,,
  \end{equation}
with $c_{\ell_{V},G}$ small enough,
and accurate when
  \begin{equation}
    \label{eq.compepstaurel2}
    \varepsilon^{2\delta}\leq C_{\ell_{V},G}\tau_{x}^{5/3}\,.
  \end{equation}
Remember that  the constants $c_{\ell_{V},G}$ and
 $C_{\ell_{V},G}$ depend also on $\delta_{0}$, when
$\delta\in (0,\delta_{0}]$\,.
\item   It is interesting to notice that the worst term in the right-hand
  side of \eqref{eq.erreurenerlin} comes from the error of order
  $\mathcal{O}(\varepsilon^{2+2\delta})$ in the Born-Oppenheimer
  approximation. There seems to be no way to get an additional factor
  $\tau_{x}^{\alpha}$ with $\alpha>0$ because the initial problem is
  rapidly oscillatory in the $y$-variable in a $\tau_{x}$-dependent
  scale. This can be seen on the gain associated with
  the metric $g_{\tau}$, for $\tau_{y}=1$ and $\tau_{x}>0$, which is simply
  $\langle \sqrt{\tau}_{x}p\rangle$ or essentially $1$ when $p$ is small.
\end{itemize}

\end{remarque}

\subsection{Existence of a minimizer for $\mathcal{E}_{\varepsilon}$}
\label{se.exisEeps}
With the choice \eqref{eq.defVepsIntro}-\eqref{eq.defvIntro} of the potential
$V_{\varepsilon,\tau}$, the linear Hamiltonian $H_{Lin}$ can be
written
\begin{eqnarray*}
  &&
H_{Lin}=-\varepsilon^{2+2\delta}\tau_{x}\Delta + u_{0}(q,\tau)
\begin{pmatrix}
  2\sqrt{1+\tau_{x}x^{2}}&0\\
0&0
\end{pmatrix}u_{0}(q,\tau)^{*}\\
&&\hspace{5cm}+\varepsilon^{2+2\delta}\frac{v(\tau_{x}^{1/2}.)}{\ell_{V}^{2}}
-\varepsilon^{2+2\delta}W_{\tau}(x,y)\,,\\
\text{with}
&&W_{\tau}(x,y)=\frac{\tau_{x}^{2}}{(1+\tau_{x}x^{2})^{2}}+
\frac{1}{1+\tau_{x}x^{2}}\,.
\end{eqnarray*}
For $t\leq \frac{\varepsilon^{2+2\delta}}{2\ell_{V}^{2}}$, the
negative part
$\left(\varepsilon^{2+2\delta}\frac{v(\sqrt{\tau_{x}}.)}{\ell_{V}^{2}}-t\right)_{-}$ is compactly supported. Set
\begin{multline*}
H_{Lin,t,+}=-\varepsilon^{2+2\delta}\tau_{x}\Delta + u_{0}(q,\tau)
\begin{pmatrix}
  2\sqrt{1+\tau_{x}x^{2}}&0\\
0&0
\end{pmatrix}u_{0}(q,\tau)^{*}
\\
+
\left(\varepsilon^{2+2\delta}\frac{v(\sqrt{\tau_{x}}.)}{\ell_{V}^{2}}-t\right)_{+}\,,
\end{multline*}
so that
\begin{multline}
  \label{eq.decompEeps}
\mathcal{E}_{\varepsilon}(\psi)-t\|\psi\|^{2}=\langle
\psi\,,H_{Lin,t,+}\psi\rangle
+\frac{G_{\varepsilon,\tau}}{2}\int|\psi|^{4}
\\
+ \langle
\psi\,,\,\varepsilon^{2+2\delta}\left[\left(\frac{v(\sqrt{\tau_{x}}.)}{\ell_{V}^{2}}-\frac{t}{\varepsilon^{2+2\delta}}\right)_{-}
-W_{\tau}\right]\psi\rangle\,.
\end{multline}
\begin{proposition}
  \label{pr.existEeps}
Assume $\varepsilon\leq\varepsilon_{\ell_{V},G}$ and $\tau_{x}\leq
\tau_{\ell_{V},G}$  with $\tau_{\ell_{V},G}$ small enough. Then the infimum
$$
\inf\{\mathcal{E}_{\varepsilon}(u), \quad u\in H^{1}(\rz^{2}), \|u\|_{L^{2}}=1\}
$$
 is achieved. Any element $\psi$ of
 $\mathrm{Argmin}~\mathcal{E}_{\varepsilon}$ solves an Euler-Lagrange
 equation
$$
H_{Lin}\psi+G_{\varepsilon,\tau}|\psi|^{2}\psi=\lambda_{\psi}\psi
$$
with the estimates
\begin{eqnarray*}
  &&
|\mathcal{E}_{\varepsilon}(\psi)|+|\lambda_{\psi}|\leq C_{\ell_{V},G}\varepsilon^{2+2\delta}
\quad,\quad
\|(1+\tau_{x}|D_{q}|^{2})^{1/2}\psi\|_{L^{2}}\leq C_{\ell_{V},G}\,,
\\
&&
\begin{array}[c]{ll}
\mathcal{E}_{\varepsilon}(\psi)=\mathcal{E}_{\varepsilon,min}&\leq
\varepsilon^{2+2\delta}\left[\mathcal{E}_{\tau,min}+C_{\ell_{V},G}\varepsilon^{2\delta})\right]
\\
&\leq
\varepsilon^{2+2\delta}\left[\tau_{x}\mathcal{E}_{H,min}+
C_{\ell_{V},G}'(\tau_{x}^{5/3}+\varepsilon^{2\delta})\right]\,.
\end{array}
\end{eqnarray*}
\end{proposition}
\begindemonstration{}
From Proposition~\ref{pr.bornesup}, we know that
$$
\inf_{\|u\|_{L^{2}}=1}\mathcal{E}_{\varepsilon}(u)\leq
\varepsilon^{2+2\delta}\tau_{x}\mathcal{E}_{\tau,min}+
C_{\ell_{V},G}\varepsilon^{2+4\delta}=:\frac{t}{2}\,.
$$
For $\varepsilon\leq \varepsilon_{\ell_{V},G}$ and $\tau_{x}\leq
\tau_{\ell_{V},G}$, $t$ smaller than
$\frac{\varepsilon^{2+2\delta}}{2\ell_{V}^{2}}$\,. Consider the
decomposition \eqref{eq.decompEeps} for the energy
$\mathcal{E}_{\varepsilon}(u)-t\|u\|_{L^{2}}^{2}$\,.
By the same argument (convexity  of the positive part and compactness of
the negative part) as we used for $\mathcal{E}_{\tau}$ in the proof of
Proposition~\ref{pr.redHR} (second step), a weak limit of an extracted
sequence of minimizers in $H^{1}(\rz^{2})$ is a minimum for
$\mathcal{E}_{\varepsilon}$ on $\left\{\|u\|_{L^{2}}=1\right\}$\,.\\
A element $\psi$ of $\mathrm{Argmin}~\mathcal{E}_{\varepsilon}$ satisfies
\begin{eqnarray*}
\varepsilon^{2+2\delta}\tau_{x}\left[\langle \psi\,,\, -\Delta
\psi\rangle+\frac{G}{2}\int|\psi|^{4}\right]
&&\leq \langle \psi\,,\, H_{Lin,t, +}\psi\rangle
+\frac{G_{\varepsilon,\tau}}{2}\int|\psi|^{4}
\\
&&\hspace{-2cm}\leq \mathcal{E}_{\varepsilon,min}-\langle
\psi\,,\,\varepsilon^{2+2\delta}\left[\left(\frac{v(\sqrt{\tau_{x}}.)}{\ell_{V}^{2}}-\frac{t}{\varepsilon^{2+2\delta}}\right)_{-}
-W_{\tau}\right]\psi\rangle
\\
&&\leq
\frac{t}{2}+t+
(1+\tau_{x}^{2})\varepsilon^{2+2\delta}
\leq C_{\ell_{V},G}'\varepsilon^{2+2\delta}\,,
\end{eqnarray*}
by recalling $v\geq 0$ for the last line.
This implies
$$
\|\psi\|_{H^{1}}^{2}\leq
C_{\ell_{V},G}'\tau_{x}^{-1}\quad,\quad
\|\psi\|_{L^{4}}^{4}\leq \frac{2C_{\ell_{V},G}'}{G}\tau_{x}^{-1}\,,
$$
and by interpolation with $\|u\|_{L^{6}}\leq C\|\nabla
u\|_{L^{2}}^{2/3}\|u\|_{L^{2}}^{1/3}$,
$\|\psi\|_{L^{6}}=\mathcal{O}(\tau_{x}^{-1/3})$\,. The first
inequality with $\|\psi\|_{L^{2}}=1$, gives
$\|(1+\tau_{x}|D_{q}|^{2})^{1/2}\psi\|_{L^{2}}=\mathcal{O}(1)$\,.\\
The Euler-Lagrange equation
$$
H_{Lin}\psi+G_{\varepsilon,\tau}|\psi|^{2}\psi=\lambda\psi
$$
implies
\begin{eqnarray*}
|\lambda|&\leq &
2\left[\langle \psi\,,\, H_{Lin,t, +}\psi\rangle
+\frac{G_{\varepsilon,\tau}}{2}\int|\psi|^{4}\right]
\\
&&-\langle
\psi\,,\,\varepsilon^{2+2\delta}\left[\left(\frac{v(\sqrt{\tau_{x}}.)}{\ell_{V}^{2}}-\frac{t}{\varepsilon^{2+2\delta}}\right)_{-}
-W_{\tau}\right]\psi\rangle
\leq C''_{\ell_{V},G}\varepsilon^{2+2\delta}\,.
\end{eqnarray*}
Similarly, the lower bound $\mathcal{E}_{\varepsilon}(\psi)\geq
-2\varepsilon^{2+2\delta}$ is due to  $W_{\tau}\geq -2$\,.
The upper bound of
$\mathcal{E}_{\varepsilon}(\psi)=\mathcal{E}_{\varepsilon,min}$ comes
from Proposition~\ref{pr.bornesup} and \eqref{eq.choixtest}.
\enddemonstration{}

\subsection{Comparison of minimal energies between
  $\mathcal{E}_{\varepsilon}$ and $\mathcal{E}_{\tau}$}
\label{se.accminen}
In this subsection, we specify a priori estimates for the minimizers of
$\mathcal{E}_{\varepsilon}$ and compare the energies
$\mathcal{E}_{\varepsilon,min}$ and $\mathcal{E}_{\tau,min}$ without
imposing relations between $\varepsilon^{2\delta}$ and
$\tau_{x}$\,. This is not necessary at this level, if one uses
carefully bootstrap arguments.
\begin{proposition}
\label{pr.compen}
Let $V_{\varepsilon,\tau}$ be given by
\eqref{eq.defVepsIntro}-\eqref{eq.defvIntro} and assume $\tau_{x}\leq
\tau_{\ell_{V},G}$ and $\varepsilon\leq \varepsilon_{\ell_{V},G}$ so
that $\mathcal{E}_{\varepsilon}$ admits a ground state according to
Proposition~\ref{pr.existEeps}\,.
The operator $\hat{U}$ is the unitary
transform provided by Theorem~\ref{th.BornOp} and an element
$\psi\in \mathrm{Argmin}~\mathcal{E}_{\varepsilon}$ is written
$\hat{U}
\begin{pmatrix}
  e^{+i\frac{y}{2\sqrt{\tau_{x}}}}a_{+}\\
e^{-i\frac{y}{2\sqrt{\tau_{x}}}}a_{-}
\end{pmatrix}$\,, with $a=
\begin{pmatrix}
  a_{+}\\a_{-}
\end{pmatrix}
$\,.
Then, the estimates
\begin{eqnarray}
 \label{eq.H1a}
 &&\|(1+\tau_{x}|D_{q}|^{2})^{1/2}a\|_{L^{2}}\leq C_{\ell_{V},G}\,,\\
\label{eq.L2a+}
&&
\|a_{+}\|_{L^{2}}^{2}+G_{\varepsilon,\tau}\int|a|^{4}\leq
  C_{\ell_{V},G}\varepsilon^{2+4\delta}+\mathcal{E}_{\varepsilon,min}\,,\\
\label{eq.ena-}
&&
|\mathcal{E}_{\varepsilon}(\psi)-\varepsilon^{2+2\delta}\tau_{x}\mathcal{E}_{\tau}(a_{-})|\leq
C_{\ell_{V},G}\varepsilon^{2+4\delta}\,,\\
&&
\label{eq.compEmin}
|\mathcal{E}_{\varepsilon,min}-\varepsilon^{2+2\delta}\tau_{x}\mathcal{E}_{\tau,min}|\leq C_{\ell_{V},G}\varepsilon^{2+4\delta}\,,
\end{eqnarray}
hold with right-hand
sides  which can be
replaced
by $C_{\ell_{V},G}\varepsilon^{2+2\delta}\tau_{x}$ when
$\varepsilon^{2\delta}\leq c_{\ell_{V},G}\tau_{x}$\,.
\end{proposition}
\begindemonstration{}
 For $\tilde{a}=
\begin{pmatrix}
  e^{i\frac{y}{2\sqrt{\tau_{x}}}}a_{+}\\
 e^{-i\frac{y}{2\sqrt{\tau_{x}}}}a_{-}
\end{pmatrix}
$ and $a=\begin{pmatrix}
  a_{+}\\
 a_{-}
\end{pmatrix}$, the norms $\|(1+\tau_{x}|D_{q}|^{2})^{1/2}a\|_{L^{2}}$ and
$\|(1+\tau_{x}|D_{q}|^{2})^{1/2}\tilde{a}\|_{L^{2}}$ are uniformly equivalent
because
$$
e^{\pm
  i\frac{y}{2\sqrt{\tau_{x}}}}(\sqrt{\tau_{x}}D_{y})e^{\mp
  i\frac{y}{2\sqrt{\tau_{x}}}}
=(\sqrt{\tau_{x}}D_{y})\mp\frac{1}{2}\,.
$$
 Hence it suffices to
estimate $\sqrt{\tau_{x}}\varepsilon D_{q}\hat{U}^{*}\psi$\,. Remember that
$\hat{U}=U(q,\varepsilon D_{q},\tau,\varepsilon)$ with $U\in
S_{u}(1,g_{\tau};\mathcal{M}_{2}(\cz))$, while $\sqrt{\tau_{x}}p\in
S_{u}(\langle\sqrt{\tau_{x}}p\rangle,g_{\tau};\rz^{2})$\,.
With $\|(1+\tau_{x}
|D_{q}|^{2})^{1/2}\psi\|_{L^{2}}\leq C_{\ell_{V},G}$ in
Proposition~\ref{pr.existEeps}, simply compute:
\begin{equation*}
  \sqrt{\tau_{x}}\varepsilon D_{q}\hat{U}^{*}\psi
=\varepsilon \hat{U}^{*}\sqrt{\tau_{x}} D_{q}\psi +
\left[\sqrt{\tau_{x}}\varepsilon D_{q}, \hat{U}^{*}\right]\psi
=\mathcal{O}(\varepsilon)\quad\text{in}~L^{2}(\rz^{2};\cz^{2})\,,
\end{equation*}
and divide by $\varepsilon$ for \eqref{eq.H1a}.\\
In order to compare the energies, consider first the linear part by writing
\begin{multline*}
\left|\langle \psi\,,
H_{Lin}\psi\rangle-\varepsilon^{2+2\delta}\tau_{x}\left[\langle
  a_{+}\,,\,\hat{H}_{+}a_{+}\rangle+ \langle a_{-}\,,\,\hat{H}_{-}a_{-}\rangle\right]\right|
=
\\
\left|
\langle \tilde{a}\,,
(\hat{U}^{*}H_{Lin}\hat{U}-\varepsilon^{2+2\delta}H_{BO})\tilde{a}\rangle
\right|\,,
\end{multline*}
with
\begin{eqnarray*}
\tilde{a}&=&
\begin{pmatrix}
  e^{i\frac{y}{2\sqrt{\tau_{x}}}}a_{+}\\
 e^{-i\frac{y}{2\sqrt{\tau_{x}}}}a_{-}
\end{pmatrix}\,,
\\
  H_{BO}&=&\begin{pmatrix}
 e^{i\frac{y}{2\sqrt{\tau_{x}}}}\hat{H}_{+}e^{-i\frac{y}{2\sqrt{\tau_{x}}}}&0\\
0&e^{-i\frac{y}{2\sqrt{\tau_{x}}}}\hat{H}_{-}
e^{i\frac{y}{2\sqrt{\tau_{x}}}}
\end{pmatrix}\,,
\\
\hat{H}_{\pm}&=&-\partial_{x}^{2}-(\partial_{y}\pm i
  \frac{x}{2\sqrt{1+\tau_{x}x^{2}}})^{2}+\frac{v(\sqrt{\tau_{x}}.)}{\ell_{V}^{2}\tau_{x}}+\frac{(1\pm
  1)}{\varepsilon^{2+2\delta}\tau_{x}} (1+\tau_{x}x^{2})^{1/2}\,.
\end{eqnarray*}
Here it is convenient to write
\begin{eqnarray}
\nonumber
&&\varepsilon^{2+2\delta}\tau_{x}H_{BO}= \varepsilon^{2\delta}
\hat{B}
+
\begin{pmatrix}
  2(1+\tau_{x}x^{2})^{1/2}\\
0
\end{pmatrix}\,,
\\
\label{eq.defm1}
&&
\hat{B}= \begin{pmatrix}
\hat{B}_{+}& 0\\
0&\hat{B}_{-}
\end{pmatrix}\quad,\quad
\hat{B}_{\pm}=
B_{\pm}(q,\varepsilon D_{q},\tau,\varepsilon)\,,\\
\label{eq.defm2}
&&
\hspace{-2cm}
B_{\pm}(q,p,\tau,\varepsilon)=\tau_{x}p_{x}^{2}+
(\sqrt{\tau_{x}}p_{y}\pm
\frac{\varepsilon}{2}(\frac{\sqrt{\tau_{x}}x}{\sqrt{1+\tau_{x}x^{2}}}
-\sqrt{\tau_{x}}))^{2}
+\frac{\varepsilon^{2}v(\sqrt{\tau_{x}}.)}{\ell_{V}^{2}}\,.
\end{eqnarray}
With the notation
$\hat{\chi}=\chi(\tau_{x}\varepsilon^{2}|D_{q}|^{2})$ of
Proposition~\ref{pr.adiabexpl2}, Lemma~\ref{le.ellipt} below says in
particular
\begin{eqnarray*}
  &&
\|(1-\hat{\chi})\psi\|_{L^{2}}^{2}\leq C_{\ell_{V},G}\left[\langle
  \tilde{a}, \hat{B}\tilde{a}\rangle + \varepsilon^{1+2\delta}\|\tilde{a}\|_{L^{2}}^{2}\right]\,,
\\
&&
\|\sqrt{\tau_{x}}\varepsilon|D_{q}|(1-\hat{\chi})\psi\|_{L^{2}}^{2}
\leq C_{\ell_{V},G}\left[\langle
  \tilde{a}, \hat{B}\tilde{a}\rangle +
  \varepsilon^{1+2\delta}\|\tilde{a}\|_{L^{2}}^{2}\right]\,.
\end{eqnarray*}
Proposition~\ref{pr.adiabexpl2} yields
$$
\left|
\langle \tilde{a}\,,
\left(\hat{U}^{*}H_{Lin}\hat{U}-\varepsilon^{2\delta}\hat{B}-
\begin{pmatrix}
  2(1+\tau_{x}x^{2})^{1/2}\\
0
\end{pmatrix}\right)\tilde{a}\rangle
\right|\leq
C_{\ell_{V},G}\left[\varepsilon^{2+4\delta}+\varepsilon^{1+2\delta}\langle\tilde{a}\,,\,
\hat{B}\tilde{a}\rangle\right]\,.
$$
For the nonlinear part of the energy,
Proposition~\ref{pr.L4err} gives
$$
(1-C\varepsilon^{1+2\delta})\int_{\rz^{2}}|\psi|^{4}\leq
\int_{\rz^{2}}|\tilde{a}|^{4}=\int_{\rz^{2}}|a|^{4}\leq (1+C\varepsilon^{1+2\delta})\int_{\rz^{2}}|\psi|^{4}
$$
with $\tilde{a}=
\begin{pmatrix}
  e^{i\frac{y}{2\sqrt{\tau_{x}}}}a_{+}\\
 e^{-i\frac{y}{2\sqrt{\tau_{x}}}}a_{-}
\end{pmatrix}
$ and $a=\begin{pmatrix}
  a_{+}\\
 a_{-}
\end{pmatrix}$\,.
The bound
$\int_{\rz^{2}}|\psi|^{4}\leq
C_{\ell_{V},G}\tau_{x}^{-1}$ coming from
$\mathcal{E}_{\varepsilon}(\psi)=\mathcal{O}(\varepsilon^{2+2\delta})$
with $\varepsilon^{2+2\delta}W_{\tau}=\mathcal{O}(\varepsilon^{2+2\delta})$, leads to
$$
\left|\frac{G\tau_{x}\varepsilon^{2+2\delta}}{2}
\int_{\rz^{2}}|\psi|^{4}-\frac{G\tau_{x}\varepsilon^{2+2\delta}}{2}\int_{\rz^{2}}|\tilde{a}|^{4}\right|
\leq
C_{\ell_{V},G}\varepsilon^{2+2\delta}\tau_{x}\times
\varepsilon^{1+2\delta}\tau_{x}^{-1}
=C_{\ell_{V},G}\varepsilon^{3+4\delta}\,,
$$
which is smaller than the error term for the linear part.
We have proved
\begin{multline*}
 \left|\mathcal{E}_{\varepsilon}(\psi)-\varepsilon^{2\delta}
\langle \tilde{a}\,,\, \hat{B}\tilde{a}\rangle
-
\langle \tilde{a}_{+}\,,\, 2(1+\tau_{x}x^{2})^{1/2}\tilde{a}_{+}\rangle
-\frac{G_{\varepsilon,\tau}}{2}\int |a|^{4}
\right|
\\
\leq C_{\ell_{V},G}\left[\varepsilon^{2+4\delta}+ \varepsilon \times
  \varepsilon^{2\delta}
\langle \tilde{a}\,,\, \hat{B}\tilde{a}\rangle\right]\,.
\end{multline*}
With
$\mathcal{E}_{\varepsilon}(\psi)=\mathcal{E}_{\varepsilon,min}$
($=\mathcal{O}(\varepsilon^{2+2\delta})$), this gives
$$
(1-C\varepsilon)\varepsilon^{2\delta}\langle
  \tilde{a}\,,\,\hat{B}\tilde{a}\rangle+
2\langle
  \tilde{a}_{+}\,,\,(1+\tau_{x}x^{2})^{1/2}\tilde{a}_{+}\rangle+\frac{G_{\varepsilon,\tau}}{2}\int|a|^{4}
\leq
C\varepsilon^{2+4\delta}+\mathcal{E}_{\varepsilon,min}\,,
$$
where all the terms of the left-hand side are now non negative. We
deduce \eqref{eq.L2a+} and by bootstrapping
\begin{multline*}
\varepsilon^{2+2\delta}\tau_{x}\langle a_{-}\,,\,
\hat{H}_{-}a_{-}\rangle +\frac{G_{\varepsilon,\tau}}{2}\int |a|^{4}
\leq \varepsilon^{2\delta}\langle
\tilde{a}\,,\,\hat{B}\tilde{a}\rangle
 +\frac{G_{\varepsilon,\tau}}{2}\int |\tilde{a}|^{4}
\\
\leq \mathcal{E}_{\varepsilon,min}+C_{\ell_{V},G}\varepsilon^{2+4\delta}\,.
\end{multline*}
Additionally,
$\int|a|^{4}=\int\left(|a_{-}|^{2}+|a_{+}|^{2}\right)^{2}\geq
\int|a_{-}|^{4}$ also gives
$$
0\leq \|a_{-}\|_{L^{2}}^{4}\mathcal{E}_{\tau,min}\leq \mathcal{E}_{\tau}(a_{-})
\leq
C_{\ell_{V},G}\frac{\varepsilon^{2\delta}}{\tau_{x}}+\frac{\mathcal{E}_{\varepsilon,min}}{\varepsilon^{2+2\delta}\tau_{x}}\leq \mathcal{E}_{\tau,min}+C'_{\ell_{V},G}\frac{\varepsilon^{2\delta}}{\tau_{x}}\,.
$$
With
$\|a_{-}\|^{2}_{L^{2}}=1-\|a_{+}\|^{2}=1+\mathcal{O}(\varepsilon^{2+2\delta})=1+\mathcal{O}(\frac{\varepsilon^{2\delta}}{\tau_{x}})$
and $\mathcal{E}_{\tau,min}=\mathcal{E}_{H,min}+\mathcal{O}(\tau_{x}^{2/3})$,
this finally leads to
$$
|\frac{\mathcal{E}_{\varepsilon,min}}{\varepsilon^{2+2\delta}\tau_{x}}
-\mathcal{E}_{\tau}(a_{-})|\leq C''_{\ell_{V},G}\frac{\varepsilon^{2\delta}}{\tau_{x}}\,.
$$
\enddemonstration{}
\begin{lemme}
\label{le.ellipt} Assume $B(q,p,\tau,\varepsilon)=\tau_{x}|p|^{2}+\varepsilon r$ with $r\in
S_{u}(\langle \sqrt{\tau_{x}}p\rangle, g_{\tau};
\mathcal{M}_{2}(\cz))$, take any $\chi\in
\mathcal{C}^{\infty}_{0}(\rz)$ and set $\hat{B}=B(q,\varepsilon
D_{q},\tau,\varepsilon)$ and
$\hat{\chi}=\chi(\tau_{x}\varepsilon|D_{q}|^{2})$\,. Then
 there exists $\varepsilon_{B,\chi}>0$ and
$C_{B,\chi}>0$ such that the estimates
\begin{eqnarray*}
&&
  \|(1-\hat{\chi})u\|_{L^{2}}^{2}+\|\sqrt{\tau_{x}}|\varepsilon
  D_{q}|(1-\hat{\chi})u\|_{L^{2}}^{2}\leq C_{B,\chi}\left[\langle u\,,\,
    \hat{B}u\rangle
+ \varepsilon^{1+2\delta}\|u\|_{L^{2}}^{2}\right]\,,
\\
\text{and}&&
  \|(1-\hat{\chi})u\|_{L^{2}}+\|\sqrt{\tau_{x}}|\varepsilon
  D_{q}|(1-\hat{\chi})u\|_{L^{2}}\leq C_{B,\chi}\left[
    \|\hat{B}u\|_{L^{2}}
+ \varepsilon^{1+2\delta}\|u\|_{L^{2}}\right]\,,
\end{eqnarray*}
hold uniformly w.r.t $\delta\in (0,\delta_{0}]$ and $\varepsilon\in (0,\varepsilon_{B,\chi})$\,.
\end{lemme}
\begindemonstration{}
For $\chi_{0}\in \mathcal{C}^{\infty}_{0}(\rz)$  such that
$\chi_{0}\geq 0$ and $\chi_{0}\equiv 1$ around $0$, the symbol
$\chi_{0}(\tau_{x}p^{2})+B$ is an elliptic symbol in $S_{u}(\langle
\sqrt{\tau_{x}}p\rangle^{2},g_{\tau}:\mathcal{M}_{2}(\cz))$\,. For
$\varepsilon>0$ small enough according to $(\chi_{0},B)$, its quantization
$\hat{\chi}_{0}+\hat{B}$ is
invertible and its inverse belongs to  $OpS_{u}(\langle
\sqrt{\tau_{x}}p\rangle^{-2},g_{\tau}, \mathcal{M}_{2}(\cz))$\,. Next we
notice
that in
$$
(1-\hat{\chi})=(1-\hat\chi)\circ\left[\hat{\chi}_{0}+\hat{B}\right]^{-1}\circ
\hat{B}
+
(1-\hat\chi)\circ\left[\hat{\chi}_{0}+\hat{B}\right]^{-1}\circ \hat{\chi}_{0}\,,
$$
the last term belongs to $Op\mathcal{N}_{u,g_{\tau}}$ if $\chi_{0}\prec
\chi$\,.
All the estimates are consequences of
$$
(1-\hat{\chi})=(1-\hat\chi)\circ\left[\hat{\chi}_{0}+\hat{B}\right]^{-1}\circ
\hat{B} +\varepsilon^{1+2\delta}\rho\,,
$$
with $\rho\in S_{u}(\frac{1}{\langle \tau_{x}p\rangle^{\infty}},g_{\tau};\mathcal{M}_{2}(\cz))$\,.
\enddemonstration{}
\subsection{Adiabatic Euler-Lagrange equation}
\label{se.adELeq}

As suggested by Proposition~\ref{pr.appH2}, or the last steps in the
proof of Proposition~\ref{pr.redHR}, an accurate comparison of
minimizers requires some comparison of the Euler-Lagrange equations.
We check here that the $a_{-}$ component in
Proposition~\ref{pr.compen} solves approximately the Euler-Lagrange
equations for minimizers of $\mathcal{E}_{\tau}$\,. Here the bootstrap
argument is made in terms of operators instead of quadratic forms.
In order to get reliable results, we now assume
$$
\varepsilon^{2\delta}\leq c_{\ell_{V},G}\tau_{x}\,,
$$
with $c_{\ell_{V},G}$ chosen small enough.\\
With such an assumption we know:
\begin{itemize}
\item from \eqref{eq.compEmin} and
  $|\mathcal{E}_{\tau,min}-\mathcal{E}_{H,min}|=\mathcal{O}(\tau_{x}^{2/3})$,
  $$
C^{-1}\varepsilon^{2+2\delta}\tau_{x}\leq
\mathcal{E}_{\varepsilon,min}\leq C\varepsilon^{2+2\delta}\tau_{x}\,.
$$
\item When $\psi\in \mathrm{Argmin}~\mathcal{E}_{\varepsilon,min}$ is
  written $\psi=\hat{U}
  \begin{pmatrix}
    e^{i\frac{y}{2\sqrt{\tau_{x}}}}a_{+}\\
 e^{-\frac{y}{2\sqrt{\tau_{x}}}}a_{+}
  \end{pmatrix}
$\, the $L^{4}$-norm of $a$ is uniformly bounded,  $\int|a|^{4}\leq
C_{\ell_{V},G}$,
according to  \eqref{eq.L2a+}.
By applying Lemma~\ref{le.Taylor} this also gives $\int|\psi|^{4}\leq
C_{\ell_{V},G}$\,.
We also have  $\|a_{+}\|_{L^{2}}^{2}\leq C_{\ell_{V},G}\varepsilon^{2+2\delta}\tau_{x}$\,.
\item The Lagrange multiplier $\lambda_{\psi}$,
  associated with $\psi\in
  \mathrm{Argmin}~\mathcal{E}_{\varepsilon,min}$, equals
  $\mathcal{E}_{\varepsilon}(\psi)+ \frac{G_{\varepsilon,\tau}}{2}\int
  |\psi|^{4}$\,. Thus it is of order $\mathcal{O}(\varepsilon^{2+2\delta}\tau_{x})$\,.
\end{itemize}
\begin{proposition}
\label{pr.adEL}
Under the same assumptions as in Proposition~\ref{pr.compen} and with
the above condition $\varepsilon^{2\delta}\leq
c_{\ell_{V},G}\tau_{x}$, write a
minimizer $\psi$ of $\mathcal{E}_{\varepsilon}$ in the form $\psi=\hat{U}
\begin{pmatrix}
  e^{i\frac{y}{2\sqrt{\tau_{x}}}}a_{+}\\
  e^{-i\frac{y}{2\sqrt{\tau_{x}}}}a_{-}
\end{pmatrix}$\,.
Then the component $a_{-}$ solves the equation
$$
\hat{H}_{-}a_{-}+G\int|a_{-}|^{2}a_{-}=\lambda_{\psi}a_{-}+r_{\varepsilon}\quad
\text{with}\quad \|r_{\varepsilon}\|_{L^{2}}\leq
C_{\ell_{V},G}(\varepsilon+\frac{\varepsilon^{2\delta}}{\tau_{x}})\,,
$$
while $\|a_{+}\|_{L^{2}}\leq
C_{\ell_{V},G}\varepsilon^{2+2\delta}\tau_{x}$\,.\\
The $L^{p}$-norms of $a$ and $\psi$ are uniformly bounded by $C_{\ell_{V},G}$
for $p\in\left[2,6\right]$\,.\\
Moreover if $\hat{B}$ is the operator defined in
\eqref{eq.defm1}-\eqref{eq.defm2}, $\tilde{a}=\begin{pmatrix}
  e^{i\frac{y}{2\sqrt{\tau_{x}}}}a_{+}\\
  e^{-i\frac{y}{2\sqrt{\tau_{x}}}}a_{-}
\end{pmatrix}$ satisfies
\begin{equation}
\label{eq.H2a}
\|\hat{B}\tilde{a}\|_{L^{2}}\leq C_{\ell_{V},G}\varepsilon^{2}\tau_{x}\,.
\end{equation}
\end{proposition}
\begin{remarque}
  The relation of $\|\hat{B}\tilde{a}\|_{L^{2}}$ with
  $\|\tilde{a}\|_{H^{2}}$ is given by
$$
C^{-1}
\left[\varepsilon^{2}\tau_{x}\|D_{q}^{2}\tilde{a}\|+ \varepsilon^{2}\|\tilde{a}\|_{L^{2}}\right]
\leq
\|\hat{B}\tilde{a}\|_{L^{2}}+ \varepsilon^{2}\|\tilde{a}\|_{L^{2}} \leq C
\left[\varepsilon^{2}\tau_{x}\|D_{q}^{2}\tilde{a}\|+ \varepsilon^{2}\|\tilde{a}\|_{L^{2}}\right]\,.
$$
Note that the $L^{2}$ remainder terms have a factor $\varepsilon^{2}$
and not $\varepsilon^{2}\tau_{x}$\,.
\end{remarque}
\begindemonstration{}
Playing with the Euler-Lagrange equation for $\psi$,
we shall first prove \eqref{eq.H2a}  by using the same argument as we did
for $\langle \tilde{a}\,,\, \hat{B}\tilde{a}\rangle$ in the variational proof
of  Proposition~\ref{pr.compen} and then use it in order
to estimate $\|r\|_{L^{2}}$\,.
 The Euler-Lagrange equation for $\psi$
$$
H_{Lin}\psi+ G_{_{\varepsilon,\tau}}|\psi|^{2}\psi=\lambda_{\psi}\psi
$$
becomes
$$
\hat{U}^{*}H_{Lin}\hat{U}\tilde{a}+G_{\varepsilon, \tau}\hat{U}^{*}\left(|\psi|^{2}\psi\right)=\lambda_{\psi} \tilde{a}\,.
$$
Remember that Born-Oppenheimer Hamiltonian is given by
\begin{eqnarray*}
\varepsilon^{2+2\delta}\tau_{x}H_{BO}&=&
\varepsilon^{2+2\delta}\tau_{x}\begin{pmatrix}
 e^{i\frac{y}{2\sqrt{\tau_{x}}}}\hat{H}_{+}e^{-i\frac{y}{2\sqrt{\tau_{x}}}}&0\\
0&e^{-i\frac{y}{2\sqrt{\tau_{x}}}}\hat{H}_{-}
e^{i\frac{y}{2\sqrt{\tau_{x}}}}
\end{pmatrix}\,,
\\
&=&\varepsilon^{2\delta}\hat{B}+
\begin{pmatrix}
  2\sqrt{1+\tau_{x}x^{2}}\\0
\end{pmatrix}\,,
\end{eqnarray*}
by using the notations of \eqref{eq.defm1}-\eqref{eq.defm2}.

Let us consider first the linear part after decomposing $\tilde{a}$
into $\tilde{a}=\hat{\chi}\tilde{a}+(1-\hat{\chi})\tilde{a}$, where
the kinetic energy cut-off operator
$\hat{\chi}=\chi(\tau_{x} |\varepsilon D_{q}|^{2})$ has been
introduced in Proposition~\ref{pr.adiabexpl2}:
\begin{multline*}
\hat{U}^{*}H_{Lin}\hat{U}\tilde{a}=\underbrace{\hat{U}^{*}H_{Lin}\hat{U}\hat{\chi}\tilde{a}}_{(I)}
\\
+\underbrace{(\hat{U}^{*}-\hat{U}_{0}^{*})H_{Lin}\hat{U}(1-\hat{\chi})\tilde{a}
+ \hat{U}_{0}^{*}H_{Lin}(\hat{U}-\hat{U}_{0})(1-\hat{\chi})\tilde{a}
}_{(II)}
\\
+ \underbrace{\hat{U}_{0}^{*}H_{Lin}\hat{U}_{0}(1-\hat{\chi})\tilde{a}}_{(III)}\,.
\end{multline*}
The three terms are treated by reconsidering the computations done for
Proposition~\ref{pr.adiabexpl2}. By inserting a cut-off
$\hat{\chi_{1}}$\,, $\chi\prec\chi_{1}$\,, in
$\chi_{1}\hat{U}^{*}H_{Lin}\hat{U}\hat{\chi}\tilde{a}$ with
$(1-\hat{\chi}_{1})\hat{U}^{*}H_{Lin}\hat{U}\hat{\chi}\in
Op\mathcal{N}_{u,g_{\tau}}$, we get
$$
(I)=\varepsilon^{2+2\delta}\tau_{x}H_{BO}\hat{\chi}\tilde{a}
+\mathcal{O}(\varepsilon^{2+4\delta})\quad\text{in}\, L^{2}(\rz^{2})\,.
$$
With $\varepsilon^{-1-2\delta}(\hat{U}-\hat{U_{0}})\in
OpS_{u}(\frac{1}{\langle \sqrt{\tau_{x}}p\rangle^{\infty}},
g_{\tau};\mathcal{M}_{2}(\cz))$ and by using Lemma~\ref{le.ellipt},
 the second term is estimated by
$$
\|(II)\|_{L^{2}}\leq
C\varepsilon^{1+2\delta}\|(1-\hat{\chi})\tilde{a}\|_{L^{2}}
\leq
C'\left[\varepsilon\|\varepsilon^{2\delta}\hat{B}\tilde{a}\|_{L^{2}}+ \varepsilon^{2+4\delta}\right]\,.
$$
We write the third term as
$$
(III)=\varepsilon^{2+2\delta}\tau_{x}H_{BO}(1-\hat{\chi})\tilde{a} +
\varepsilon^{2+2\delta}\tau_{x}\mathcal{D}_{kin}(1-\hat{\chi})\tilde{a}+\mathcal{D}_{pot}(1-\hat{\chi})\tilde{a}\,,
$$
where $\mathcal{D}_{kin}$ and $\mathcal{D}_{pot}$ are defined by
\eqref{eq.Dkin}-\eqref{eq.Dpot}. Following the arguments given in the
proof of Proposition~\ref{pr.adiabexpl2} after these definitions, we
get
\begin{eqnarray*}
\|\mathcal{D}_{pot}(1-\hat{\chi})\tilde{a}\|_{L^{2}}
&\leq &C\varepsilon^{1+2\delta}\|(1-\hat{\chi})\tilde{a}\|_{L^{2}}
\leq
\\
&\leq&
C'\left[\varepsilon\|\varepsilon^{2\delta}\hat{B}\tilde{a}\|_{L^{2}}+ \varepsilon^{2+4\delta}\right]\,,
\\
\|\mathcal{D}_{kin}(1-\hat{\chi})\tilde{a}\|_{L^{2}}
&\leq &
C\left[\varepsilon^{1+4\delta}\|(1-\hat{\chi})\tilde{a}\|_{L^{2}}+
  \varepsilon^{1+2\delta}\|(\varepsilon
  \sqrt{\tau_{x}}|D_{q}|)(1-\hat{\chi})\tilde{a}\|_{L^{2}}\right]
\\
&\leq&
C'\left[\varepsilon\|\varepsilon^{2\delta}\hat{B}\tilde{a}\|_{L^{2}}+ \varepsilon^{2+4\delta}\right]\,.
\end{eqnarray*}
Hence the Euler-Lagrange equation can be written
\begin{eqnarray*}
&&
\varepsilon^{2+2\delta}\tau_{x}H_{BO}\tilde{a}=\lambda_{\psi}\tilde{a}-\varepsilon^{2+2\delta}\tau_{x}G
\hat{U}^{*}\left(|\psi|^{2}\psi\right)+r
\\
\text{with}
&&
\|r\|_{L^{2}}\leq C
\left[
\varepsilon\|\varepsilon^{2\delta}\hat{B}\tilde{a}\|_{L^{2}}
+\varepsilon^{2+4\delta}\right]\,.
\end{eqnarray*}
From Proposition~\ref{pr.existEeps}, with $\varepsilon^{2\delta}\leq
c\tau_{x}$, we know $\|\psi\|_{L^{4}}\leq C$\,. It can be transformed
into $\|\tilde{a}\|_{L^{4}}=\|\hat{U}^{*}\psi\|_{L^{4}}\leq C$ with
$$
\||\psi|^{2}\psi\|_{L^{2}}\leq C\|\psi\|_{L^{6}}^{3}\leq C'\|\tilde{a}\|_{L^{6}}^{3}\,,
$$
by  applying
Lemma~\ref{le.Taylor}, adapted for the metric $g_{\tau}$ like
in the proof Proposition~\ref{pr.L4err}.
We start from the interpolation inequality
$$
\|f\|_{L^{6}}\leq
C\|D_{q}^{2}f\|_{L^{2}}^{\frac{1}{9}}\|f\|_{L^{4}}^{\frac{8}{9}}
\leq \frac{C}{(\tau_{x}\varepsilon)^{\frac{2}{9}}}
\|\tau_{x}^{2}\varepsilon^{2}D_{q}^{2}f\|_{L^{2}}^{\frac{1}{9}}\|f\|_{L^{4}}^{\frac{8}{9}}.
$$
By introducing the operator $\tilde{B}=B(\tau_{x}^{-1/2}\varepsilon q,
\tau_{x}^{\frac{1}{2}}D_{q},\tau,\varepsilon)$ with
\begin{multline*}
B_{\pm}(\tau_{x}^{-1/2} q, \tau_{x}^{1/2}D_{q},\tau,\varepsilon)=(\tau_{x}\varepsilon)^{2}D_{q}^{2}\pm
\frac{\varepsilon}{2}\left(\frac{x}{\sqrt{1+x^{2}}}-\sqrt{\tau_{x}}\right)(\tau_{x}\varepsilon)D_{y}
\\
+\frac{\varepsilon^{2}}{4}\left[\left(\frac{x}{\sqrt{1+x^{2}}}-\sqrt{\tau_{x}}\right)^{2}+\frac{4v}{\ell_{V}^{2}}\right]\,,
\end{multline*}
it becomes
$$
\|f\|_{L^{6}}\leq\frac{C_{\ell_{V}}}{(\tau_{x}\varepsilon)^{\frac{2}{9}}}
\left[\|\tilde{B}f\|_{L^{2}}^{\frac{1}{9}}\|f\|_{L^{4}}^{\frac{8}{9}}
+ \varepsilon^{\frac{2}{9}}\|f\|_{L^{2}}^{\frac{1}{9}}\|f\|_{L^{4}}^{\frac{8}{9}}\right]\,.
$$
The above relation with $f=\tau_{x}^{-\frac{1}{2}}\tilde{a}(\tau_{x}^{-\frac{1}{2}}.)$ which satisfies
$$
\|f\|_{L^{6}}=\tau_{x}^{-\frac{1}{3}}\|\tilde{a}\|_{L^{6}}\quad,\quad
\|f\|_{L^{4}}=\tau_{x}^{-\frac{1}{2}}\|\tilde{a}\|_{L^{4}}
\quad\text{and}\quad
\|\tilde{B}f\|_{L^{2}}=\|\hat{B}\tilde{a}\|_{L^{2}}\,
$$
leads to
\begin{eqnarray*}
\tau_{x}^{-\frac{1}{3}}\|\tilde{a}\|_{L^{6}}
&\leq&
\frac{C_{\ell_{V}}}{(\tau_{x}\varepsilon)^{\frac{2}{9}}}
\left[\|\hat{B}\tilde{a}\|_{L^{2}}^{\frac{1}{9}}\tau_{x}^{-\frac{4}{9}}\|\tilde{a}\|_{L^{4}}^{\frac{8}{9}}
+
\varepsilon^{\frac{2}{9}}\|\tilde{a}\|_{L^{2}}^{\frac{1}{9}}\tau_{x}^{-\frac{4}{9}}\|\tilde{a}\|_{L^{4}}^{\frac{8}{9}}\right]
\\
&\leq &
C_{\ell_{V},G}\tau_{x}^{-\frac{1}{3}}
\left[\varepsilon^{-\frac{2}{9}}\|\hat{B}\tilde{a}\|_{L^{2}}^{\frac{1}{9}}+1\right]\,,
\end{eqnarray*}
and finally
$$
G_{\varepsilon,\tau}\||\psi|^{2}\psi\|_{L^{2}}\leq
G\varepsilon^{2+2\delta}\tau_{x}\|\psi\|_{L^{6}}^{3}
\leq
C_{\ell_{V},G}\left[\varepsilon^{4/3}\tau_{x}\|\varepsilon^{2\delta}\hat{B}\tilde{a}\|_{L^{2}}+
\varepsilon^{2+2\delta}\tau_{x}\right]
$$
With $|\lambda_{\psi}|=\mathcal{O}(\varepsilon^{2+2\delta}\tau_{x})$, we obtain
$$
\|\varepsilon^{2+2\delta}\tau_{x}H_{BO}\tilde{a}\|_{L^{2}}
\leq
C
\left[\varepsilon\|\varepsilon^{2\delta}\hat{B}\tilde{a}\|_{L^{2}}+\varepsilon^{2+4\delta}+\varepsilon^{2+2\delta}\tau_{x}\right]\,,
$$
and recall
$$
\varepsilon^{2+2\delta\tau_{x}}H_{BO}=
\begin{pmatrix}
  \varepsilon^{2\delta}\hat{B}_{+}+2\sqrt{1+\tau_{x}x^{2}}\\
\varepsilon^{2\delta}\hat{B}_{-}
\end{pmatrix}\,.
$$
According to Lemma~\ref{le.ellipt2} below
$$
\|\varepsilon^{2\delta}\hat{B}_{+}\tilde{a}_{+}\|_{L^{2}}+
\|2\sqrt{1+\tau_{x}x^{2}}\tilde{a}_{+}\|_{L^{2}}
\leq
C\|(\varepsilon^{2\delta}\hat{B}_{+}+2\sqrt{1+\tau_{x}x^{2}})\tilde{a}_{+}\|_{L^{2}}\,.
$$
We deduce
$$
\|\varepsilon^{2\delta}\hat{B}\tilde{a}\|_{L^{2}}\leq C\varepsilon^{2+2\delta}\tau_{x}\,.
$$
Plugging this result into the  estimate of the remainder $r$ gives
$$
\varepsilon^{2+2\delta}\tau_{x}H_{BO}\tilde{a}=\lambda_{\tau}\tilde{a}+
G_{\varepsilon,\tau}\hat{U}^{*}(|\psi|^{2}\psi)+ \mathcal{O}(\varepsilon^{3+2\delta}\tau_{x}+\varepsilon^{2+4\delta})\,.
$$
Consider now more carefully the nonlinear term
\begin{multline*}
\hat{U}^{*}\left(|\psi|^{2}\psi\right)
=
\underbrace{
\hat{U}_{0}^{*}\left[|\hat{U}_{0}\tilde{a}|^{2}(\hat{U}_{0}\tilde{a})\right]}_{(1)}
+
\underbrace{
(\hat{U}^{*}-\hat{U}_{0}^{*})\left[|\hat{U}_{0}\tilde{a}|^{2}
(\hat{U}_{0}\tilde{a})\right]}_{(2)}
\\
+
\underbrace{
\hat{U}^{*}\left[|\hat{U}\tilde{a}|^{2}
(\hat{U}\tilde{a})
-
|\hat{U}_{0}\tilde{a}|^{2}
(\hat{U}_{0}\tilde{a})\right]
}_{(3)}\,.
\end{multline*}
By differentiating the relation $\int |f|^{4}=\int
|\hat{U}_{0}f|^{4}$ w.r.t $f$, the first term equals
$$
(1)=|\tilde{a}|^{4}\tilde{a}\,.
$$
By semiclassical calculus in the metric $g_{\tau}$, the operator
$ \sqrt{\tau_{x}}\varepsilon D_{q}\hat{U}_{0}$  equals
$$
 \sqrt{\tau_{x}}\varepsilon D_{q}\hat{U}_{0}
=\hat{U}_{0}\sqrt{\tau_{x}}\varepsilon D_{q} +
\left[\sqrt{\tau_{x}}\varepsilon D_{q}, \hat{U}_{0}\right]\psi
=\hat{U}_{0}\sqrt{\tau_{x}}\varepsilon D_{q}+
\mathcal{O}(\varepsilon)\,,
$$
where the remainder estimate holds in $\mathcal{L}(L^{2}(\rz^{2};\cz^{2}))$\,.
We have already proved $\|\tilde{a}\|_{L^{6}}\leq C$ and
Lemma~\ref{le.Taylor} leads again to
$$
\|\hat{U}\tilde{a}\|_{L^{6}}+
\|\hat{U}_{0}\tilde{a}\|_{L^{6}}\leq
C\,.
$$
With $\varepsilon^{-1-2\delta}(\hat{U}-\hat{U}_{0})\in
OpS_{u}(1,g_{\tau};\mathcal{M}_{2}(\cz))$, this gives
$$
\|(2)\|\leq C\varepsilon^{1+2\delta}\,.
$$
For the third term, we use
\begin{eqnarray*}
\|\left[|\hat{U}\tilde{a}|^{2}
(\hat{U}\tilde{a})
-
|\hat{U}_{0}\tilde{a}|^{2}
(\hat{U}_{0}\tilde{a})\right]\|_{L^{2}}
&&
\leq
C_{0}\|(\hat{U}-\hat{U}_{0})\tilde{a}\|_{L^{6}}
\left[\|\hat{U}\tilde{a}\|_{L^{6}}^{2}+\|\hat{U}_{0}\tilde{a}\|_{L^{6}}^{2}\right]
\\
&&\leq C\varepsilon^{1+2\delta}\,.
\end{eqnarray*}
We have proved
$$
G_{\tau,\varepsilon}\hat{U}^{*}(|\psi|^{2}\psi)-G_{\tau,\varepsilon}|\tilde{a}|^{2}\tilde{a}
=\mathcal{O}(\varepsilon^{3+4\delta}\tau_{x})\,,
$$
which is even better than the estimate for the linear part.\\
For the $\tilde{a}_{-}$ component, we get
$$
H_{BO}
\begin{pmatrix}
  0\\\tilde{a}_{-}
\end{pmatrix}
+G(|\tilde{a}_{-}|^{2}+|\tilde{a}_{+}|^{2})\tilde{a}_{-}
=\frac{\lambda_{\psi}}{\varepsilon^{2+2\delta}\tau_{x}}\tilde{a}_{-}+\mathcal{O}(\varepsilon+
\frac{\varepsilon^{2\delta}}{\tau_{x}})\,,
$$
and it remains to estimate the term
$|\tilde{a}_{+}|^{2}\tilde{a}_{-}$\,.
The first line of the system may be written
$$
\varepsilon^{2\delta}\hat{B}_{+}\tilde{a}_{+}+
2\sqrt{1+\tau_{x}x^{2}}\tilde{a_{+}}
+
G_{\varepsilon,\tau}\left(|\tilde{a}_{+}|^{2}+|\tilde{a}_{-}|^{2}\right)
\tilde{a}_{+}
=\lambda_{\psi}\tilde{a}_{+}+
\mathcal{O}(\varepsilon^{3+2\delta}\tau_{x}+ \varepsilon^{2+2\delta}\tau_{x})\,.
$$
Taking the scalar product with $\tilde{a}_{+}$, with
$\lambda_{\tau}=\mathcal{O}(\varepsilon^{2+2\delta}\tau_{x})$ and
$\hat{B}_{+}\geq 0$, gives
$$
\int |\tilde{a}_{+}|^{4}\leq
C[\|\tilde{a}_{+}\|_{L^{2}}^{2}+(\varepsilon+
\frac{\varepsilon^{2\delta}}{\tau_{x}})\|a_{+}\|_{L^{2}}]
\leq C'\varepsilon^{2+2\delta}\tau_{x}\,.
$$
Then  the same argument as in the estimate of
$\|\tilde{a}\|_{L^{6}}$ gives
\begin{eqnarray*}
\|\tilde{a}_{+}\|_{L^{6}}
&\leq&
C\left[\varepsilon^{-\frac{2}{9}}\|\hat{B}_{+}\tilde{a}_{+}\|_{L^{2}}^{\frac{1}{9}}\|\tilde{a}_{+}\|_{L^{4}}^{\frac{8}{9}}
+
\|\tilde{a}_{+}\|_{L^{2}}^{\frac{1}{9}}\|\tilde{a}_{+}\|_{L^{4}}^{\frac{8}{9}}\right]
\\
&\leq &
C'
\left[(\varepsilon^{-2-2\delta}\|\varepsilon^{2\delta}\hat{B}\tilde{a}\|_{L^{2}})^{\frac{1}{9}}
+
(\varepsilon^{2+2\delta}\tau_{x})^{\frac{1}{18}}
\right](\varepsilon^{2+2\delta}\tau_{x})^{\frac{2}{9}}\\
&\leq&
C''\tau_{x}^{\frac{1}{9}}(\varepsilon^{2+2\delta}\tau_{x})^{\frac{2}{9}}
\leq C'' \varepsilon^{\frac{4+4\delta}{9}}\tau_{x}^{\frac{1}{3}}\,.
\end{eqnarray*}
Again with $\|\tilde{a}\|_{L^{6}}\leq C$, we deduce
$$
\||\tilde{a}_{+}|^{2}\tilde{a}_{-}\|_{L^{2}}\leq
C_{\ell_{V},G}\varepsilon^{\frac{8+8\delta}{3}}\tau_{x}^{\frac{2}{3}}\leq
C_{\ell_{V},G}'\varepsilon\,.
$$
The final result is just a transcription in terms of $a$\,.
\enddemonstration{}
\begin{lemme}
\label{le.ellipt2}
Let $B_{+}$ be the symbol
$$
B_{+}(q,p,\tau_{x},\varepsilon)=\tau_{x}p_{x}^{2}+(\sqrt{\tau_{x}}p_{y}+\frac{\varepsilon}{2}(\frac{\sqrt{\tau_{x}}x}{\sqrt{1+\tau_{x}x^{2}}}-\sqrt{\tau_{x}}))^{2}
+\frac{\varepsilon^{2}v(\sqrt{\tau_{x}}.)}{\ell_{V}^{2}}
$$
introduced in \eqref{eq.defm2}. By setting $\hat{B}_{+}=B_{+}(q,\varepsilon
D_{q},\tau_{x},\varepsilon)$, the operator
$A=\varepsilon^{2\delta}\hat{B}_{+}+2\sqrt{1+\tau_{x}x^{2}}$
is self-adjoint with $D(A)=\left\{u\in L^{2}(\rz^{2}),\, Au\in
  L^{2}(\rz^{2})\right\}$
as soon as $\varepsilon\leq \varepsilon_{0}$, with $\varepsilon_{0}$
independent of $\tau_{x}$\,. Moreover the inequality
$$
\forall u\in D(A),\quad\|\varepsilon^{2\delta}\hat{B}_{+}u\|_{L^{2}}+\|2\sqrt{1+\tau_{x}x^{2}}u\|_{L^{2}}
\leq C \|Au\|_{L^{2}}
$$
holds with a constant $C$ independent of $(\varepsilon,\tau_{x})$\,.
\end{lemme}
\begindemonstration{}
The operator $A$ can be written
$$
A=a_{0}(q,\varepsilon^{1+\delta}D_{q},\tau_{x})+
\varepsilon^{1+\delta}a_{1}(q,\varepsilon^{1+\delta}D_{q},\tau_{x})+\varepsilon^{2+2\delta}a_{2}(q,\tau_{x})\,,
$$
with
$a_{k}\in S_{u}(\langle \sqrt{\tau_{x}}p\rangle^{2-k}+\langle
\sqrt{\tau_{x}}x\rangle^{(1-k)_{+}}, g_{\tau})$ and
$$
a_{0}(q,p,\tau_{x})=p^{2}+2\sqrt{1+\tau_{x}x^{2}}\,.
$$
Therefore the operator $A$ is elliptic in $OpS_{u}(\langle \sqrt{\tau_{x}}p\rangle^{2}+\langle
\sqrt{\tau_{x}}x\rangle, g_{\tau})$ and the result about the domain
follows with
$$
\|\varepsilon^{2\delta}\hat{B}_{+}u\|_{L^{2}}+\|2\sqrt{1+\tau_{x}x^{2}}u\|_{L^{2}}
\leq C (\|Au\|_{L^{2}}+\|u\|_{L^{2}})
$$
for all $u\in D(A)$\,. We conclude with
$$
2\|u\|_{L^{2}}^{2}\leq \langle u\,,\, Au\rangle\leq \|u\|_{L^{2}}\|Au\|_{L^{2}}
$$
due to $\hat{B}_{+}\geq 0$\,.
\enddemonstration{}
\subsection{End of the proof of Theorem \ref{th.mintot}}
\label{se.mintot}
Assume $G,\ell_{V}>0$ be fixed and $\delta\in
(0,\delta_{0}]$\,. Although we dropped $\delta_{0}$ in our notations,
all the constants  in the previous inequalities depend on
$(G,\ell_{V},\delta_{0})$\,.
We assume now  $\tau_{x}\leq \tau_{\ell_{V},G,\delta_{0}}$,
$\varepsilon\leq \varepsilon_{\ell_{V},G,\delta_{0}}$  and
$$
\varepsilon^{2\delta}\leq \tau_{x}^{\frac{5}{3}}\,.
$$
 Proposition~\ref{pr.compen} says
$$
|\mathcal{E}_{\varepsilon,min}-\varepsilon^{2+2\delta}\tau_{x}\mathcal{E}_{\tau,min}|\leq
C\varepsilon^{2+4\delta}\leq C'\varepsilon^{2+2\delta}\tau_{x}^{\frac{5}{3}}
$$
and we recall
$$
|\mathcal{E}_{\tau,min}-\mathcal{E}_{H,min}|\leq C\tau_{x}^{\frac{2}{3}}
$$
by Proposition~\ref{pr.redHR}.\\
When $\psi=\hat{U}
\begin{pmatrix}
  e^{i\frac{y}{2\sqrt{\tau_{x}}}}a_{+}\\
 e^{-i\frac{y}{2\sqrt{\tau_{x}}}}a_{-}
\end{pmatrix}$\,, Proposition~\ref{pr.adEL} says
$$
\|a_{+}\|_{L^{2}}\leq C\varepsilon^{2+2\delta}\tau_{x}
\quad,\quad
\|a\|_{H^{2}}\leq \frac{C}{\tau_{x}}
\quad\text{and}\quad
\|a\|_{L^{4}}+\|a\|_{L^{6}}\leq C\,,
$$
while $a_{-}$ solves the approximate Euler-Lagrange equation
$$
\hat{H}_{-}a_{-}+G\int|a_{-}|^{2}a_{-}=\lambda_{\psi}a_{-}+r_{\varepsilon}\quad
\text{with}\quad \|r_{\varepsilon}\|_{L^{2}}\leq
C_{\ell_{V},G}(\varepsilon+\frac{\varepsilon^{2\delta}}{\tau_{x}})\,.
$$
For $u=\|a_{-}\|_{L^{2}}^{-1}a_{-}$ the energy $\mathcal{E}_{\tau}(u)$ satisfy
$$
|\mathcal{E}_{\tau}(u)-\mathcal{E}_{\tau,min}|
\leq C\tau_{x}^{\frac{2}{3}}
$$
and the above equation becomes
$$
\hat{H}_{-}u+G\int|u|^{2}u=\lambda_{\psi}u+
\mathcal{O}(\varepsilon+\frac{\varepsilon^{2\delta}}{\tau_{x}})\,.
$$
We conclude by referring to Proposition~\ref{pr.appH2} applied to $u$
and then renormalizing for $a_{-}$: For
$\chi=(\chi_{1},\chi_{2})$ with $\chi_{1}\in
\mathcal{C}^{\infty}_{0}(\rz^{2})$, $\chi_{1}^{2}+\chi_{2}^{2}\equiv
1$, $\chi_{1}=1$ in a neighborhood of $0$
\begin{eqnarray*}
  && \|\chi_{2}(\tau_{x}^{\frac{1}{9}}.)a_{-}\|_{L^{2}}\leq
  C\tau_{x}^{\frac{1}{3}}\\
&&
|\mathcal{E}_{\tau}(\chi_{1}(\tau_{x}^{\frac{1}{9}}a_{-})-\mathcal{E}_{H,min})|\leq
C\tau_{x}^{\frac{2}{3}}
\\
&&
d_{\mathcal{H}_{2}}(\chi_{1}(\tau_{x}^{\frac{1}{9}}.)a_{-},
\mathrm{Argmin}~{\mathcal{E}_{H}})\leq
C(\tau_{x}^{\frac{2\nu_{\ell_{V},G}}{3}}+\varepsilon)\,,\quad
\nu_{\ell_{V},G}\in (0,\frac{1}{2}]\,.
\end{eqnarray*}
For the $L^{\infty}$-estimates of $a_{+}$ and $d_{L^{\infty}}(\chi_{1}(\tau_{x}^{\frac{1}{9}}.)a_{-},
\mathrm{Argmin}~{\mathcal{E}_{H}})$, we simply use the interpolation
inequality
$$
\|u\|_{L^{\infty}}\leq C\|u\|_{L^{2}}^{\frac{1}{2}}\|\Delta
u\|_{L^{2}}^{\frac{1}{2}}\leq C\|u\|_{L^{2}}^{\frac{1}{2}}\|u\|_{H^{2}}^{\frac{1}{2}}\,,
$$
valid for any $u\in H^{2}(\rz^{2})$ (write
$u(0)=\int_{\rz^{2}}\hat{u}(\xi)\frac{d\xi}{(2\pi)^{2}}$, cut the
integral according to $|\xi|\lessgtr R$, estimate both term by
Cauchy-Schwartz with $\hat{u}=\frac{1}{|\xi|^{2}}(|\xi|^{2}\hat{u})$
when $|\xi|\geq R$, and then optimize w.r.t $R$).

\section{Additional comments}
\label{se.complements}

We briefly discuss and sketch how our analysis could be adapted to
other problems. No definite statement
is given. Complete proofs require additional work, which may be
done in the future.

\subsection{About the smallness condition of $\varepsilon$ w.r.t $\tau_{x}$}
\label{se.epstau}
In our main result, Theorem~\ref{th.mintot}, condition \eqref{eq.condtaueps} is used, namely
$$
\varepsilon^{2\delta}\leq \frac{\tau_{x}^{\frac{5}{3}}}{C_{\delta}}\,.
$$
One may wonder whether such a condition is necessary in order to
compare the minimization problems for $\mathcal{E}_{\varepsilon}$ and
$\varepsilon^{2+2\delta}\tau_{x}\mathcal{E}_{H}$\,.
When comparing the minimal energies in
Proposition~\ref{pr.compen}, we found
$$
|\mathcal{E}_{\varepsilon,min}-\varepsilon^{2+2\delta}\mathcal{E}_{H,min}|\leq C\varepsilon^{2+4\delta}\,,
$$
while we know that
$\mathcal{E}_{H,min}=\mathcal{E}_{H,min}(\ell_{V},G)$ is a positive
number independent of $\tau_{x}$ and $\varepsilon$\,. Hence it seems
natural to say that $\varepsilon^{2\delta}\ll \tau_{x}$, at least, is
required to ensure that
$\varepsilon^{2+2\delta}\tau_{x}\mathcal{E}_{H,min}$ is a good
approximation of $\mathcal{E}_{\varepsilon,min}$\,. The error is made of three parts:
\begin{itemize}
\item the error term for the Born-Oppenheimer approximation in the
  low-frequency range given in Theorem~\ref{th.BornOp};
\item the error term coming from the truncated high frequency part;
\item the non linear term.
\end{itemize}
The non linear term is
$\frac{G\varepsilon^{2+2\delta}\tau_{x}}{2}\int|\psi|^{4}$, so that a
small error in $\|\psi\|_{L^{4}}^{4}$ will give a negligible term w.r.t
$\varepsilon^{2+2\delta}\tau_{x}\mathcal{E}_{H,min}$\,.
The question is thus mainly about the linear problem. If one looks
more carefully at the error term of Theorem~\ref{th.BornOp}, it is
made of the term
\begin{equation}
  \label{eq.correc}
-\varepsilon^{2}\frac{(\partial_{p_{k}}f_{\varepsilon})(\partial_{p_{\ell}}f_{\varepsilon})}{E_{+}-E_{-}}
\overline{X}_{k}X_{\ell}\,,
\end{equation}
according to Proposition~\ref{pr.hameff}, and of terms coming from the
third order term of Moyal products.
The function $f_{\varepsilon}$ is in our case
$f_{\varepsilon}(p)=\varepsilon^{2\delta}\tau_{x}p^{2}\gamma(\tau_{x}p^{2})$,
$p=(p_{x},p_{y})$, $k,\ell\in \left\{x,y\right\}$\,, and the factors
$X_{x}$ and $X_{y}$ computed in the proof of
Proposition~\ref{pr.adiabexpl1} are at most of order
$\frac{1}{\sqrt{\tau_{x}}}$\,. Hence the quantity \eqref{eq.correc} is
an $\mathcal{O}(\varepsilon^{2+4\delta}\tau_{x})$ which is again
negligible w.r.t
$\varepsilon^{2+2\delta}\tau_{x}\mathcal{E}_{H,min}$\,.
By considering the higher order terms in the Moyal product, the
fast oscillating part of the symbol w.r.t $y$, at the frequency
$\frac{1}{\sqrt{\tau_{x}}}$, deteriorates the estimates:
although there are compensations with the slow variations w.r.t $p_{y}$,
 always multiplied by $\sqrt{\tau_{x}}$, only an $\varepsilon^{k}$ factor without
$\tau_{x}$ appears in the $k$-th order term.\\
Hence, computing the higher order terms,
 at least up to order $3$, in the adiabatic approximation
 and then considering the question
of the high-frequency truncation, is a way to
understand whether the
smallness of $\varepsilon$ w.r.t $\tau_{x}$ is necessary.

\subsection{Anisotropic nonlinearity}
\label{se.aniso}
Our work assumes an isotropic nonlinearity. A more general nonlinear term
 would be
$$
\frac{G\varepsilon^{2+2\delta}\tau_{x}}{2}\int
\alpha_{1}|\psi_{1}|^{4}+2\alpha_{12}|\psi_{1}|^{2}|\psi_{2}|^{2}+
\alpha_{2}|\psi_{2}|^{4}~dxdy\,,\quad
\psi=
\begin{pmatrix}
  \psi_{1}\\\psi_{2}
\end{pmatrix}\,.
$$ Our case is $\alpha_{1}=\alpha_{2}=\alpha_{12}=1$.
 Let
$\psi=\hat{U}\phi$, with
 the unitary transform
$\hat{U}=\hat{U_{0}}+\varepsilon^{1+2\delta}\hat{R}$, (or conversely
$\phi=\hat{U}^{*}\psi$). Then the same arguments as in
Subsection~\ref{se.nonlin} will lead to
\begin{multline*}
\int \alpha_{1}|\psi_{1}|^{4}+2\alpha_{12}|\psi_{1}|^{2}|\psi_{2}|^{2}
+ \alpha_{2}|\psi_{2}|^{4}~dxdy
\\
=
\left (\int \alpha_{1}|\psi_{1}^{0}|^{4}+2\alpha_{12}|\psi_{1}^{0}|^{2}|\psi_{2}^{0}|^{2}
+ \alpha_{2}|\psi_{2}^{0}|^{4}~dxdy\right )
(1+
 \mathcal{O}(\varepsilon^{1+2\delta}))\,,
\end{multline*}
after setting $\psi^{0}(q)=u_{0}(q)\phi(q)$ at every $q=(x,y)$\,.
In our case
$$
u_{0}(x,y)=
\begin{pmatrix}
  C&Se^{i\varphi}\\
Se^{-i\varphi}&-C
\end{pmatrix}
\quad, \quad C=\cos\left(\frac{\theta}{2}\right)
\,,\,
S=\sin\left(\frac{\theta}{2}\right)\,,
$$
with $\theta=\underline{\theta}(\sqrt{\tau_{x}}x)\,,\, \varphi=\frac{y}{\sqrt{\tau_{x}}}$\,.
The point-wise identities
\begin{eqnarray*}
&&
|\psi_{1}^{0}|^{2}+|\psi_{2}^{0}|^{2}
= |\phi_{1}|^{2}+|\phi_{2}|^{2}
\\
&&
|\psi_{1}^{0}|^{2}-|\psi_{2}^{0}|^{2}
=
\cos(\theta)\left(|\phi_{1}|^{2}-|\phi_{2}|^{2}\right)
+2\sin(\theta)\Real(\overline{\phi_{1}e^{i\varphi}}\phi_{2})\,,
\end{eqnarray*}
lead to
\begin{eqnarray*}
&&
\alpha_{1}|\psi_{1}^{0}|^4+2\alpha_{12}|\psi_{1}^{0}|^{2}|\psi_{2}^{0}|^{2}+\alpha_{2}|\psi_{2}^{0}|^{4}
=
\frac{\alpha_{1}+2\alpha_{12}+\alpha_{2}}{4}\left(|\phi_{1}|^{2}+|\phi_{2}|^{2}\right)^{2}
\\
&&\;+ \frac{\alpha_{1}-\alpha_{2}}{2}\left(|\phi_{1}|^{2}+|\phi_{2}|^{2}\right)
\left(\cos(\theta)(|\phi_{1}|^{2}-|\phi_{2}|^{2})+2\sin(\theta)
\Real(\overline{\phi_{1}e^{i\varphi}}
  \phi_{2})\right)
\\
&&\;+
\frac{\alpha_{1}-2\alpha_{12}+\alpha_{2}}{4}
\left(\cos(\theta)(|\phi_{1}|^{2}-|\phi_{2}|^{2})+2\sin(\theta)
\Real(\overline{\phi_{1}e^{i\varphi}}
  \phi_{2})\right)^{2}
\\
&&
\hspace{2cm}
=
\left[\alpha_{1}\cos^{4}(\frac{\theta}{2})+\alpha_{2}\sin^{4}(\frac{\theta}{2})+
\frac{\alpha_{12}}{2}\sin^{2}(\theta)\right]|\phi_{1}|^{4}
\\
&&
\hspace{3cm}
+
\left[\alpha_{1}\sin^{4}(\frac{\theta}{2})+\alpha_{2}\cos^{4}(\frac{\theta}{2})+
\frac{\alpha_{12}}{2}\sin^{2}(\theta)\right]|\phi_{2}|^{4}
\\
&&
\hspace{3cm}
+ \left[\frac{(\alpha_{1}+\alpha_{2})}{2}\sin^{2}(\theta)+
  \alpha_{12}(1+\cos^{2}(\theta))\right]|\phi_{1}|^{2}|\phi_{2}|^{2}
\\
&&
\hspace{1cm}
+\left[(\alpha_{1}-\alpha_{2})+(\alpha_{1}-2\alpha_{12}+
  \alpha_{2})\cos(\theta)\right]\sin(\theta)|\phi_{1}|^{2}
\Real(\overline{\phi_{1}e^{i\varphi}}\phi_{2})
\\
&&
\hspace{1cm} +\left[(\alpha_{1}-\alpha_{2})-(\alpha_{1}-2\alpha_{12}+
  \alpha_{2})\cos(\theta)\right]\sin(\theta)|\phi_{2}|^{2}
\Real(\overline{\phi_{1}e^{i\varphi}}\phi_{2})
\\
&&
\hspace{3cm}
+(\alpha_{1}-2\alpha_{12}+\alpha_{2})\sin^{2}(\theta)\Real\left[(\overline{\phi_{2}e^{i\varphi}}\phi_{2})\right]\,.
\end{eqnarray*}
At least three points have to be adapted from the previous analysis:
\begin{description}
\item[1)] When we take a test function $\psi=\hat{U}
  \begin{pmatrix}
    0\\e^{-i\frac{y}{\sqrt{\tau_{x}}}}a_{-}
  \end{pmatrix}$
the energy $\mathcal{E}_{\varepsilon}(\psi)$ will be close to
$\varepsilon^{2+2\delta}\tau_{x}\mathcal{E}_{\tau}(a_{-})$, with
\begin{multline*}
\mathcal{E}_{\tau}(a_{-})=\langle a_{-}\,,\,\hat{H}_{-}a_{-}\rangle
\\
+\frac{G}{2}
\int
\left[\alpha_{1}\cos^{4}(\frac{\theta}{2})+\alpha_{2}\sin^{4}(\frac{\theta}{2})+
\frac{\alpha_{12}}{2}\sin^{2}(\theta)\right]|a_{-}|^{4}~dxdy\,,
\end{multline*}
and
$\cos(\theta)=\frac{\sqrt{\tau_{x}}x}{\sqrt{1+\tau_{x}x^{2}}}$\,. Hence
before taking the limit $\tau_{x}\to 0$ we have a position dependent
nonlinearity.
This will induce another error term when comparing with the
energy $\mathcal{E}_{H}(a_{-})=\langle a_{-}\,,\,
H_{\ell_{V}} a_{-}\rangle+\frac{G (\alpha_{1}+\alpha_{2}+2\alpha_{12})}{8}\int|a_{-}|^{4}$\,, in the limit
$\tau_{x}\to 0$\,.\\
Another possibility consists in considering the case when
$|\alpha_{2}-\alpha_{1}|+|\alpha_{12}-\alpha_{1}|$ is  small as
$(\varepsilon,\tau_{x})\to (0,0)$\,. The energy
$\mathcal{E}_{\tau}(\psi)$\,, written as,
\begin{multline*}
\langle a_{-}\,,\,\hat{H}_{-}a_{-}\rangle
\\+\frac{G}{2}
\int
\left[\alpha_{1}+(\alpha_{2}-\alpha_{1})\sin^{4}(\frac{\theta}{2})+
\frac{(\alpha_{12}-\alpha_{1})}{2}\sin^{2}(\theta)\right]|a_{-}|^{4}~dxdy\,,
\end{multline*}
will converge to $\mathcal{E}_{H}(a_{-})=\langle a_{-}\,,\,
H_{\ell_{V}} a_{-}\rangle+\frac{G\alpha_{1}}{2}\int|a_{-}|^{4}$\,.
\item[2)] The existence of a minimizer for
  $\mathcal{E}_{\varepsilon}$ and  the variational argument showing that
  $\|a_{+}\|_{L^{2}}^{2}=\mathcal{O}(\varepsilon^{2+4\delta})+\mathcal{E}_{\varepsilon,min}$ when
$\psi=\hat{U} \begin{pmatrix}
    e^{i\frac{y}{\sqrt{\tau_{x}}}}a_{+}\\e^{-i\frac{y}{\sqrt{\tau_{x}}}}a_{-}
  \end{pmatrix}$ is a ground state for $\mathcal{E}_{\varepsilon}$
  will be essentially the same as in the isotropic case.
\item[3)] The analysis and the use of the Euler-Lagrange equation for
  ground states of $\mathcal{E}_{\varepsilon}$, like in
  Subsection~\ref{se.adELeq},
 will certainly be more
  delicate because it will be a system, and the vanishing of the
  crossing terms have to be considered more carefully.
\end{description}

\subsection{Minimization for excited states}
\label{se.minexc}

One may consider like  in \cite{DGJO} the question of minimizing the
energy, for states prepared according to $\psi_{+}$ local eigenvector
of the potential.
Two things have to be modified in order to adapt the previous
analysis:
\begin{description}
\item[1)]  The space of states, on which the energy is
minimized has to be specified. The unitary transform $\hat{U}$
introduced in Theorem~\ref{th.BornOp} provides a simple way to
formulate this minimization problem:  set
$\mathcal{F}_{+}=\Ran\hat{U}P_{+}$ with $P_{+}=
\begin{pmatrix}
  1&0\\
0&0
\end{pmatrix}
$ and consider
$$
\inf_{\psi\in \mathcal{F}_{+}\,,\, \|\psi\|=1}\mathcal{E}_{\varepsilon}(\psi)\,.
$$
\item[2)] In order to get asymptotically as $\varepsilon\to 0$, the
  same scalar minimization problems with $\mathcal{E}_{\tau}$ and
  $\mathcal{E}_{H}$, the external potential $V_{\varepsilon,\tau}$ has
  to be changed. It must be now
  \begin{multline*}
V_{\varepsilon,\tau} (x,y)=
\frac{\varepsilon^{2+2\delta}}{\ell_{V}^{2}}v(\sqrt{\tau_{x}}x,\sqrt{\tau_{x}}y)
\\-
\sqrt{1+\tau_{x}x^{2}}-\varepsilon^{2+2\delta}
\left[
\frac{\tau_{x}^{2}}{(1+\tau_{x}x^{2})^{2}}+
\frac{1}{1+\tau_{x}x^{2}}
\right]\,.
\end{multline*}
\end{description}
The analysis of this problem is essentially the same as for the
complete minimization problem. It is even simpler because the unitary
$\hat{U}$ is directly introduced. A slightly different question is
about the minimization of the energy $\mathcal{E}_{\varepsilon}$
in the space $\mathcal{F}_{+}^{0}=\Ran
\hat{U}_{0}P_{+}$, but the accurate comparison between $\hat{U}$ and
$\hat{U}_{0}$ widely used through this article would lead to similar
results.\\
Possibly  this extension can even be generalized to
higher rank matricial potentials with eigenvectors
$\psi_{1},\ldots, \psi_{N}$, for states modeled on any given
$\psi_{k}$\,.

\subsection{Time dynamics of adiabatically prepared states}
\label{se.timedyn}
Nonlinear adiabatic time evolution has been considered recently in
\cite{CaFe}.
Note that our problem is slightly different because, we are
considering a spatial adiabatic problem, but some techniques may be
related.\\
When $\psi_{0}=\hat{U}
\begin{pmatrix}
  0\\a_{-,0}
\end{pmatrix}
$
(resp. $\psi_{0}=\hat{U}
\begin{pmatrix}
  a_{+,0}\\0
\end{pmatrix}
$ ), the question is whether the solution $\psi(t)$ to
$$
i\partial_{t}\psi=H_{Lin}\psi+ G_{\varepsilon,\tau}|\psi|^{2}\psi\quad,\quad \psi(t=0)=\psi_{0}
$$
remains close to $\hat{U}
\begin{pmatrix}
  0\\ a_{-}(t)
\end{pmatrix}
$
 (resp. $
 \begin{pmatrix}
   a_{+}(t)\\0
 \end{pmatrix}
$)
with
\begin{eqnarray*}
  &&i\partial_{t}a_{-}=\varepsilon^{2+2\delta}\tau_{x}\hat{H}_{-}a_{-}+\frac{G\varepsilon^{2+2\delta}\tau_{x}}{2}|a_{-}|^{2}a_{-}\quad,\quad
  a_{-}(t=0)=a_{-,0}\,,\\
\text{resp.}
&&
i\partial_{t}a_{+}=\varepsilon^{2+2\delta}\tau_{x}\hat{H}_{+}a_{+}+\frac{G\varepsilon^{2+2\delta}\tau_{x}}{2}|a_{+}|^{2}a_{+}\quad,\quad
  a_{+}(t=0)=a_{+,0}\,.
\end{eqnarray*}
More precisely, the question is about the range of time where this
approximation is valid:  what is the size of $T_{\varepsilon}$ w.r.t
$\varepsilon$ such that $\sup_{t\in [-T_{\varepsilon},T_{\varepsilon}]}\|\psi(t)-\hat{U}
\begin{pmatrix}
  0\\a_{-}(t)
\end{pmatrix}
\|$ (resp. $\sup_{t\in [-T_{\varepsilon},T_{\varepsilon}]}\|\psi(t)-
\begin{pmatrix}
  a_{+}(t)\\0
\end{pmatrix}
\|$) remains small.\\
Since the approximation of $\hat{U}^{*}H_{Lin}\hat{U}$ by
$\begin{pmatrix}
  \hat{H}_{+}&0\\
0 &\hat{H_{-}}
\end{pmatrix}$ is good in the low frequency range, a natural
assumption will be that the initial data are supported in the low
frequency region
$$
a_{0,\pm}=\chi(\tau_{x}|\varepsilon D_{q}|^{2})a_{0,\pm}\,,
$$
for some compactly supported $\chi$\,. Then the question is whether the
norm of $(1-\tilde{\chi}(\varepsilon^{2}\tau_{x}D_{q}^{2}))\psi(t)$
remains small for $t\in [0,T_{\varepsilon}]$ for $\tilde{\chi}\in
\mathcal{C}^{\infty}_{0}(\rz)$, $\chi\prec \tilde{\chi}$\,.
 Then the two last parts of Section~\ref{se.minipb}, concerned with
  the high frequency part and the effect of $\hat{U}$ on the
  nonlinear term, have to be reconsidered.
\\
Note that the adiabatically prepared state with $\psi_{0}=\hat{U}
\begin{pmatrix}
  a_{0,+}\\0
\end{pmatrix}
$
are probably not stable for
very long time,
$T_{\varepsilon}=\mathcal{O}(e^{\frac{c}{\varepsilon}})$, because the
characteristic set
$$
\mathcal{C}_{\lambda}=
\left\{(q,p)\in \rz^{4}, \quad
\det\left(\varepsilon^{2\delta}\tau_{x}p^{2}+V_{\varepsilon,\tau}(q)+M(q)-\lambda\right)
=0\right\}
$$
contains two components when $\lambda\geq \min E_{+}$, one
corresponding to the higher level of $M(q)$ with $|p|^{2}\leq
C(\lambda)$,  and another one for the
lower level of $M(q)$ but with large $p$'s. This means that a tunnel
effect will occur between the two levels, so that  adiabatically
prepared states, with energies close to $\lambda\geq \min E_{+}$, will
not remain in this state for (very) large times.


\appendix

\section{Semiclassical calculus}
\label{se.semiclass}
\subsection{Short review in the scalar case}
\label{se.shscal}
Consider a H{\"o}rmander metric, $g$, that is a metric on
$\rz^{2d}_{q,p}$, which satisfies the uncertainty principle, the
slowness and temperance conditions (see \cite{Hor,BoLe})  and
consider $g$-weights (slow and tempered for $g$) $M$, $M_{1}$, $M_{2}$\,.
The symplectic form on $\rz^{2d}_{q,p}$ is denoted by $\sigma$~:
$$
\sigma( X, X')=\sum_{j=1}^{d}p^{j}q'_{j}-q_{j}{p'}^{j}\,,\quad
X=(q,p)\,, X'=(q',p')\,.
$$
The dual metric $g^{\sigma}$ is given by
$g_{X}^{\sigma}(T)=\sup_{T'\neq 0}\frac{|\sigma(T,T')|^{2}}{g_{X}(T')}$ and the
gain associated with $g$ is
\begin{equation}
  \label{eq.uncertainty}
\lambda(X)=\inf_{T\neq 0}\left(\frac{g^{\sigma}_{X}(T)}{g_{X}(T)}\right)^{1/2} \geq 1 \quad(\text{uncertainty})\,.
\end{equation}
In the scalar case,  the space $S(M,g)$ is then the subspace of ${\mathcal
  C}^{\infty}(\rz^{2d}_{p,q};\cz)$ of functions  such that
$$
\forall N\in\nz,\exists
C_{N}\,\quad\sup_{g_{X}(T_{\ell})=1}|(T_{1}\ldots
T_{N}a)(X)|\leq C_{N}M(X)\,,
$$
after identifying a vector field $T_{\ell}$ with a first-order
differentiation operator. For a given $N\in\nz$,  the system of seminorms
$(p_{N,M,g})_{N\in\nz}$ defined by
$$
p_{N,M,g}(a)=\sup_{X\in
  \rz^{2d}}\sup_{g_{X}(T_{\ell})=1}M(X)^{-1}|(T_{1},\ldots,T_{N}a)(X)|
$$
 makes $S(M,g)$ a Fr{\'e}chet space.\\
For $a\in {\mathcal S}'(\rz^{2d}_{q,p})$ and $\varepsilon>0$, the Weyl quantized operator
$a^{W}(q,\varepsilon D_{q}):{\mathcal S}(\rz^{d})\to {\mathcal S}'(\rz^{d})$
is given by its kernel
$$
a^{W}(q,\varepsilon
D_{q})(q,q')=\int_{\rz^{d}}e^{i\frac{(q-q').p}{\varepsilon}}a(\frac{q+q'}{2},
p)~\frac{dp}{(2\pi\varepsilon)^{d}}\,.
$$
When no confusion is possible, we shall use the shortest notations
$$
a^{W}(q,\varepsilon D_{q})=a^{W}= a(q,\varepsilon D_{q})\,.
$$
When $a\in S(M,g)$, it sends ${\mathcal S}(\rz^{d})$ (resp. ${\mathcal
  S}'(\rz^{d})$) into itself and the composition $a_{1}\circ a_{2}$
makes sense for $a_{j}\in S(M_{j},g)$\,.
The Moyal product $a_{1}\sharp^{\varepsilon}a_{2}$ is then defined as
the Weyl-symbol of $a_{1}^{W}(q,\varepsilon D_{q})\circ
a_{2}^{W}(q,\varepsilon D_{q})$:
\begin{eqnarray}
\nonumber
 &&\hspace{-0.6cm} a_{1}\sharp^{\varepsilon}a_{2}(X)=\left(e^{\frac{i\varepsilon}{2}\sigma(D_{X_{1}},D_{X_{2}})}
    a_{1}(X_{1})a_{2}(X_{2})\right)\Big|_{X_{1}=X_{2}=X}
\\
\nonumber
&&
\hspace{-0.3cm}=\sum_{j=0}^{J-1}
    \frac{\left(\frac{i\varepsilon}{2}\sigma(D_{X_{1}},D_{X_{2}})\right)^{j}}{j!}a_{1}(X_{1})
a_{2}(X_{2})\Big|_{X_{1}=X_{2}=X}\\
\nonumber
&&+\int_{0}^{1}\frac{(1-\theta)^{J-1}}{(J-1)!}
e^{\frac{i\varepsilon}{2}\theta\sigma(D_{X_{1}},D_{X_{2}})}
\left(\frac{i\varepsilon}{2}\sigma(D_{X_{1}},D_{X_{2}})\right)^{J}a_{1}(X_{1})a_{2}(X_{2})
\Big|_{X_{1}=X_{2}=X}\\
\label{eq.Wexpan3}
&&\hspace{-0.3cm}=\sum_{j=0}^{J-1}
    \frac{\left(\frac{i\varepsilon}{2}\sigma(D_{X_{1}},D_{X_{2}})\right)^{j}}{j!}a_{1}(X_{1})
a_{2}(X_{2})\Big|_{X_{1}=X_{2}=X} + \varepsilon^{J} R_{J}(a_{1},a_{2},\varepsilon)(X)\;,
\end{eqnarray}
where $R_{J}(.,.,\varepsilon)$ is a uniformly continuous bilinear operator from $S(M_{1},g)\times
S(M_{2},g)$ into $S(M_{1}M_{2}\lambda^{-J},g)$ (i.e. any seminorm of $R_{J}(a_{1},a_{2},\varepsilon)$
is uniformly controlled by some bilinear expression of a finite number of
seminorms of $a_{1}$ and $a_{2}$).\\
The three first terms of the previous expansion are given by
\begin{eqnarray}
\nonumber
&&\hspace{-1cm}a_{1}\sharp^{\varepsilon}a_{2}=
 a_{1}a_{2}+
\frac{\varepsilon}{2i}\left[\partial_{p_{k}}a_{1}\partial_{q^{k}}a_{2}-\partial_{q^{k}}a_{1}\partial_{p_{k}}a_{2}\right]
\\
\label{eq.Wexpanexp}
&&\hspace{-1cm}-\frac{\varepsilon^{2}}{8}
\left[(\partial^{2}_{p_{k},p_{\ell}}a_{1})(\partial^{2}_{q^{k}q^{\ell}}a_{2})
+(\partial^{2}_{q^{k},q^{\ell}}a_{1})(\partial^{2}_{p_{k},p_{\ell}}a_{2})-2
(\partial^{2}_{p_{k},q^{\ell}}a_{1})(\partial^{2}_{q^{k},p_{\ell}}a_{2})\right]
\\
\nonumber
&&\hspace{6cm}
 +\varepsilon^{3}R_{3}(a_{1},a_{2},\varepsilon)\,,
\end{eqnarray}
by making use of the Einstein convention
$s^{j}t_{j}=\sum_{j}s^{j}t_{j}$\,.\\

\begin{definition}
With the small parameter $\varepsilon\in (0,\varepsilon_{0})$, it is
more convenient to consider the Fr{\'e}chet-space
$S_{u}(M,g)$ of bounded functions from $(0,\varepsilon_{0})$ to
$S(M,g)$ endowed with the seminorms:
$$
P_{N,M,g}(a)=\sup_{\varepsilon\in (0,\varepsilon_{0})}p_{N,M,g}(a(\varepsilon))\,.
$$
The subscript $_{u}$ stands for \underline{u}niform seminorm
estimates.\\
The space of $\varepsilon$-quantized family of symbols will be denoted
by $OpS_{u}(M,g)$:
$$
\left(a(\varepsilon)\in S_{u}(M,g)\right)\Leftrightarrow
\left(a^{W}(q,\varepsilon D_{q},\varepsilon)\in OpS_{u}(M,g)\right)\,.
$$
\end{definition}
We shall give a variation of this definition in
Subsection~\ref{se.variametric} below for parameter dependent metrics.

We recall the Beals-criterion proved in
\cite{BoCh} (see \cite{NaNi} for the $\varepsilon$-dependent version)
for diagonal  H{\"o}rmander metrics of the form
$g=\sum_{j=1}^{d}\frac{d{q^{j}}^{2}}{\varphi_{j}(X)^{2}}+\frac{d{p_{j}}^{2}}{\psi_{j}(X)^{2}}$\,.\\
Set $\mathfrak{D}=(D_{q^{1}},\ldots, D_{q^{d}}, q^{1},\ldots, q^{d})$ and for
$E\in \nz^{2d}$ introduce the multi-commutator
$\ad_{\mathfrak{D}}^{E}=\prod_{j=1}^{|E|}\ad_{\mathfrak{D}_{j}}^{E_{j}}$
acting on the continuous operators from $\mathcal{S}(\rz^{d})$ to
$\mathcal{S}'(\rz^{d})$ and the weight $M_{E}(X)=(\varphi(X),\psi(X))^{E}$\,.\\
Then the Beals criterion says that an operator
$A:\mathcal{S}(\rz^{d})\to \mathcal{S}'(\rz^{d})$ belongs to
$OpS(1,g)$ if and only if
$$
\ad_{\mathfrak{D}}^{E}A\in \mathcal{L}(L^{2}(\rz^{d}), H^{\varepsilon}(M_{E}))
$$
for all $E\in \nz^{2d}$\,, when the Sobolev space
$H^{\varepsilon}(M_{E})$ is given by
$$
\|u\|_{H^{\varepsilon}(M_{E})}=\|M_{E}(q,\varepsilon D_{q})u\|_{L^{2}}\,.
$$
Moreover the family of seminorms
$(Q_{N})_{N\in\nz}$, defined by,
\begin{equation}
  \label{eq.semnBeals}
Q_{N}(A)=\sup_{\varepsilon\in (0,\varepsilon_{0})}\max_{|E|=N}\varepsilon^{-N}\|\ad_{\mathfrak{D}}^{E}A(\varepsilon)\|_{\mathcal{L}(L^{2}(\rz^{d}),  H^{\varepsilon}(M_{E},g))}
\end{equation}
on $OpS_{u}(M,g)$, is uniformly equivalent to the family
$(P_{N}(a))_{N\in\nz}$ on $S_{u}(M,g)$ after the
identification $A(\varepsilon)=a^{W}(q,\varepsilon D_{q},\varepsilon)$\,.
``Uniformly'' means here that the comparison of the two
topologies is expressed with constants independent of $\varepsilon\in
(0,\varepsilon_{0})$\,.\\
Below is an example of a metric which satisfies all the assumptions and
which is used in our computations
\begin{equation}
  \label{eq.exemple}
g=\frac{d{q'}^{2}}{(1+|q'|^{2})^{\varrho}}+d{q''}^{2}+\frac{dp^{2}}{(1+|p|^{2})^{\varrho'}}\,,
\quad
m\geq 0\,,\quad
\varrho,\varrho'\in [0,1]\,,
\end{equation}
with the gain function $\lambda(q,p)=\langle p\rangle^{\varrho'}$\,.\\
It is convenient to introduce the class of negligible symbols and
operators.
\begin{definition}
\label{de.neg}
An element $a\in S_{u}(1,g)$, belongs to
$\mathcal{N}_{u,g}$ if
$$
\forall N,N'\in\nz,\,\exists C_{N,N'}>0,\quad
\varepsilon^{-N}a(\varepsilon)\in S_{u}(\lambda^{-N'},g)\,.
$$
Similarly, $Op\mathcal{N}_{u,g}$ denotes the
$\varepsilon$-quantized version:
$$
\left(
a^{W}(q,\varepsilon D_{q},\varepsilon)\in Op\mathcal{N}_{u,g}
\right)
\Leftrightarrow
\left(a\in \mathcal{N}_{u,g}\right)
\,.
$$
\end{definition}
Combined with the following relation between cut-off functions it
provides easy estimates from phase-space localization.
\begin{definition}
\label{de.ordercutoff}
  For two cut-off functions $\chi_{1},
\chi_{2}\in \mathcal{C}^{\infty}_{0}(\rz^{2d})$, $0\leq \chi_{1,2}\leq
1$, the notation $\chi_{1}\prec \chi_{2}$ means that
$\chi_{2}\equiv 1 $ in a neighborhood of $\supp \chi_{1}$\,.
\end{definition}
For example the pseudodifferential calculus leads to
\begin{equation}
  \label{eq.truncneg}
\left(\chi_{1}\prec \chi_{2}\right)\Rightarrow
\left(\forall a\in S_{u}(M,g),\quad
  (1-\chi_{2})\sharp^{\varepsilon}a\sharp^{\varepsilon}\chi_{1}\in \mathcal{N}_{u,g}\right)\,,
\end{equation}
when the weight $M$ satisfies $M\leq C\lambda^{C}$ for some $C>0$\,.
\subsection{Applications to matricial operators}
\label{se.matsemi}
Operator valued pseudodifferential calculus has been studied in
\cite{Bak}. When $\mathfrak{h}$ is a Hilbert space it suffices to
 tensorize the previous calculus with
$\mathcal{L}(\mathfrak{h})$, which corresponds to the
 componentwise definition in $\mathcal{M}_{n}(\cz)$ when
 $\mathfrak{h}=\cz^{n}$\,.
The corresponding class of symbols
 associated with a H{\"o}rmander metric $g$ and a $g$-tempered weight $M$
 is denoted by  $S(M,g;\mathcal{L}(\mathfrak{h}))$
while the set of bounded families in
$S(M,g;\mathcal{L}(\mathfrak{h}))$
 parametrized by $\varepsilon\in (0,\varepsilon_{0})$ is denoted by
$S_{u}(M,g;\mathcal{L}(\mathfrak{h}))$\,.\\
The asymptotic expansions
\eqref{eq.Wexpan3}-\eqref{eq.Wexpanexp} of the Moyal
product clearly holds (see \cite{PST} for a presentation without
specifying the remainder terms) if one takes care of the order of the symbols. For
example, the two first terms of the expansion of
 a commutator $[a_{1}(q,\varepsilon D_{q}),a_{2}(q,\varepsilon
D_{q})]$ are
\begin{multline}
  \label{eq.commmat}
a_{1}\sharp^{\varepsilon} a_{2}-a_{2}\sharp^{\varepsilon}a_{1}= \left[a_{1}(q,p),
  a_{2}(q,p)\right]
\\
+
\frac{\varepsilon}{2i}(\partial_{p}a_{1}\partial_{q}a_{2}-\partial_{q}a_{1}\partial_{p}a_{2})
-
\frac{\varepsilon}{2i}(\partial_{p}a_{2}\partial_{q}a_{1}-\partial_{q}a_{2}\partial_{p}a_{1})
\\+
\varepsilon^{2}\left[R_{2}(a_{1},a_{2},\varepsilon)-R_{2}(a_{2},a_{1},\varepsilon)\right]
\end{multline}
with no simpler expression when the matricial symbols do not
commute.\\
When the H{\"o}rmander metric $g$ has the form
$\sum_{j}\frac{d{q^{j}}^{2}}{\varphi_{j}^{2}}+\frac{dp_{j}^{2}}{\psi_{j}^{2}}$ the uniform Beals
criterion also holds: In the seminorms $Q(a)$ defined in
\eqref{eq.semnBeals} simply replace $L^{2}(\rz^{d})$ by the Hilbert tensor
product $L^{2}(\rz^{d})\otimes \mathfrak{h}$ and consider the operators
$\mathfrak{D}$ and $M_{E}(q,\varepsilon D_{q})$ as the diagonal ones
 $\mathfrak{D}\otimes
\Id_{\mathfrak{h}}$ and  $M_{E}(,\varepsilon D_{q})\otimes
\Id_{\mathfrak{h}}$\,.\\
Finally the Definition~\ref{de.neg} of negligible symbols also makes
sense for matricial symbols after replacing
$OpS(\lambda^{-N'},g)$ by
$OpS(\lambda^{-N'},g;\mathcal{L}(\mathfrak{h}))$\,.\\

We end this section with  standard applications of the Beals
criterion.
\begin{proposition}
\label{pr.beals}
Consider a diagonal H{\"o}rmander metric
$g=\sum_{j=1}^{d}\frac{d{q^{j}}^{2}}{\varphi_{j}(X)^{2}}+
\frac{dp_{j}^{2}}{\psi_{j}(X)^{2}}$
and the constant metric $g_{0}=dq^{2}+dp^{2}$\,.
Assume that
any $f$ chosen in $\left\{\varphi_{j},\psi_{j},\, j\in
  \left\{1,\ldots,d\right\}\right\}$
is a $g_{0}$-slow and
-tempered weight such that
$f(q,\varepsilon D_{q})^{s}$ belongs
$OpS_{u}(f^{s},g_{0})$ for any $s\in\nz$\,.~\footnote{This
 last condition is redundant after a possible modification of
 $\varphi_{j}$ and $\psi_{j}$ if one refers to \cite{BoCh}, but easier
 to check directly in our examples than giving the general proof.}\\
Assume that $A\in OpS_{u}(1,g;\mathcal{L}(\mathfrak{h}))$ is
a family of invertible operators in $\mathcal{L}(L^{2}(\rz^{d})\otimes
\mathfrak{h})$ and such that
$\|A^{-1}(\varepsilon)\|$ is uniformly bounded in w.r.t
$\varepsilon\in (0,\varepsilon_{0})$\,. Then $A^{-1}$ belongs to
$OpS_{u}(1,g;\mathcal{L}(\mathfrak{h}))$\,.
\end{proposition}
\begindemonstration{}
We start from the relation
\begin{equation}
  \label{eq.commutB1}
\ad_{\mathfrak{D}}^{E}A^{-1}
=\sum_{|E'|=|E|, E'\in \nz^{|E|}}c_{E,E'}A^{-1}
\prod_{j=1}^{|E|}\left[\left(\ad_{\mathfrak{D}}^{E_{j}'}A\right)A^{-1}\right]\,.
\end{equation}
which holds in $\mathcal{L}(L^{2}(\rz^{d})\otimes \mathfrak{h})$ for
all $E\in \nz^{2d}$\,. Hence the Beals criterion in the metric $g_{0}$
says that $A^{-1}$ belongs to $OpS_{u}(1,g; \mathcal{L}(\mathfrak{h}))$\,.
In particular $A^{-1}$ belongs to
$\mathcal{L}(H^{\varepsilon}(M_{E"}))$ for any $E"\in \nz^{2d}$\,.
This allows to apply the
 Beals criterion in the metric $g$ and yields the result.
\qed
\subsection{Pseudodifferential projections}

An application of the Beals criterion says that a true
pseudodifferential projection can be made from an approximate one
at the principal symbol level. This holds for matricial symbols.
\begin{proposition}
\label{pr.appproj}  Consider a diagonal H{\"o}rmander metric
$g=\sum_{j=1}^{d}\frac{d{q^{j}}^{2}}{\varphi_{j}(X)^{2}}+
\frac{dp_{j}^{2}}{\psi_{j}(X)^{2}}$  with the same properties as in Proposition~\ref{pr.beals}.
Assume that the  operator
$\hat{\Pi}\in OpS_{u}(1,g;\mathcal{L}(\mathfrak{h}))$
satisfies
\begin{eqnarray*}
&& \left(\hat{\Pi}\circ
  \hat{\Pi}-\hat{\Pi}\right)=\varepsilon^{\mu}\hat{R}_{\mu} +
\varepsilon^{\nu}\hat{R}_{\nu}
\,,
\end{eqnarray*}
with $R_{\mu}\in S_{u}(M,g;\mathcal{L}(\mathfrak{h}))$, $R_{\nu}\in
S_{u}(N,g;\mathcal{L}(\mathfrak{h}))$ $\nu >\mu>0$ and $M, N\leq 1$\,.
Assume additionally that there exist $\chi, \chi' \in
S_{u}(1,g;\mathcal{L}(\mathfrak{h}))$, $0\leq \chi\leq \chi'\leq 1$ such
that $\chi\prec \chi'$
and $\chi' R_{\mu}=0$\,.\\
Then for $\varepsilon_{1}\leq \varepsilon_{0}$ small enough, the
operator
$$
\hat{P}=\frac{1}{2i\pi}\int_{|z-1|=1/2}(z-\hat{\Pi})^{-1}~dz
$$
is well defined for $\varepsilon\in (0,\varepsilon_{1})$ and
satisfies
\begin{eqnarray*}
  && \hat{P}\circ \hat{P}= \hat{P}\quad \text{in}\quad
\mathcal{L}(L^{2}(\rz^{d})\otimes \mathfrak{h})\,,\\
\text{and}&&
\left(\hat{P}-\hat{\Pi}\right)
=\varepsilon^{\mu}\hat{\Pi}_{\mu}+\varepsilon^{\nu}\hat{\Pi}_{\nu}\,.
\end{eqnarray*}
with $\Pi_{\mu}\in S_{u}(M,g;\mathcal{L}(\mathfrak{h}))$, $\Pi_{\nu}\in
S_{u}(N,g;\mathcal{L}(\mathfrak{h}))$  and $\hat{\Pi}_{\mu}\circ
\hat{\chi}\in Op\mathcal{N}_{u,g}$\,.
\end{proposition}
\begindemonstration{}
The first result concerned with the definition of $\hat{P}$ is a
direct application of the simple general result in Lemma~\ref{le.TP}
applied with $T=\hat{\Pi}$ and $\mathcal{H}=L^{2}(\rz^{d})\otimes
\mathfrak{h}$\,.
Note that our assumptions $\nu>1$ and $M\leq 1$ ensure that
$\|\hat{\Pi}^{2}-\hat{\Pi}\|\leq \frac{1}{8}$ as soon as
$\varepsilon\leq \varepsilon_{1}$ with $\varepsilon_{1}$ small
enough.\\
Writing
$$
\hat{\Pi}-\hat{P}
=(\hat{\Pi}^{2}-\hat{\Pi})(A_{1}-A_{0})
$$
with
$$
A_{x}=\frac{1}{2i\pi}\int_{|z-x|}(z-\hat{\Pi})^{-1}~dz\,,\quad x\in \{0,1\}\,,
$$
reduces the problem to proving
$$
A_{x}\in
OpS_{u}(1,g;\mathcal{L}(\mathfrak{h}))
\quad\text{or}\quad(z-\hat{\Pi})^{-1}\in
OpS_{u}(1,g;\mathcal{L}(\mathfrak{h}))\,,
$$
when $\left\{|z-x|=1/2\right\}$ with uniform bounds. But this was
proved in Proposition~\ref{pr.beals} as a consequence of the Beals
criterion. This ends the proof.
\enddemonstration{}
\begin{lemme}
  \label{le.TP}
Assume that in an Hilbert space $\mathcal{H}$, the operator
 $T\in \mathcal{L}(\mathcal{H})$ satisfies $\|T^{2}-T\|\leq
\delta<1/4$ and $\|T\|\leq C$. Then there exists
$c_{\delta}<\frac{1}{2}$ such that
\begin{eqnarray*}
&&\sigma(T)\subset \left\{z\in\cz\,,
|z(z-1)|\leq \delta\right\}\subset \left\{z\in \cz\,,|z|\leq
c_{\delta}\right\}\cup
\left\{z\in \cz\,,|z-1|\leq
c_{\delta}\right\}\,,\\
\\
&&
\max\left\{\|(z-T)^{-1}\|, |z-1|=\frac{1}{2}\right\}\leq 2\frac{2C+1}{1-4\delta}\,.
\end{eqnarray*}
Moreover, the operator
$$
P=\frac{1}{2i\pi}\int_{|z-1|=1/2}(z-T)^{-1}~dz\,.
$$
differs from $T$ according to
\begin{eqnarray}
  \label{eq.diffTP}
&&  T-P=(T^{2}-T)(A_{1}-A_{0})=(A_{1}-A_{0})(T^{2}-T)\,,\\
\text{with}
&&
A_{1}=\frac{1}{2i\pi}\int_{|z-1|=\frac{1}{2}}(T-z)^{-1}(1-z)^{-1}~dz\\
\text{and}&&
A_{0}=\frac{1}{2i\pi}\int_{|z|=\frac{1}{2}}(T-z)^{-1}z^{-1}~dz\,.
\end{eqnarray}
\end{lemme}
\begindemonstration{}
 If $z\in \sigma(T)$ then $z(z-1)\in\sigma(T(T-1))\subset \left\{z\in\cz, |z(z-1)|\leq
  \frac{1}{4}\right\}$ (Remember that $|z(z-1)|=\frac{1}{4}$ means
$|Z-\frac{1}{4}|=\frac{1}{4}$
with $Z=(z-\frac{1}{2})^{2}$).
Consider $z\in \cz$ such that $|z-1|=\frac{1}{2}$, then the relation
$$
(T-z)(T-(1-z))=z(1-z)+(T^{2}-T)\,,
$$
with $|z(1-z)|\geq \frac{1}{4}$ and $\|T^{2}-T\|\leq \delta<\frac{1}{4}$, implies
$$
\|(T-z)^{-1}\|\leq \|T-(1-z)\|\|[z(1-z)+T^{2}-T]^{-1}\|\leq
\frac{C+\frac{1}{2}}{\frac{1}{4}-\delta}\quad\text{for}~|z-1|=\frac{1}{2}\,.
$$
The symmetry with respect to $z=\frac{1}{2}$ due to
$(1-T)(1-T)-(1-T)=T^{2}-T$ implies also
$$
\|(T-z)^{-1}\|\leq \frac{C+\frac{1}{2}}{\frac{1}{4}-\delta}\quad\text{for}~|z|=\frac{1}{2}\,.
$$
Compute
\begin{eqnarray*}
  T-P &=& \frac{1}
  {2i\pi}\int_{|z-1|=\frac{1}{2}}\left[T(z-1)^{-1}-(z-T)^{-1}\right]~dz\\
&=& (T-1)P+(T^{2}-T)A_{1}\\
\text{with}\quad
A_{1}&=&\frac{1}{2i\pi}\int_{|z-1|=\frac{1}{2}}(T-z)^{-1}(1-z)^{-1}~dz\,.
\end{eqnarray*}
In particular this implies to
$T(1-P)=(T^{2}-T)A_{1}$
while replacing $T$ with $(1-T)$ and $P$ with $1-P$ leads to
$$
P-T=-T(1-P)+(T^{2}-T)A_{0} \quad\text{with}\quad
A_{0}=\frac{1}{2i\pi}\int_{|z|=\frac{1}{2}}(T-z)^{-1}z^{-1}~dz\,.
$$
Summing the two previous identities yields the result.
\enddemonstration{}

\subsection{Extension to parameter dependent metrics}
\label{se.variametric}
Additionally to the semiclassical (or adiabatic) parameter, we need
other parameters $\tau=(\tau',\tau'')\in (0,1]^{2}$ on which the
metric $g=g_{\tau}$ depends. In general consider $\tau\in \mathcal{T}\subset\rz^{\nu}$
and a family of H{\"o}rmander metrics $(g_{\tau})_{\tau\in
  \mathcal{T}}$ defined on $\rz^{2d}_{q,p}$\,.
\begin{definition}
\label{de.admimetric} The family of metrics
$(g_{\tau})_{\tau\in\mathcal{T}}$ is said admissible if the
uncertainty principle
\eqref{eq.uncertainty} is satisfied and if the slowness
and temperance constants $C_{1},C_{2}, N_{2}$ involved in
\begin{eqnarray*}
  (\text{slowness})&&
\left(g_{\tau,X}(X-Y)\leq \frac{1}{C_{1}}\right)
\Rightarrow
\left( \left(\frac{g_{\tau,X}}{g_{\tau,Y}}\right)^{\pm 1}\leq C_{1}\right)\,,\\
(\text{temperance})
&&
\left(\frac{g_{\tau,X}}{g_{Y}}\right)^{\pm 1}\leq
C_{2}(1+g_{\tau,X}^{\sigma}(X-Y))^{N_{2}}\,,
\end{eqnarray*}
can be chosen uniformly w.r.t $\tau\in \mathcal{T}$\,.\\
Accordingly a family of weights $(M_{\tau})_{\tau\in \mathcal{T}}$
will be admissible if the slowness and temperance constants of
$M_{\tau}$ w.r.t $g_{\tau}$ can be chosen uniform w.r.t $\tau\in \mathcal{T}$\,.
\end{definition}
The important point is that all the estimates of the Weyl-H{\"o}rmander
pseudodifferential calculus (see \cite{Hor,BoLe}), including the equivalence of norms in the
Beals criterion of \cite{BoCh}, occur with constants which are
determined by the dimension $d$, the uncertainty lower bound (which is
$1$ here),  the slowness and temperance constants.
Hence all the pseudodifferential and semiclassical estimates,
(operator norms or seminorms of remainder terms) are uniform w.r.t to
$(\varepsilon,\tau)$ as long as the symbols,
$a_{k}(\varepsilon,\tau)$, $k=1,2$, have
uniformly controlled seminorm in $S(M_{\tau,k},g_{\tau})$,
w.r.t $(\varepsilon,\tau)\in (0,\varepsilon_{0}]\times
\mathcal{T}$\,.\\
The definition of symbol classes $S_{u}(M,g)$ with uniform control of
seminorms w.r.t $\varepsilon\in (0,\varepsilon_{0}]$ can be extended
to admissible families $(M_{\tau},g_{\tau})_{\tau\in \mathcal{T}}$\,.
\begin{definition}
\label{de.Su}
For an admissible family $(M_{\tau},g_{\tau})_{\tau\in \mathcal{T}}$,
the set of parameter dependent symbol $a(X,\varepsilon,\tau)$,
$X\in \rz^{2d}$, $(\varepsilon,\tau)\in (0,\varepsilon_{0}]\times
\tau$ with uniform estimates
$$
\forall N\in\nz,\; \exists C_{N}>0,\;\forall (\varepsilon,\tau)\in
(0,\varepsilon_{0}]
\times \mathcal{T}\,,\quad
p_{N,M_{\tau},g_{\tau}}(a(\varepsilon,\tau))\leq C_{N}\,,
$$
is denoted by $S_{u}(M_{\tau},g_{\tau}, \mathcal{L}(\mathfrak{h}))$, $\tau\in \mathcal{T}$\,.\\
Equivalently the set of semiclassically quantized operators
$a(q,\varepsilon D_{q},\varepsilon,\tau)$ when the symbol $a$ belongs
to
$S_{u}(M_{\tau},g_{\tau};\mathcal{L}(\mathfrak{h}))$ is denoted by
$OpS_{u}(M_{\tau},g_{\tau};\mathcal{L}(\mathfrak{h}))$\,.\\
The set of negligible symbols and operators associated with
$(g_{\tau})_{\tau\in \mathcal{T}}$ with uniform estimates in
Definition~\ref{de.neg} w.r.t $\tau\in \mathcal{T}$ is denoted by
$\mathcal{N}_{u,g_{\tau}}$ and $Op\mathcal{N}_{u,g_{\tau}}$\,.
\end{definition}
\begin{proposition}
  \label{pr.metadm}
For $\tau=(\tau',\tau'')\in (0,1]^{2}$ the family $(g_{\tau})_{\tau\in
  (0,1]^{2}}$ defined on $\rz^{2d}= \rz^{2(d'+d'')}$ by
$$
g_{\tau}=\frac{\tau'dq'^{2}}{\langle \sqrt{\tau'}q'\rangle^{2}}+ \tau''d{q''}^{2} +
\frac{\tau'\tau''dp^{2}}{\langle
  \sqrt{\tau'\tau''}p\rangle^{2}}
$$
is admissible.
\end{proposition}
\begindemonstration{}
It is easier to consider the symplectically equivalent metric
(use the transform $(q',q'',p',p'')\to ({\tau'}^{\frac{1}{2}}q',
{\tau''}^{\frac{1}{2}}q'', {\tau'}^{-\frac{1}{2}}p',
{\tau''}^{-\frac{1}{2}}p'')$)
$$
\tilde{g}_{\tau}=\frac{dq'^{2}}{\langle q'\rangle^{2}}+dq'' +
\frac{{\tau'}^{2}dp'^{2}+{\tau''}^{2}dp''^{2}}{\langle p\rangle_{\tau}^{2}}
$$
after setting $\langle p\rangle_{\tau}^{2}=1+{\tau'}^{2}{p'}^{2}+{\tau''}^{2}{p''}^{2}$\,.
Firstly remember that the metric
$$
g_{(1,1)}=\frac{dq'^{2}}{\langle q'\rangle^{2}}+ dq'' +
\frac{dp^{2}}{\langle
  p\rangle^{2}}\,.
$$
is a H{\"o}rmander metric. The metric $\tilde{g}_{\tau}^{\sigma}$ is given by
$$
\tilde{g}^{\sigma}_{\tau}= \frac{\langle
  p\rangle_{\tau}^{2}}{\tau'^{2}}dq'^{2}+
\frac{\langle p\rangle_{\tau}^{2}}{\tau''^{2}}dq''^{2}
+\langle q'\rangle^{2}dp'^{2}+ dp''^{2}\,.
$$
Hence the uncertainty principle \eqref{eq.uncertainty} is satisfied
with
$$
\lambda_{\tau}(q',q'',p',p'')= \min\left\{\langle
  q'\rangle\frac{\langle p\rangle_{\tau}}{\tau'}, \frac{\langle
    p\rangle_{\tau}}{\tau''}\right\}
\geq \langle p\rangle_{\tau}\geq 1\,.
$$
In order to  check the uniform slowness and temperance of $\tilde{g}_{\tau}$,
introduce the new variables $X_{\tau}=(q',q'',\tau' p',\tau'' p'')$
when $X=(q',q'',p',p'')$ with a similar definition for $Y_{\tau}$ and $T_{\tau}$\,.
\begin{description}
\item[Slowness:] Write
$$
\left(\tilde{g}_{\tau,X}(X-Y)\leq \frac{1}{C_{1}}\right)\Leftrightarrow
\left(\tilde{g}_{(1,1), X_{\tau}}(X_{\tau}-Y_{\tau})\leq \frac{1}{C_{1}}\right)\,.
$$
 When $C_{1}$  is slowness constant  of
  $\tilde{g}_{(1,1)}$, this implies
$$
\left(\frac{\tilde{g}_{(1,1),X_{\tau}}}{\tilde{g}_{(1,1), Y_{\tau}}}\right)^{\pm
  1}\leq C_{1}
$$
which is nothing but
$$
\left(\frac{\tilde{g}_{\tau,X}}{\tilde{g}_{\tau, Y}}\right)^{\pm
  1}\leq C_{1}\,.
$$
\item[Temperance:] Write
$$
\left(\frac{\tilde{g}_{\tau,X}(T)}{\tilde{g}_{\tau, Y}(T)}\right)^{\pm 1}=
\left(\frac{\tilde{g}_{(1,1),X_{\tau}}(T_{\tau})}{\tilde{g}_{(1,1), Y_{\tau}}(T_{\tau})}\right)^{\pm
  1}
\leq C_{2}(1+\tilde{g}_{(1,1),X_{\tau}}^{\sigma}(X_{\tau}-Y_{\tau}))^{N_{2}}\,,
$$
when $C_{2}$ and $N_{2}$ are the temperance constants for
$\tilde{g}_{(1,1)}$\,.
The problem is reduced to showing
$$
\forall X,T\in \rz^{2d}\,,\quad
\tilde{g}_{(1,1),X_{\tau}}^{\sigma}(T_{\tau})\leq \tilde{g}^{\sigma}_{\tau,X}(T)\,.
$$
The expression of $\tilde{g}^{\sigma}_{\tau}$ gives with
$X=(q',q'',p',p'')$ and
$T=(\theta',\theta'',\pi',\pi'')$
$$
\tilde{g}_{(1,1)X_{\tau}}^{\sigma}(T_{\tau})= \langle p\rangle_{\tau}^{2}\theta'^{2}+ \langle
p\rangle^{2}_{\tau}\theta''^{2}+\langle q'\rangle^{2}\tau'^{2}\pi'^{2}
+ \tau''^{2}\pi''^{2}
\leq \tilde{g}_{\tau, X}^{\sigma}(T)
$$
owing to $\max\{\tau',\tau''\}\leq 1$\,.
\end{description}
\enddemonstration{}

\hfill

\noindent {\bf Acknowledgements:}
The authors are very grateful to  Jean Dalibard for explaining the mathematical questions arising
   from a two level atom in a light beam and for many detailed discussions
that took place for the duration of this work. They also wish to thank
G.~Panati, C.~Lebris and S.V.~Ngoc for their occasional help.  This
work was initiated while the second author had a
 CNRS-sabbatical semester  at CMAP in Ecole Polytechnique.
They acknowledge
support from the French ministry Grant ANR-BLAN-0238, VoLQuan.

\bibliographystyle{plain}

\end{document}